\documentclass[11pt,english]{amsart}

\usepackage[T2A,T1]{fontenc}
\usepackage[latin9]{inputenc}
\usepackage{geometry}
\usepackage{xcolor}
\usepackage{babel}
\usepackage{verbatim}
\usepackage{textcomp}
\usepackage{mathrsfs}
\usepackage{amsbsy}
\usepackage{amstext}
\usepackage{amsmath}
\usepackage{amsthm}
\usepackage{amssymb}
\usepackage{setspace}
\usepackage{wasysym}
\usepackage[all]{xy}
\usepackage{tikz-cd}

\usepackage[unicode=true,pdfusetitle,bookmarks=true,bookmarksnumbered=false,bookmarksopen=false,breaklinks=false,pdfborder={0 0 1},backref=false,colorlinks=false]{hyperref}
\usepackage{enumerate}

\makeatletter

\newcommand{\sa}[1]{\overset{#1}{\longrightarrow}}
\newcommand{\ra}{\longrightarrow}
\newcommand{\hra}{\hookrightarrow}
\newcommand{\iso}{\simeq}
\newcommand{\Al}{\mathbb A^1}
\newcommand{\Gm}{\mathbb G_{\mathrm m}}
\newcommand{\cM}{\mathcal M^{\mathrm{gr}}}

\newcommand{\bE}{\mathbb E}
\newcommand{\bC}{\mathbb C}

\newcommand{\bQ}{\mathbb Q}
\newcommand{\bZ}{\mathbb Z}
\newcommand{\bR}{\mathbb R}
\newcommand{\cC}{\mathcal C}
\newcommand{\cO}{\mathcal O}
\newcommand{\uX}{\underline{X}}
\renewcommand{\proof}{\textit{Proof: }}
\newcommand{\disk}{\mathring \Delta}
\newcommand{\dloc}{DLocSys(\mathring \Delta)_\bQ}
\newcommand{\uni}{DLocSys^{uni}(\mathring \Delta)_\bQ}
\newcommand{\cw}{\mathcal{W}}
\newcommand{\tcM}{\tilde{\mathcal M}^{gr}}

\numberwithin{equation}{section}
\numberwithin{figure}{section}

\theoremstyle{plain}
\newtheorem{thm}{Theorem}[section]

\theoremstyle{remark}
\newtheorem{rem}[thm]{Remark}

\theoremstyle{plain}
\newtheorem{lem}[thm]{Lemma}

\theoremstyle{definition}
\newtheorem{defn}[thm]{Definition}

\theoremstyle{definition}
\newtheorem{example}[thm]{Example}
\theoremstyle{plain}
\newtheorem{cor}[thm]{Corollary}
\theoremstyle{plain}
\newtheorem{prop}[thm]{Proposition}

\theoremstyle{plain}
\newtheorem{conj}[thm]{Conjecture}

\RequirePackage{cmap}
\usepackage[all]{xy}

\makeatother

\begin{document}
	
\title[]{A Voevodsky motive associated to a log scheme}
%The motive of a log scheme. Voevodsky's motive assosiated to a log scheme

\author[]{Georgii Shuklin}

\begin{abstract} 
	For each fs log scheme $(X,\mathcal M_X)$ over a field $k$ we construct a geometrical Voevodsky motive $[X]^{log}\in DM_{gm}(k,\bQ)$.
	We prove that, for $k=\bC$, the Betti realization of $[X]^{log}$ is the log Betti cohomology of $(X, \mathcal M_X)$. We
	give applications to  motivic tubular neighborhoods, limit motives and the monodromy filtrations.
	%For each fs log scheme $(X,\mathcal M_X)$ over a field $k$ we construct the geometrical Voevodsky motive $[X]^{log}\in DM_{gm}(k,\bQ)$ and
	%give application to the construction of motivic tubular neighborhood, limit motives and monodromy filtration.
	%We prove that the Betti realization of $[X]^{log}$ is quasi-isomorphic to the complex of singular cochains $Sing^*(X^{log},\bQ)$ on the Kato-Nakayama space $X^{log}.$ We also show that the diagonal $\Delta: (X,\mathcal M_X)\to (X,\mathcal M_X)\times (X,\mathcal M_X)$ induces on $[X]^{log}$ the structure of motivic $\mathbb E_\infty$-coalgebra.
\end{abstract}

\maketitle	
		\vspace{-9mm}	
	\setcounter{tocdepth}{1}
	\tableofcontents
	
\begin{sloppypar}	
	\vspace{-8mm}
	\section{Introduction}
	Let $X$ be a scheme over a field $k$. Recall that a \textit{log structure} on $X$ is a sheaf of commutative monoids $\mathcal M$ together with a homomorphism $\alpha:\mathcal M\to \cO_X$ such that $\alpha^{-1}(\cO^*)\iso \cO^*.$ A pair $(X,\mathcal M)$ 
	consisting of a scheme and a log structure 
	is called a \textit{log scheme}. Usually, requires some finiteness conditions on a monoid $\mathcal M$ so we will always assume that all log structures are \textit{fine and saturated} \footnote{In particular, $\mathcal M$ is finitely generated and $\mathcal M/\cO^*$ has no torsion. See \cite{'key-33} for more details.}. Fine and saturated (fs) log schemes form the category $Sch^{Log}_{fs}/k$ which contains $Sch/k$ as a full subcategory. 
	\\It is well known that several Weil cohomology theories can be extended to the category of log schemes. The key result in this area is the paper of Kato and Nakayama \cite{key-4}. The authors define for each fs log scheme $(X, \mathcal M)$ over $\bC$ a certain topological space $X^{log}$ which is now called the \textit{Kato-Nakayama space}. Then the log Betti cohomology of $X$ can be defined as $H^*_{Betti}(X^{log},\bZ)$. Kato and Nakayama also proved several comparison theorems between log étale, log De Rham and log Betti cohomology that generalize the classical comparison results.
	
	One can expect the existence of the motive whom's realizations recover various log cohomology theories. Unfortunately, at the moment, there is no such construction.
	%\\In contrast, there is no construction that allows to define the motive of a log scheme. 
	There are several variants of the construction of the logarithmic motivic category (for example \cite{P16}, \cite{'key-25} or \cite{IKNU19}) but in all cases the connection with the Voevodsky category remains at the level of conjectures.
	% or work only with smooth log schemes.
	\\In this article we define a geometrical Voevodsky motive for each fs log scheme. Namely, let $[-]: Sch/k\ra DM_{gm}(k,\bQ)$\footnote{here $DM_{gm}(k,\bQ)$ is the category of geometrical Voevodsky motives}
	be the functor which maps $X\in Sch/k$ to the homological motive $[X].$ We extend [-] to the $(\infty,1)$-functor
	$$[-]^{log}:Sch^{Log}_{fs}/k\ra DM_{gm}(k,\bQ), \ \ \ \ \ \ \ (X,\mathcal M_X)\mapsto [X]^{log},$$
	which has the following basic properties:  
	%\vspace*{-4mm}
	\vspace{-\topsep}
	\begin{enumerate}
		%\vspace*{-7mm}
		\item 
		Let us endow the category $Sch^{Log}_{fs}/k$ with Cartesian monoidal structure. Then $[-]^{log}$ is monoidal. That is
		for each log scheme $(X,\mathcal M)$ the motive $[X]^{log}$ has the 
		structure of $\mathbb E^\infty$-coalgebra.
		%\vspace*{-2mm}
		\item
		for $k\subset \mathbb C$ the Betti realization of $[X]^{log}$ is quasi-isomorphic to the cohomology of Kato-Nakayama space:
		$$R_{Betti}(\mathrm{[X]^{log}})\simeq \mathrm{Sing^*}(X^{\mathrm{log}})_\mathbb Q.$$
		In particular, the log Betti cohomology $H_{Betti}^*(X^{log},\bQ)$ can be endowed with the canonical mixed Hodge structure.
	\end{enumerate}
	\begin{example}
		Let $i:D\hookrightarrow X$ be a normal crossing divisor in a smooth variety $X$. Then $X$ can be endowed with the canonical log structure $\mathcal M_X:=\mathcal O_X\cap j_*\cO^*_{X-D}.$ Let $\mathcal M_D$ be the restriction of the log structure to $D.$ Let us define \textit{the motivic punctured tubular neighborhood} $[PTN_XD]$ as the homotopy fiber
		$$Cone([D]\oplus[X-D]\to [X])[-1].$$ We prove that
		\begin{equation}\label{iso1}
			[(X,\mathcal M_X)]^{log}\iso [X-D]\,\,\,\,\mathrm{and}\,\,\,\,\,\,  [(D,\mathcal M_D)]^{log}\iso [PTN_XD].
		\end{equation}
	\end{example}
	We expect that our construction will be useful for computing the limit motive functor. That is, let $f:X\ra C$ be a proper semi-stable degeneration over a smooth curve $C/k$. Let $X_s$ be the special fiber. Let us fix a local parameter $t$ around $s\in C.$ %which is defined over $\Q$.Let $X_0$ be the special fiber. Let us fix the local parameter $t$ around $0\in C.$ 
	By the properties of $[-]^{log}$  the map of the motivic punctured tubular neighborhood
	$$[PTN_X(X_s)]\ra [PTN_C(s)]\iso [\Gm]$$
	is a homomorphism of $\mathbb E_\infty$-coalgebras. We define %\textit{the Voevodsky limit motive of $f$} by the formula 
	$$ \mathrm{LM}_f:=[PTN_X(X_s)]^*\otimes^L_{[\Gm]^*}\bQ.$$
	Recall that \textit{the limit motive} is $\Upsilon_{\mathrm{id}_C}([X_\eta]^*)$ where $$\Upsilon_{\mathrm{id}_C}:DM_{et}(k(\eta),\bQ)\ra DM_{et}(k(s),\bQ)$$
	is \textit{the functor of  motivic unipotent nearby cycles}\footnote{here $\eta$ is a generic point of $C$.} (see \cite{'key-18}). We prove the following
	\begin{prop}
		There is a canonical quasi-isomorphism 
		$$\mathrm{ML}_f\iso \Upsilon_{\mathrm{id}_C}([X_\eta]^*). $$
	\end{prop}
	Let us explain the generic picture behind Proposition 1.2 in the case $k=\bC$. Let ${\psi_f:D(Sh(X^{an}_\eta,\bQ))\ra D(Sh(X_s^{an},\bQ))}$ be the functor of nearby cycles and $\Delta$ be a small disk around $s\in C$. Let us choose the punctured tubular neighborhood $PTN_{X^{an}}(X_s^{an})$ such that the map $f:PTN_{X^{an}}(X_s^{an})\to \mathring \Delta$ is a fiber bundle. %Our construction of the Voevodsky limit motive is inspired by the quasi-isomorphism
	Then there is the quasi-isomorphism
	\begin{equation}\label{QI3}
		\mathrm{Sing}^*(PTN_{X^{an}}(X_s^{an}),\bQ)\otimes^L_{\mathrm{Sing}^*(\mathring \Delta,\bQ)}\bQ \iso R\Gamma(X_s,\psi_f\bQ).
	\end{equation}
	One can think about the proposition as the motivic version of \ref{QI3}.
	
	In \cite{St76}, Steenbrink endowed $\psi_f\bQ$ with the structure of a mixed Hodge complex. Let us denote this mixed complex by $\psi_f\bQ^{Hdg}.$
	It follows from the property (2) that the Betti realization of $\mathrm{LM}_f$ is quasi-isomorphic to $R\Gamma(X_0,\psi_f\bQ).$ We propose that this quasi-isomorphism can be lifted on the level of  Hodge complexes.
	\begin{conj}
		There is a quasi-isomorphism of mixed Hodge complexes
		$$R_{Hodge}(\mathrm{LM}_f)\iso R\Gamma(X_s,\psi_f\bQ^{Hdg}).$$
	\end{conj}
Recall that the cohomology of $R\Gamma(X_s,\psi_f\bQ^{Hdg})$ is \textit{the limit Hodge structure} of the family $f$.
	%Using the isomorphism $\mathbb H^n(X_0,\psi_f \bQ)\iso H^n_{Betti}(X_t,\bZ)$, one can endow $ H^n_{Betti}(X_t,\bZ)$ 
Let us denote this Hodge structure by $(H^n(X_t,\bQ),W_\bullet,F_{\mathrm{lim}}^\bullet).$ \footnote{here $X_t$ is a fiber of $f$ at the complex point $t\neq s.$}. 
	
Now, suppose that the monodromy  $M:H^n_{Betti}(X_t,\bQ)\to H^n_{Betti}(X_t,\bQ)$ has only one Jordan block. Then the limit Hodge structure has a Hodge-Tate type. 
	\begin{conj}
		(Kerr, Griffiths, Green) Suppose that the degeneration $f$ is defined over $k\subset \bC$. Then the limit Hodge structure $(H^n(X_t,\bQ),W_\bullet,F_{\mathrm{lim}}^\bullet) $
		can be lifted to a mixed Tate motive over $k$.
	\end{conj}
	Conjecturally above mentioned  motive is a cohomology of $\mathrm{LM}_f$ with respect to motivic t-structure.

	The conjecture is based on the following arithmetic considerations. Let $(V,W_\bullet,F^\bullet) $ be a Hodge-Tate structure. Following Goncharov \cite{G99}, one can define the \textit{period operator} ${P:V_\bC\to V_\bC}$ and, after choosing a basis with respect to $W_\bullet$, the \textit{period matrix}
	\[
	\left(\begin{array}{cccc}
		1 & p_{01} &  p_{02} & ... \\
		0 & 1 &  p_{12} & ... \\
		... & ... & ... & ...\\
		0 & ... & 0 & 1
	\end{array}\right).
	\]
	Let $\bQ(p_{ij})$ be the subfield of $\bC$ generated by the entries of the period matrix. Note that the matrix is defined only up to a multiplication by a non-degenerate rational matrix. So the extension $\bQ(p_{ij})$ is well defined.
	
	%Assume that a Tate-Hodge structure arises from the family $f$ with indecomposable monodromy. The Conjecture 04 implies that $\bQ(p_{ij})$ is a subfield of $\bQ(\frac{log(\alpha)}{2\pi i},\frac{\zeta(3)}{(2\pi i)^3},...,\frac{\zeta(2n-1)}{(2\pi i)^{2n-1}},...)$\footnote{here $\alpha$ runs over all invertible rational number}.
	Now, let $\{Y_t\}$ be the family of Calabi-Yau varieties which is mirror dual to the universal family of quintic 3-folds.
	In \cite{GGK}, Kerr, Griffiths and Green   considered as an example of $f$ the family $\{Y_t\}$ together with the canonical parameter $t=t_{can}$
	\footnote{note that the degeneration $\{Y_t\}$ is defined over $\bQ.$}. Firstly, they deduced from  Conjecture 1.4 that $\bQ(p_{ij})=\bQ(\frac{\zeta(3)}{(2\pi i)^3}).$ Then the authors computed the periods using mirror symmetry. More concretely, they proved that, in a suitable basis, the period matrix has the form
	\[
	\left(\begin{array}{cccc}
		1 & 0 & \frac{25}{12} & \frac{-200\zeta(3)}{(2\pi i)^3} \\
		0 & 1 & \frac{-5}{2} & \frac{25}{12} \\
		0 & 0 & 1 & 0\\
		0 & 0 & 0 & 1
	\end{array}\right).
	\]

	We hope that our construction is related to the Conjecture 1.4. In particular, Conjecture 1.4 follows from Conjecture 1.3 in the case when the motive of all irreducible components of $X_s$ together with all intersections are mixed Tate motives.
	
	Now, let us sketch the construction of the functor $[-]^{log}.$
	Suppose $(X,\mathcal M\sa{\alpha}\cO_X)$ is a fs log scheme over $k=\bC.$ Firstly, let us explain how to compute rational Betti cohomology of $X^{log}$ in terms of the sheaf $\mathcal M.$ Let us denote by $exp(\alpha)$ the complex 
	\[
	...\longrightarrow\underset{(-1)}{0_{\,_{\,}}}\longrightarrow\underset{(0)}{\mathcal{O}_{X, \, an_{\,}}}\longrightarrow\underset{(1)}{\mathcal{M}_{an}^{gr}}\longrightarrow\underset{(2)}{0_{\,_{\,}}}\longrightarrow...
	\]
	Here $gr$ means group completion, $an$ means analytification and the differential is given by the composition $\mathcal{O}_{X,  \, an}\!\overset{\mathrm{exp}}{\rightarrow}\mathcal{O}_{X,  \, an}^{*}\!\to\mathcal{M}_{an}^{gr}$. The following fact was discovered  by Steenbrink (in the case of a log smooth scheme). 
	\begin{prop}
		Let $\pi: X^{log}\ra X^{an}$ be the canonical map. Then there is the quasi-isomorphism 
		$$R\pi_*\bQ\iso \mathrm{colim}_n S^n(exp(\alpha)_\bQ).$$
	\end{prop}
	Here the colimit is taken along the maps $S^{n-1}(exp(\alpha)_\bQ)\to S^n(exp(\alpha)_\bQ)$ induced by the inclusion $\bZ\to exp(\alpha).$  
	
	Now, we define the motivic analog of the complex $exp(\alpha).$ Let us extend the sheaf $\mathcal M^{gr}$ on the big étale site of $X$ by the rule 
	$$\mathcal M^{gr}(Y\sa{g} X):=\Gamma (Y,(f^*_{log}\mathcal M)^{gr})$$
	where $f_{log}^*\mathcal M$ is the pullback of log structure. Let us take $\Al$-invariant replacement and stabilize by $\Omega_{(\Gm,1)}$. As the result we get the object of $DA_{et}(X,\bZ)$ which we will denote by $\mathrm{M^{gr}_\bZ}$ \footnote{here $DA_{et}(X,\bZ)$ is the category of étale motivic sheaves \cite{'key-19}}. This construction turns out to be functorial. So the inclusion $\cO^*\hookrightarrow \cM$ induces the map $\bZ(1)[1]\to \mathrm{M^{gr}_\bZ}$ and we can define the motivic sheaf
	$$\mathrm{colim}_n S^n(\mathrm{M^{gr}}(-1)[-1])\in DA_{et}(X,\bQ).$$
	Here $\mathrm{M^{gr}}:= \mathrm{M^{gr}_\bZ}\otimes_\bZ \bQ$ and $S^n$ are motivic symmetric powers \cite{'key-28}. Let ${f \!:\! X\to \mathrm{Spec(k)}}$ be the  canonical morphism. We define the motive $[X]^{log}$ by the formula 
	$$[X]^{log}:=\mathbb DRf_*(\mathrm{colim}_n S^n(\mathrm{M^{gr}}(-1)[-1])).$$
	Note that the construction makes sense for any based field k.
	
Observe that the motive of a log scheme depends only on the group completion $\cM$. This observation allows to extend the functor $[-]^{log}$ on the category of virtual logarithmic schemes (which was introduced in \cite{'key-25}). A \textit{virtual log structure} on a scheme $X$ is a sheaf of abelian group $\mathcal L$ together with an inclusion $\cO^*\hookrightarrow \mathcal L$. Of course, the most important example is $\mathcal L=\cM.$  Another series of examples arises from virtual log schemes over log point $pt_{log}.$ Namely, let $Y\to pt_{log}$ be a morphism of log schemes. In contrast to  $Sch^{Log}_{fs}/k$, in the category of virtual log schemes there is the morphism $Spec(k)\to pt_{log}$ . So we can define the fiber product $\tilde Y:=Y\times_{pt_{log}}\mathrm{Spec(k)}.$ In particular, in the case of a degeneration $f$, the choosing of parameter gives the morphism $(X_s,\mathcal M_{X_s}) \to pt_{log}.$ Then the motive $[\tilde X_s]^{log}$ is dual to $\mathrm{LM_f}$ (see Proposition \ref{P16.7}).
%One of the benefits of looking at virtual log schemes is that we have, in a sense, the resolution of singularities (see Theorem ?) which is crucial to prove of the isomorphisms \ref{iso1}.

\textbf{Acknowledgments.} I would like to thank my supervisor Prof. Vadim Vologodsky for consistent support and guidance during the running of this project. 
I am also grateful to Artem Prikhodko and Vova Shaidurov for the discussions that helped
me a lot. Finally, I should separately thank Ravil Gabdurakhmanov and Vova Shaidurov for
helping me edit this text. The work was supported by the Russian Science Foundation,
grant 21-11-00153.
\section{Presheaf of virtual log structures}

\subsection*{Definition of virtual log structures}

Let $X\!\in\!\mathrm{Sch/k}$ be a scheme of finite type. Let us define
the category $\mathrm{\mathbf{vLog_{\mathit{X}}}}$ in the following
way: the objects of $\mathrm{\mathbf{vLog_{\mathit{X}}}}$ are pairs
$(\!\mathcal{F},\!i\!:\!\mathcal{O}^{*}\!\hookrightarrow\mathcal{F})$
where $\mathcal{F}$ is an étale sheaf of abelian groups on $X$ and
$i$ is an inclusion, the morphisms of $\mathrm{\mathbf{vLog_{\mathit{X}}}}$
are maps $\mu\!:\!\mathcal{F}_{1}\!\rightarrow\mathcal{F}_{2}$ which
commutes with the inclusions of $\mathcal{O}_{X}^{*}$.

For the pair $(\!\mathcal{F}\!,i)$ let us define the ghost sheaf
as the factor $\overline{\mathcal{F}}\!=\!\mathcal{F}\!/\mathcal{O}^{*}$.
We will call the category $\mathrm{\mathbf{vLog_{\mathit{X}}}}$ \textit{the
	category of all virtual logarithmic structures on} $X$.

\subsection*{Functoriality. }

Let $f\!:\!X\!\rightarrow Y$ be a morphism of schemes. For a virtual
log structure\linebreak{}
$\mathcal{F}_{Y}\!\in\!\mathbf{vLog_{\mathit{Y}}}$ let us define
the pullback of $\mathcal{F}_{Y}$ as a colimit of the diagram 
\[
\xymatrix{f^{-1}\!\mathcal{O}_{Y}^{*}\ar[d]\ar[r]^{f^{-1}(i)} & f^{-1}\!\mathcal{F}_{Y}\ar[d]\\
	\mathcal{O}_{X}^{*}\ar[r]^{\tilde{i}} & f^{*}\!\mathcal{F}_{Y}
}
\]
Notice that $f^{-1}$ is exact. So $f^{-1}\!(i)$ and $\tilde{i}$
are inclusions. Observe that this construction is functorial by $\mathcal{F}_{Y}$.
So we constructed the functor $f^{*}\!\!:\!\mathbf{vLog_{\mathit{Y}}}\longrightarrow\mathbf{vLog_{\mathit{X}}}$.
Moreover, by construction for any $f\!:\!X\!\rightarrow Y$ and $g:Z\longrightarrow X$
we have the natural isomorphism $(fg)^{*}\simeq g^{*}f^{*}$. So we
can define the presheaf of categories (contravariant pseudofunctor)
\[
\mathbf{vLog}:\mathrm{Sch/k\longrightarrow Cat}
\]

\subsection*{Coproduct. }

The category of virtual log structures $\mathbf{vLog_{\mathit{X}}}$
admits a coproduct of any objects. For two virtual log structures
$(\mathcal{F}_{1},i_{1})$ and $(\mathcal{F}_{1},i_{1})$ let us denote
their coproduct by ${(\mathcal{F}_{1}\ast\mathcal{F}_{2},i_{1}\ast i_{2})}$.
It is easy to see that the sheaf $\mathcal{F}_{1}\!\ast\!\mathcal{F}_{2}\in\mathrm{Sh}(X_{et},\bZ)$
can be computed as a pushout: 
\[
\mathcal{F}_{1}\!\ast\!\mathcal{F}_{2}:=\mathcal{F}_{1}\oplus_{\mathcal{O}^{*}}\mathcal{F}_{2}.
\]

\begin{rem}
	\label{'R2.1} Notice that for any morphism $f\!:\!X\!\longrightarrow Y$ the pullback
	functor $f^{*}$ commute with coproduct. Indeed, let us consider the
	diagram
	\[
	\xymatrix@=7pt{ & f^{-1}\mathcal{F}_{2}\ar@{->}[rr]\ar@{-}[d] &  & f^{-1}(\mathcal{F}_{1}\!\ast\!\mathcal{F}_{2})\ar@{->}[dd]\\
		f^{-1}\mathcal{O}^{*}\ar@{->}[rr]\ar@{->}[dd]\ar@{->}[ru] &  & \,\,f^{-1}\mathcal{F}_{1}\,\,\ar@{->}[dd]\ar@{->}[ru]\\
		& f^{*}\mathcal{F}_{2}\ar@{-}[r]\ar@{<-}[u] &  & f^{*}(\mathcal{F}_{1}\!\ast\!\mathcal{F}_{2})\ar@{<-}[l]\\
		\mathcal{O}^{*}\ar@{->}[rr]\ar@{->}[ru] &  & \,\,f^{*}\mathcal{F}_{1}\,\,\ar@{->}[ru]
	}
	\]
	\begin{comment}
	\textbackslash xymatrix@=7pt \{ \& f\textasciicircum\{-1\}\textbackslash mathcal\{F\}\_\{2\}\textbackslash ar@\{->\}{[}rr{]}\textbackslash ar@\{-\}{[}d{]}
	\& \& f\textasciicircum\{-1\}(\textbackslash mathcal\{F\}\_\{1\}\textbackslash !\textbackslash ast\textbackslash !\textbackslash mathcal\{F\}\_\{2\})\textbackslash ar@\{->\}{[}dd{]}\textbackslash\textbackslash{}
	f\textasciicircum\{-1\}\textbackslash mathcal\{O\}\textasciicircum\{{*}\}\textbackslash ar@\{->\}{[}rr{]}\textbackslash ar@\{->\}{[}dd{]}\textbackslash ar@\{->\}{[}ru{]}
	\& \& \textbackslash ,\textbackslash ,f\textasciicircum\{-1\}\textbackslash mathcal\{F\}\_\{1\}\textbackslash ,\textbackslash ,\textbackslash ar@\{->\}{[}dd{]}\textbackslash ar@\{->\}{[}ru{]}\textbackslash\textbackslash{}
	\& f\textasciicircum\{{*}\}\textbackslash mathcal\{F\}\_\{2\}\textbackslash ar@\{-\}{[}r{]}\textbackslash ar@\{<-\}{[}u{]}
	\& \& f\textasciicircum\{{*}\}(\textbackslash mathcal\{F\}\_\{1\}\textbackslash !\textbackslash ast\textbackslash !\textbackslash mathcal\{F\}\_\{2\})\textbackslash ar@\{<-\}{[}l{]}\textbackslash\textbackslash{}
	\textbackslash mathcal\{O\}\textasciicircum\{{*}\}\textbackslash ar@\{->\}{[}rr{]}\textbackslash ar@\{->\}{[}ru{]}
	\& \& \textbackslash ,\textbackslash ,f\textasciicircum\{{*}\}\textbackslash mathcal\{F\}\_\{1\}\textbackslash ,\textbackslash ,\textbackslash ar@\{->\}{[}ru{]}
	\}
	\end{comment}
	Notice that the union of front and right squares is cocartesian square.
	Moreover the front square is a pushout by definition. Using the pushout
	lemma we conclude that the right square is cocartesian. Observe that
	the upper square is cocartesian because $f^{-1}$ is exact. Hence
	the union of left and bottom squares is cocartesian. But the left
	square is a pushout by definition. So using the pushout lemma again
	we get that the bottom square is cocartesian.
\end{rem}

\begin{rem}
	\label{'R2.2} (\textit{Group completion of a sheaf}). Let $\mathcal{X}$ be a site
	and $\mathrm{Sh}(\!\mathcal{X}\!,\mathbf{cMon})$ be a category of
	sheaves of commutative monoids on $\mathcal{X}$. The natural inclusion
	$inc\!:\!\mathrm{Sh}(\!\mathcal{X}\!,\mathbf{Ab})\!\hookrightarrow\mathrm{Sh}(\!\mathcal{X}\!,\mathbf{cMon})$
	admits the left adjoint called group completion $gr\!:\!\mathrm{Sh}(\!\mathcal{X}\!,\mathbf{cMon})\longrightarrow\mathrm{Sh}(\!\mathcal{X}\!,\mathbf{Ab})$.
	This functor can also be described in the following way. Notice that
	the diagram 
	\[
	\xymatrix{\mathrm{Sh}(\!\mathcal{X}\!,\mathbf{Ab})\ar@{^{(}->}[r]^{i}\ar[d]^{inc} & \mathrm{PSh}(\!\mathcal{X}\!,\mathbf{Ab})\ar[d]^{inc}\\
		\mathrm{Sh}(\!\mathcal{X}\!,\mathbf{cMon})\ar@{^{(}->}[r]^{i} & \mathrm{PSh}(\!\mathcal{X}\!,\mathbf{cMon})
	}
	\]
	is commutative. So for any presheaf of monoids $M$ we have $(M^{gr})^{\#}\!\simeq\!(M^{\#})^{gr}$
	where $\#$ is a sheafification. On the other hand, it is easy to
	see that the group completion\linebreak{}
	$gr\!:\!\mathrm{PSh}(\!\mathcal{X}\!,\mathbf{cMon})\!\longrightarrow\mathrm{PSh}(\mathcal{X},\mathbf{Ab})$
	can be computed pointwise: $M^{gr}(U)\!\!=\!\!(M(U))^{gr}$. So for
	any sheaf of monoids $\mathcal{M}$ we have
	\[
	\mathcal{M}^{gr}\!\simeq\!(i(\mathcal{M})^{\#}\!)^{gr}\!\!\simeq\!(i(\mathcal{M})^{gr})^{\#}.
	\]
\end{rem}
\subsection*{The virtualization of a log structure.}
Let $\mathcal{M}\!\overset{\alpha}{\longrightarrow}\mathcal{O}$ be
an integral log structure on $X$ (see \cite{'key-33} for the definition).
By the definition of log structures we have the inclusion $i\!:\!\mathcal{O}^{*}\!\hookrightarrow\mathcal{M}$.
Moreover, any map of log structures $f\!:\!\mathcal{M}_{1}\!\rightarrow\mathcal{M}_{2}$
preserves such inclusions. For the log structure $\mathcal{M}\!\overset{\alpha}{\longrightarrow}\mathcal{O}$
let us define \textit{the virtualization} as a pair $(\mathcal{M}^{gr},i^{gr})$.
By Proposition 1.1.3. of \cite{'key-33} the sheaf of monoids $\mathcal{M}$
is integral iff $\mathcal{M}(U)$ is integral for any $U.$ So the
composition $\mathcal{O}^{*}\!(U)\hookrightarrow\mathcal{M}(U)\hookrightarrow(\!\mathcal{M}(U))^{gr}$
is a monomorphism. Finally, the sheafification is exact. So $i^{gr}$
is an inclusion. 

Observe that we constructed the functor $\mathrm{v}_{X}:\,\mathbf{Log_{\mathit{X}}^{\mathrm{int}}}\rightarrow\mathbf{vLog_{\mathit{X}}}$,
$(\!\mathcal{M}\overset{\alpha}{\!\longrightarrow}\mathcal{O})\longmapsto(\!\mathcal{M}^{gr}\!,i^{gr})$.
We will call this functor \emph{the virtualization functor}. 

\textit{\negmedspace{}\negmedspace{}\negmedspace{}\negmedspace{}\negmedspace{}}\textbf{\textit{Warning.}}
Most of the virtual log structures in we are interested are virtualization
of ordinary log structures. So further we will abuse the notation
and will denote by $\mathcal{M}_{X}^{gr}$ an \textit{arbitrary} virtual
log structure on $X$.

\subsection*{Fine and saturated virtual log structures }

The following lemma shows what natural analogs of fs log schemes exist in virtual log geometry.
\begin{lem}\label{'L2.3}
	Let $X$ be a fine logarithmic scheme. Then the sheaf $\overline{\mathcal{M}}_{X}^{gr}$
	is constructible. Moreover if $X$ is also saturated then $\overline{\mathcal{M}}_{X}^{gr}$
	is torsion-free.
\end{lem}
\textit{\negmedspace{}\negmedspace{}\negmedspace{}\negmedspace{}\negmedspace{}Proof:}
First of all, let us show that $\overline{\mathcal{M}}_{X}^{gr}$
is constructible. The proof  is almost the same as the proof of constructability of the monoidal sheaf $\overline{\mathcal{M}}_{X}$ (see \cite{'key-11}).
%This is a direct consequence of the facts that group completion commute with inverse image and the sheaf $\overline{\mathcal{M}}_{X}$ is a constructible monoid sheaf (the prove of second fact can be found in \cite{'key-11}). Moreover, we can give the direct prove by thesame way as in \cite{'key-11}. 
The problem is local so we
may assume that $X=Spec(k[P])$ together with the canonical log structure.
Let $p_{i}$ be generators of $P$ for $1\leq i\leq r$. For $J\subset\{1,...,r\}$
let $R_{J}:=k[P][p_{i_{k}}^{-1}]/I$ where the set $\{i_{k}\}:=\{1,...,r\}\setminus J$
and the ideal is generated by all $p_{j}$ with $j\in J$. Then each
$X_{J}:=Spec(R_{J})$ is a locally closed subset and $X\simeq\bigcup_{J}X_{J}$.
Moreover, for each geometric point $x_{J}\in X_{J}$ the stack $\overline{\mathcal{M}}_{X,x}$
equals the sharp monoid $(P+\Sigma_{i\notin J}\mathbb{Z}p_{i})/G$,
where $G\subset(P+\Sigma_{i\notin J}\mathbb{Z}p_{i})$ is the subgroup
of invertible elements. Finally, it follows from the Remark \ref{'R2.2} and left adjointness of group completion that $\overline{\mathcal{M}}_{X,x}^{gr}\simeq(\overline{\mathcal{M}}_{X,x})^{gr}$. 

Now, suppose that $X$ is saturated. Let us use two facts:
\begin{enumerate}[(i)]
	\item for the canonical log structure on $Spec(k[P])$ the map
	$\mathcal{M}\longrightarrow\mathcal{O}_{X}$ is an inclusion.
	\item for f.s. log structure $\mathcal{M}\longrightarrow\mathcal{O}_{X}$
	on any $X$ and geometric point $x\in X$ the monoid $\overline{\mathcal{M}}_{x}$
	is fine and saturated.
\end{enumerate}
Again, we may assume that $X=Spec(k[P])$. Observe that a sheaf of
abelian group $F$ is torsion-free if $F_{\overline{x}}$ is torsion-free
for any geometric point $\overline{x}\in X$. %
\begin{comment}
easy, suppose that $a\in F(U)$ is a torsion element. so for any geometric
$x\in X$ we have $a_{x}\in F_{x}$ is a torsion element. then $a_{x}=0$
for any $x$.
\end{comment}
{} Let $a\in(\overline{\mathcal{M}}_{x})^{gr}$ and $a^{n}=1$. Notice
that $1\in\overline{\mathcal{M}}_{x}$. So $a\in\overline{\mathcal{M}}_{x}$
by (ii). Let $b\in\mathcal{M}_{x}$ be a lift of $a$ and $g\in\mathcal{M}(U)$
be a corresponding section. Then $g^{n}=f$ for some invertible function
$f$ on $U$. By (i) we know that $g\in\mathcal{O}(U)$. But $g\cdot(g^{n-1}f^{-1})=1$.
So $g\in\mathcal{O}^{*}(U)$ and $a=1$. $\square$
\begin{defn}
	Let us define \textit{the category of fs virtual log structures}
	$\mathbf{vLog_{\mathit{X}}^{\mathrm{fs}}}$  as the full subcategory
	of $\mathbf{vLog_{\mathit{X}}}$ contained all virtual log structures
	$(\mathcal{F}_{X},i)$ with constructible and
	torsion-free ghost sheaf $\mathcal{\overline{F}}_{X}$ . 
\end{defn}
\begin{rem}
	Observe that $\overline{f^{*}\mathcal{F}}_{X}\simeq f^{-1}\overline{\mathcal{F}}_{X}$.
	So for any morphism $f:Y\longrightarrow X$ the pullback functor $f^{*}$
	maps $\mathbf{vLog_{\mathit{X}}^{\mathrm{fs}}}$ to $\mathbf{vLog_{\mathit{X}}^{\mathrm{fs}}}.$ On the other hand, by construction of coproduct we have $\overline{\mathcal{F}_{1}\ast\mathcal{F}_{2}}\simeq\overline{\mathcal{F}_{1}}\oplus\overline{\mathcal{F}_{2}}$.
	So $\mathbf{vLog_{\mathit{X}}^{\mathrm{fs}}}$
	is closed under coproduct into $\mathbf{vLog_{\mathit{X}}}$. Both
	of these facts allow us to consider the collections of categories
	$\mathbf{vLog_{\mathit{X}}^{\mathrm{fs}}}$
	as the presheaf of monoidal categories (with cocartesian monoidal
	structure) on $\mathrm{Sch/k}$. We will denote this presheaf by
	$\mathbf{vLog^{\mathrm{fs}}}$.
\end{rem}

\section{Category of virtual log schemes}

\subsection*{Grothendieck construction. }

Let $\mathcal{C}:\mathrm{Sch/k}^{\mathrm{op}}\rightarrow\mathbf{Cat}$
be a presheaf of categories. Let us denote by $\int_{\mathrm{Sch/k}}\mathcal{C}$
the corresponding Grothendieck construction %
\begin{comment}
I will use (contravariant) nLab's definition which agree with defintion
from the book ``Peter T. Johnstone - Sketches of an Elephant''
\end{comment}
. Objects of $\int_{\mathrm{Sch/k}}\mathcal{C}$ are pairs $(X,c_{X})$
with $c_{X}\!\in\!\mathcal{C}(X)$, the morphisms of $\int_{\mathrm{Sch/k}}\mathcal{C}$
are pairs $(f:X\!\rightarrow Y,\,f^{\#}:c_{X}\rightarrow f^{*}c_{Y})$.
We will also assume that the following two conditions are satisfied:

\noindent \begin{center}
	\begin{minipage}[t]{0.92\columnwidth}%
		($\ast$) For each $X\in\mathrm{Sch/k}$ the terminal object $\ast_{X}$
		of $\mathcal{C}(X)$ exists and $f^{*}\ast{}_{Y}=\ast_{X}$ for any
		$f:X\rightarrow Y$.
		
		($\times$) For each $X\in\mathrm{Sch/k}$ the product of any two
		objects of $\mathcal{C}(X)$ exists and commutes
		with $f^{*}$ .
	\end{minipage}
	\par\end{center}
\begin{rem}\label{'R3.1}
	Notice that $(\mathbf{vLog})^{\mathrm{op}}$
	and ($\mathbf{vLog^{\mathrm{fs}}})^{\mathrm{op}}$ are presheaves
	of categories satisfying properties ($\ast$) and ($\times$). Indeed,
	the property ($\times$) satisfied by Remark \ref{'R2.1}.
	On the other hand, in both cases we have $\ast_{X}=(\mathcal{O}_{X}^{*},id)$.
\end{rem}
Note \textcolor{black}{that due to the property} ($\ast$) the canonical
projection
\[
\mathcal{U}:\int_{\mathrm{Sch/k}}\mathcal{C}\!\rightarrow\mathrm{Sch/k;\,\,\,\,\,\,}(X,c_{X})\!\longmapsto X,
\]
admits the right adjoint
\[
T:\mathrm{Sch/k}\!\rightarrow\int_{\mathrm{Sch/k}}\mathcal{C};\,\,\,\,\,X\!\longmapsto(X,\ast_{X}).
\]
Indeed, if $(f,f^{\#})\!:\!(X,c_{X})\!\rightarrow T(Y)$ is a morphism
then $f^{\#}\!:\!\!c_{X}\rightarrow\ast_{X}$. So $f^{\#}$ is unique
and $\mathrm{Hom_{\int_{\mathrm{Sch/k}}\mathcal{C}}}(X,T(Y))=\mathrm{Hom_{\mathrm{Sch/k}}}(\mathcal{U}(X),Y)$.
Also notice that $T$ is fully faithful and $\mathcal{U}\!\cdot\!T\!=\!Id_{\mathrm{Sch/k}}$.
\begin{example}
	(\textit{The category of log schemes}) Let $\mathbf{Log}$ be the
	presheaf of all logarithmic structures on $\mathrm{Sch/k}$. Then
	$\mathrm{Sch^{Log}\!/k}\!:=\!\int_{\mathrm{Sch/k}}\mathbf{Log}{}^{\mathrm{op}}$
	is the categories of all log schemes.\linebreak{}
	In this case we have the pair of adjoint functors $T\!:\!\mathrm{Sch/k}\rightleftarrows\mathrm{Sch^{Log}\!/k}:\mathcal{U}$
	where $T$ maps $X$ to trivial log schemes $\mathcal{O}_{X}^{*}\!\overset{\alpha}{\hookrightarrow}\mathcal{O}_{X}$.
	Notice that $f^{*}$ commutes with coproducts in $\mathbf{Log}_{X}$
	because pullbacks and coproducts are the particular
	cases of fiber products in the category $\mathrm{Sch^{Log}/k}$ (see
	\cite{'key-33}). 
\end{example}

\subsection*{Category of fs virtual logarithmic schemes}

Let us define \textit{the category of fs virtual logarithmic schemes}
$\mathrm{Sch_{fs}^{vLog}/k}$ as the Grothendieck construction $\int_{\mathrm{Sch/k}}(\mathbf{vLog^{\mathrm{fs}}}){}^{\mathrm{op}}$.
More explicitly this category can be described as the category of
triples $(X,\mathcal{M}_{X}^{gr},\mathcal{O}_{X}^{*}\overset{i}{\hookrightarrow}\mathcal{M}_{X}^{gr})$
where $\mathcal{\overline{M}}_{X}^{gr}$ is torsion-free and constructible. The morphisms can be defined as pairs $(f,f^{\#})$ where $f^{\#}:f^{*}\mathcal{M}_{Y}^{gr}\rightarrow\mathcal{M}_{X}^{gr}$
commutes with the inclusions of $\mathcal{O}_{X}^{*}$. Again, we have the pair of adjoint functors $T:\mathrm{Sch/k}\rightleftarrows\mathrm{Sch_{fs}^{vLog}/k}:\mathcal{U}$
and $T(X)=(X,\mathcal{O}_{X}^{*},id_{\mathcal{O}_{X}^{*}})$.
\begin{rem}
	(\textit{alternative description of morphisms}) One can also describe morphisms of $\mathrm{Sch_{fs}^{vLog}/k}$ as pairs
	of maps $f:X\longrightarrow Y$ and $f^{\flat}:\mathcal{M}_{Y}^{gr}\longrightarrow f_{*}\mathcal{M}_{X}^{gr}$
	such that the diagram 
	\[
	\xymatrix{\mathcal{O}_{Y}^{*}\ar[d]\ar[r] & \mathcal{M}_{Y}^{gr}\ar[d]^{f^{\flat}}\\
		f_{*}\mathcal{O}_{X}^{*}\ar[r] & f_{*}\mathcal{M}_{X}^{gr}
	}
	\]
	is commutative. Indeed, if we have the pair $(f,f^{\#})$ then we
	can compose $f^{\#}$ with the canonical map $f^{-1}\mathcal{M}_{Y}^{gr}\longrightarrow f^{*}\mathcal{M}_{Y}^{gr}$
	and use the $f^{-1}-f_{*}$ adjunction. On the other hand, suppose
	we have a pair $(f,f^{\flat})$. Then again by $f^{-1}-f_{*}$ adjunction
	there is the commutative diagram 
	
	\[
	\xymatrix{f^{-1}\mathcal{O}_{Y}^{*}\ar[d]\ar[r] & f^{-1}\mathcal{M}_{Y}^{gr}\ar[d]\\
		\mathcal{O}_{X}^{*}\ar[r] & \mathcal{M}_{X}^{gr}
	}
	\]
	So we get the map $f^{*}\mathcal{M}_{Y}^{gr}\longrightarrow\mathcal{M}_{X}^{gr}$
	satisfying the required properties.
\end{rem}

\subsection*{Cartesian product. }

Let $\mathcal{C}:\mathrm{Sch/k}^{\mathrm{op}}\rightarrow\mathbf{Cat}$
be a presheaf of categories satisfying the property ($\times$) .
Then the product of any two objects $(X,c_{X})$ and $(Y,c_{Y})$
of $\int_{\mathrm{Sch/k}}\mathcal{C}$ exists and can be computed
as a pair $(X\times_{k}Y,\pi_{X}^{*}c_{X}\times\pi_{Y}^{*}c_{Y})$.
Indeed, the pair of maps $(f,f^{\#}):(Z,c_{Z})\rightarrow(X,c_{X})$
and $(g,g^{\#}):(Z,c_{Z})\rightarrow(Y,c_{Y})$ gives rise to the
unique maps $h:Z\longrightarrow X\times_{k}Y$ and $\mu:\,c_{Z}\longrightarrow f^{*}c_{X}\times g^{*}c_{Y}$.
But $f^{*}c_{X}\times g^{*}c_{Y}\simeq((h^{*}\pi_{X}^{*}c_{X})\times(h^{*}\pi_{Y}^{*}c_{Y})$
and we can use the property ($\times$). So we get 
\begin{lem}
	The Cartesian product of any objects exists in the category $\mathrm{Sch_{fs}^{vLog}/k}.$
\end{lem}

\subsection*{The virtualization of a schemes.}

Let $\mathrm{Sch_{int}^{Log}\!/k}$ be the category of integral log
schemes and $\mathrm{Sch^{vLog}\!/k}:=\int_{\mathrm{Sch/k}}\mathbf{vLog}{}^{\mathrm{op}}$
be the category of all virtual log schemes. Then we can define the
(schematic) virtualization functor
\begin{equation}\label{'F3.1}
	v:\mathrm{Sch_{int}^{Log}\!/k\longrightarrow Sch^{vLog}\!/k}
\end{equation}
by the following way: an integral log scheme $(X,\mathcal{M}_{X}\longrightarrow\mathcal{O}_{X})$
maps to $(X,\mathcal{O}_{X}\hookrightarrow\mathcal{M}_{X}^{gr})$
and the morphism $(X\overset{f}{\longrightarrow}Y,\mathcal{M}_{Y}\overset{\alpha}{\longrightarrow}f_{*}\mathcal{M}_{X})$
maps to the pair $f$ and $\alpha^{\prime}:\mathcal{M}_{Y}^{gr}\longrightarrow f_{*}\mathcal{M}_{X}^{gr}$.
Here $\alpha^{\prime}$ can be defined by the following way. Let $i_{X}$
(resp. $i_{Y}$) be a canonical inclusion of category presheaves of
monoids to the category of sheaves of monoids on $X$ (resp. $Y$).
Then the map $\alpha$ gives rise to the map of presheaves of group
$(i_{Y}\mathcal{M}_{Y})^{gr}\longrightarrow f_{*}(i_{X}\mathcal{M}_{X})^{gr}$.
Notice that $\mathcal{M}_{X}^{gr}$ is the sheafification of $(i_{X}\mathcal{M}_{X})^{gr}$.
So we can define $\alpha^{\prime}$ as the sheafification of the map$(i_{Y}\mathcal{M}_{Y})^{gr}\longrightarrow f_{*}\mathcal{M}_{X}{}^{gr}$. 

\subsection*{Relation with fs log schemes.}

The category  $\mathrm{Sch_{fs}^{vLog}/k}$
is closely related with the category of fs log schemes $\mathrm{Sch_{fine}^{Log}\!/k}$
and $\mathrm{Sch_{fs}^{vLog}/k}$. Namely, thanks to Lemma \ref{'L2.3}, the restriction of (\ref{'F3.1})
on the category $\mathrm{Sch_{fs}^{vLog}/k}$
gives the virtualization functor 
\begin{equation}\label{'F3.2}
	v:\mathrm{Sch_{fs}^{Log}\!/k\longrightarrow Sch_{fs}^{vLog}\!/k}\,\,\,\,\,\,v:\mathrm{Sch_{fine}^{Log}\!/k\longrightarrow Sch_{fine}^{vLog}\!/k}
\end{equation}
that maps an fs log schemes $(X,\mathcal{M}\overset{\alpha}{\rightarrow}\mathcal{O}_{X})$
to the fs virtual log scheme $(X,\mathcal{O}_{X}^{*}\hookrightarrow\mathcal{M}^{gr})$.
\begin{rem}
	We will not use this further, but it should be noted that all the constructed schematic virtualization functors commute with products.
	\begin{comment}
	Firstly, let us explain this for the functor $\mathrm{v}:\mathrm{Sch_{int}^{Log}\!/k\longrightarrow Sch^{vLog}\!/k}$.
	It can be verified that the collection of functors $\mathrm{v}_{X}:\,\mathbf{Log_{\mathit{X}}^{\mathrm{int}}}\rightarrow\mathbf{vLog_{\mathit{X}}}$
	give rise to the natural transformation 
	\[
	\mathrm{v}:\,\mathbf{Log^{\mathrm{int}}}\rightarrow\mathbf{vLog}
	\]
	(in other words, the virtualization commutes with pullbacks). So we
	have the functor 
	\[
	\int\mathrm{v}:\mathrm{Sch_{int}^{Log}\!/k\longrightarrow Sch^{vLog}\!/k}.
	\]
	Moreover, one can prove that $\mathrm{v}_{X}$ commutes with coproducts
	for any $X$. Hence we conclude that $\int\mathrm{v}$ commute with
	product. It only remains to check that $\int\mathrm{v}$ and $v$
	are coincide. Finally, using charts, one can check that any product
	of fs ( resp. fine) log schemes is again fs (resp. fine) log schemes
	(the similar statement for virtual fs and fine log schemes are follows
	from \textcolor{red}{Remark 2.3}).
	\end{comment}
\end{rem}
%In \textcolor{red}{Section 7} we will define the motivic cohomology
%as the functors
%\[
%H^{p,q}(-,\mathbb{Q})\!:\!\mathrm{Sch_{fs}^{vLog}\!/k}\longrightarrow\mathrm{Vect}_{\mathbb{Q}}.
%\]
%So the existences of the functor $\mathrm{v}:\mathrm{Sch_{fs}^{Log}\!/k\longrightarrow Sch_{fs}^{vLog}\!/k}$
%allows us to define the motivic cohomology of fs log schemes as a
%cohomology of their virtualization.

\section{Free $\mathbb E_\infty$-algebras}
Suppose we have a monoidal $\infty$-category $\mathcal C^\otimes$ with monoidal unit $I$. Let $\mathcal C_{I/}$ be the under category. One can define on $\mathcal C_{I/}$ the monoidal structure (see formula \ref{F2.1}). Let $I\to A$ be an object of   $\mathcal C_{I/}$. In this section we explain how to construct a free $\mathbb E_\infty$-algebra $\mathrm{Sym}^*(A)$. Throughout this  section we will use the notion of $\infty$-operads which was introduced by Lurie in \cite{key-1}. We will denote by $\mathrm{Op(\infty)}$ the quasi-category of $\infty$-operads.
\\Let $\mathrm{Fin_*}$ be the category of pointed finite sets. Recall that a symmetric monoidal ${\rm \infty-category}$ is a quasi-category $\mathcal C^\otimes$ together with a cocartesian fibration $q:\mathcal C^\otimes \rightarrow \mathrm{N(Fin_*)}$ such that for any n we have the natural equivalence 
$$ \mathcal C^\otimes_{[n]}\simeq \mathcal (C^\otimes_{[1]})^n $$
induced by the set of maps $\{ \rho^i : [n] \rightarrow [1] \}_{1\leq i\leq n} $.
For any $\infty$-operad $\mathcal O$ let us denote by $\mathrm{Alg}\mathcal{_O (C)}=\mathrm{Map_{Op(\infty)}(\mathcal{O,C^\otimes})}$ the $\infty$-category of $\mathcal O$-algebras. Note that by definition $\mathrm{Alg}\mathcal{_O (C)}$ is a full subcategory of $\mathrm{Fun_{N(Fin_*)}(\mathcal{O,C^\otimes})}.$
\begin{rem}
	Let $\mathcal O=\mathrm{N(Fin_*)}=\mathbb E_\infty$. 
	Let  $\mathcal C:= \mathcal C^\otimes_{[1]}\simeq \mathrm{Fun_{N(Fin_*)}}([1],\mathcal C^\otimes)$ be the underlying $\infty$-category of $\mathcal C^\otimes$. Then the map of simplicial sets $[1]\rightarrow \mathrm{N(Fin_*)}$ induces the forgetful functor 
	$$U:\mathrm{CAlg}\mathcal {(C)} \longrightarrow \mathcal C.$$
	Suppose that infinity direct sums are exist in $\mathcal{C}$ and for any object $c$ the functor $\otimes c$ commutes with colimits. Then by Proposition 3.1.3.13 of \cite{key-1} the forgetful functor admits left-adjoint $S^*$ which can be constructed using coinvariants of symmetric groups $\Sigma_n$. Namely, let $A\in Ob(\mathcal{C}).$ Then 
	$$S^*(A)=\bigsqcup_n (A^{\otimes n})_{h\Sigma_n}$$ together with the multiplication induced by the isomorphism 
	$$(A^{\otimes n})\otimes (A^{\otimes m})\simeq A^{\otimes n+m}. $$
\end{rem}
\subsection*{The category of pointed objects} Note that the inclusion of marked point to [1] induce the functor 
$$\mathcal C^\otimes_{[0]}\simeq \Delta^0 \rightarrow \mathcal C^\otimes_{[1]}.$$ Will will called the associated object the monoidal unit and denote by $I$.
\\Let $\mathcal C_{I/}$ be the under category. In the case then $\mathcal C$ is the nerve of a 1-category one can define the monoidal structure on $\mathcal C_{I/}$ by the rule 
\begin{equation}\label{F2.1}
	(I\rightarrow A)\otimes (I\rightarrow B):= (I\simeq I\otimes I \rightarrow A\otimes B).	
\end{equation}
In general, we can induce the monoidal structure from the category of $\mathbb E_0$-algebras of $\mathcal C$. Following \cite{key-1} let us denote by $\mathrm{Fin_*^{inj}}$ the subcategory of $\mathrm{Fin_*}$ spanned by all objects together with
those morphisms $f : [m] \rightarrow [n]$ such that $f^{-1}([i])$ has at most one element for $1 \leq i \leq n$. The nerve $\mathrm{N(Fin^{inj}_*)}$ is an $\infty$-operad, which we will denote by $\mathbb E_0$.
Note that we have the unique map  $[0]\rightarrow[1]$. The corresponding restriction 
$$\mathrm{Alg_{\mathbb E_0}}(\mathcal C)\longrightarrow \mathrm{Fun_{N(Fin_*)}(\Delta^1,\mathcal C^\otimes)}$$ is an equivalence of categories by Proposition 2.1.3.9 of \cite{key-1}.
%Remark. Atucally statement from HA about $Alg_{/E_0}$ and $Fun_{E_0}$ but both categories are equivalnt to our. Indeed, $Alg_{/E_0}(C):=fib(Alg_{E_0}(C)\rightarrow Alg_{E_0})(E_0)$ but the functor $N(Fin^{inj}_*)\rightarrow N(Fin_*)$ is an inclusion. So   $Alg_{E_0})(E_0)=pt$.
So $\mathrm{Alg_{\mathbb E_0}}(\mathcal C)\simeq \mathcal C_{I/}$. Note that the category of $\infty$-operads can be equipped with natural tensor product $\odot$ (see Section 2.2.5 of \cite{key-1}). The corresponding monoidal structure is closed so we can equipped   $\mathrm{Alg_{\mathbb E_0}}(\mathcal C)$ with the structure of $\infty$-operads. By Proposition 1.8.19 and Example 1.8.20  of \cite{key-2} $\mathrm{Alg_{\mathbb E_0}}(\mathcal C)^\otimes$ is the symmetric monoidal category.  It follows from Lemma \ref{L2.3} that the corresponding tensor product on $\mathrm{Alg_{\mathbb E_0}}(\mathcal C)^\otimes$ coincide with (\ref{F2.1}).
\begin{rem}
	(Forgetful functors). %Following \cite{key-1}
	Let $\mathrm {Triv}$ be the subcategory of $\mathrm{Fin_*}$ whose objects are the objects of $\mathrm{Fin_*}$ ,
	and whose morphisms are the inert morphisms in $\mathrm{Fin_*}$ . Let $\mathrm{Triv^\otimes = N(Triv)}$. Then the inclusion
	$\mathrm{Triv\subseteq Fin_*}$ unduces the functor $\mathrm{Triv^\otimes \rightarrow  N(Fin_*)}$ which exhibits $\mathrm{Triv^\otimes}$  as an $\infty$-operad.
\end{rem}
Let $\Delta^0\rightarrow \mathrm{Triv}$ be the functor corresponding to $[1]\in \mathrm{Ob(Triv)}$. Then by Remark 2.1.3.6 of \cite{key-1} the restriction functor 
$$\mathrm{Alg_{Triv}(\mathcal C^\otimes)\longrightarrow Fun_{N(Fin_*)}([1],C^\otimes)}$$
is a trivial Kan fibration. So we can consider any map of $\infty$-operads as a forgetful functor.
\begin{lem}\label{L2.3}
	The forgetful functor $F:\mathcal{C_{I/}}\rightarrow \mathcal C$ is monoidal.
\end{lem}

\proof : The natural restriction 
$$\mathrm{Map_{Op(\infty)}(\mathbb E_0,\mathcal C^\otimes)\rightarrow Fun_{N(Fin_*)}([1],\mathcal C^\otimes)}$$ 
factors through $\mathrm{Fun_{N(Fin_*)}(\Delta^1,\mathcal C^\otimes)}$.  So it suffice to check that the forgetful functor $\mathrm{Alg_{\mathbb E_0}(\mathcal C)^\otimes \rightarrow \mathcal C^\otimes}$ is monoidal. Now, we can use Proposition 1.8.19 of \cite{key-2}. $\square$
%Put $O=O''=N(Fin_*)$ and $X=[1]$.
\\Let us denote by $\mathcal U$ the forgetful functor $\mathrm{CAlg(\mathcal C_{I/})\rightarrow \mathcal C_{I/}}$.
\begin{lem}\label{L2.4}
	The diagram 
	\begin{center}
		\begin{tikzcd}
			\mathrm{CAlg(}\mathcal C_{I/}) \arrow[r, "\mathrm{CAlg}(F)"] \arrow[d, "\mathcal U"] & \mathrm{CAlg(\mathcal C)} \arrow[d, "U"] \\
			\mathcal C_{I/} \arrow[r, "F"]                                                       & \mathcal C                              
		\end{tikzcd}
	\end{center}
	is commutative. Moreover, the functor $\mathrm{CAlg}(F)$ is an equivalence of categories.
\end{lem}
\proof: By Lemma \ref{L2.3} $F$ induced by the monoidal functor $F^\otimes :\mathcal C^\otimes_{I/}\rightarrow \mathcal C^\otimes$. So the diagram can be rewritten as 
\begin{center}
	\begin{tikzcd}
		{\mathrm{Map_{Op(\infty)}}(\mathbb E_0,\mathcal C_{I/}^\otimes)} \arrow[d, "i_*"] \arrow[rr, "(F^\otimes)^*"] &  & {\mathrm{Map_{Op(\infty)}}(\mathbb E_0,\mathcal C^\otimes)} \arrow[d, "i_*"] \\
		{\mathrm{Map_{Op(\infty)}}(\mathrm{Triv^\otimes},\mathcal C_{I/}^\otimes)} \arrow[rr, "(F^\otimes)^*"]               &  & {\mathrm{Map_{Op(\infty)}}(\mathrm{Triv^\otimes},\mathcal C^\otimes)}               
	\end{tikzcd}
\end{center}
where $i$ is the inclusion $\mathrm{Triv\hookrightarrow N(Fin_*)}$. 
\\Let us prove that $\mathrm{CAlg}(F)$ is an equivalence. Note that the monoidal structure on $\mathrm{Op(\infty)}$ is closed and $\mathrm{Triv}^\otimes$ is the monoidal unit of $\mathrm{Op(\infty)}$. So 
$$\mathrm{CAlg(\mathcal C}_{I/})\simeq \mathrm{CAlg(\mathcal C)}_{I/}$$ and $\mathrm{CAlg}(F)$ coincides with the forgetful functor
$$\mathrm{CAlg(\mathcal C)}_{I/}\longrightarrow \mathrm{CAlg(\mathcal C)}.$$ 
But $I$ is an initial object  of $ \mathrm{CAlg(\mathcal C)}$. $\square$
\subsection*{Monoidal adjunction.} 
Let $C^\otimes$ and $D^\otimes$ be symmetric monoidal $\infty$-categories and ${F:C\to D}$ be a functor. Suppose that $F$ is monoidal i. e. $F$ can be extended to the cocartesian fibration of $\infty$-operads
\begin{equation}\label{key}
	\begin{tikzcd}
		C^\otimes \arrow[rd] \arrow[rr, "F^\otimes"] &        & D^\otimes \arrow[ld] \\
		& \Delta &                     
	\end{tikzcd}
\end{equation}
Suppose that $F$ admits right-adjoint $R$. Then by Corollary 7.3.2.7. of \cite{key-1} $F^\otimes$ is also admits right-adjoint $G^\otimes$. Moreover, $G^\otimes$ is a map  of $\infty$-operads (so $G^\otimes$ is lax-monoidal). 
\\Now, let $\cO^\otimes$ be an $\infty$-operad. Then applying $\rm Alg_\cO(-)$ we get the pair of functors 
$$F:\mathrm{Alg}_\cO(C)\leftrightarrows \mathrm{Alg}_\cO(D):G.$$
We claim that $F$ is left-adjoint to $G$. Indeed, by  Joyal (Definition 4.0.1 of \cite{key-6}), we can define the adjunction for an  $(\infty,2)$-category
$\cC$ using a unit and counit satisfying the triangle identities. By Remark 4.4.5. of \cite{key-6}, in the case
$\cC=\rm Cat_\infty$ we get the Lurie's definition of adjunction \cite{key-3}. So applying $\rm  Alg_\cO(-)$ to the unit and counit of adjunction $(F^\otimes,G^\otimes)$ we get the unit and counit of adjunction $(F,G)$ which satisfy the triangle identities by functoriality.
Suppose that $f:\cO_1^\otimes \ra \cO_2^\otimes$ is a map of operads. Then we get the natural transformation $f^*:\rm Alg_{\cO_1} \ra Alg_{\cO_2}$. So $F$ and $G$ commute with $f^*$ and $f^*$ maps the (co)unit of adjunction to the (co)unit of adjunction.  
\begin{rem} \label{MonAdj}
	In particular, for  $\cO^\otimes:= \bE_\infty$ we get the adjunction 
	$$F:\mathrm{ CAlg}(F)\leftrightarrows \mathrm{CAlg}(D):G.$$
	Applying $\rm Alg_{-}(-)$ to the map of operads $\rm Triv^\otimes \ra \bE_\infty$ we conclude that the forgetful functors $U_C:\rm CAlg\it (C)\ra C$ and $U_D:\rm CAlg \it (D)\ra D$ commutes with $F$ and $G$. So if $\mathcal A$ be a $\bE_\infty$-algebra then the unit of adjunction $U_C\mathcal A\ra GF(U_C\mathcal A)$ can be extended to the homomorphism of algebras. The similar statement is true for the counit.  
\end{rem}
\subsection*{Free algebra of pointed object}
\begin{thm}\label{T2.6}
	Suppose that infinity direct sums is exist in $\mathcal{C}$ and for any object $c$ the functor $\otimes c$ commutes with colimits. Then 
	\begin{itemize}
		\item[1)] The forgetful functor $\mathcal U$ admits left-adjoint $\mathrm{Sym}^*$. Moreover, for any object $I\overset{v}{\rightarrow}A$ we have $$\mathrm{Sym}^*(A)=I\otimes_{S^*(I)}S^*(A).$$
		\item[2)] let us consider the diagram 
		\begin{equation}\label{F2.2}
			...\rightarrow (A^{\otimes n})_{h\Sigma_n} \rightarrow (A^{\otimes n+1})_{h\Sigma_{n+1}}\rightarrow ...	
		\end{equation} 
		where the maps induced by $v$. Then there is the weak equivalence in $\mathcal{C}$:
		$$\mathrm{Sym}^*(A)\simeq \mathrm{colim}_n (A^{\otimes n})_{h\Sigma_n}.$$
	\end{itemize}
\end{thm} 
\proof: Let us prove the first statement. Let us denote by $x$ the map induced by the unit of adjunction $I\rightarrow S^*(I)$. Note $I$ has the natural structure of $\mathbb E_\infty$-algebra. Let us denote by $\mathrm{ev_1}$ the map $S^*(I)\rightarrow I$ adjoint to $\mathrm{id}_I$. Observe that the composition $$I\overset{x}{\rightarrow}S^*(I)\overset{\mathrm{ev_1}}{\rightarrow}I$$ is equal to $\mathrm{id}_I$.
\\Now, let $CAlg(S^*(I)-mod)$ be the category of $S^*(I)$-algebras. Note that the functor $\mathcal U$ factors through $CAlg(S^*(I)-mod)$. Indeed, we have $CAlg(\mathrm{C}_{I/})\simeq CAlg(\mathrm{C})$ by Lemma \ref{L2.4} and $CAlg(S^*(I)-mod)\simeq CAlg(\mathrm{C})_{S^*(I)/}$ by Corollary 3.4.1.7 of \cite{key-1}. The map $\mathrm{ev_1}$ induce the functor $$CAlg(\mathcal{C})\overset{\mathrm{{ev_1}^*}}{\longrightarrow}CAlg(\mathcal{C})_{S^*(I)/}$$ $$(I\overset{1}{\rightarrow}A)\mapsto(S^*(I)\overset{\mathrm{ev_1}}{\rightarrow} I\overset{1}{\rightarrow}A)$$
The map $x$ induce the functor $u'$:$$CAlg(\mathcal{C})_{S^*(I)/}\overset{U_{S^*(I)/}}{\longrightarrow}\mathcal{C}_{S^*(I)/}\overset{x^*}{\longrightarrow}\mathcal{C}_{I/}.$$ $$(S^*(I)\rightarrow B)\mapsto (I\overset{x}{\rightarrow}S^*(I)\rightarrow B).$$ But $\mathrm{ev_1}\cdot x=\mathrm{id}_I$ so  $u'\cdot \mathrm{ev_1}^*\simeq \mathcal{U}$. 
\\On the other hand, by results of \cite{key-1}(see Theorem 4.5.3.1 and Corollary 4.2.3.7) and Remark \ref{MonAdj}  the functor $\mathrm{ev_1}$ is right-adjoint to relative tensor product $I\otimes_{S^*(I)}-$. The functor $u'$ is also admits left-adjoint which given by $S^*$. This is the particular case of slicing of adjoint functors (see Proposition 5.2.5.1 of \cite{key-1}). So the functor $ I\otimes_{S^*(I)}S^*(-)$ is the left-adjoint to $\mathcal{U}$.
\\Now, let us prove the second statement. Let $I\overset{v}{\rightarrow}A$ be an object of $\mathcal{C}_{I/}$. Note that the colimit of any diagram can be presented as the cofiber. In the case of the diagram \ref{F2.2} we get $$\mathrm{cofib}(S^*(A)\overset{\cdot (1-v)}{\longrightarrow} S^*(A)).$$ Let us rewrite $S^*(I)$ as $S^*(I)\otimes_{S*(I)} S^*(A)$. Then we have $$\mathrm{id}\otimes (1-v)=(1-x)\otimes \mathrm{id}.$$ So it suffice to check that the triangle $$S^*(I)\otimes_{S*(I)} S^*(A)\overset{(1-x)\otimes \mathrm{id}}{\longrightarrow} S^*(I)\otimes_{S*(I)} S^*(A) \longrightarrow I\otimes_{S*(I)} S^*(A)$$ is exact. 
\\Note that the functor $$-\otimes_{S^*(I)} S^*(A):\mathcal{C}\rightarrow \mathcal{C}$$ preserves colimits. Indeed, for any $S^*(I)$-algebra  $B$ the relative tensor product is given by the geometric realization of the two-sides bar construction $\mathrm{Bar}_{S^*(I)}(B,S^*(A))$ (see Section 4.4 of \cite{key-1}). Geometric realizations commute with colimits. On the other hand, colimits in the category of simplicial objects can be computed pointwise. Note that $$ \mathrm{Bar}_{S^*(I)}(B,S^*(A))_n=B\otimes S^*(I)^{\otimes n} \otimes S^*(A).$$ So the commutation of $-\otimes_{S^*(I)} S^*(A)$ with colimits follows from the commutation of $\otimes$ with colimits. 
\\Finally, observe that $\mathrm{cofib}(S^*(I)\overset{1-x}{\longrightarrow}S^*(I))\simeq I$. Indeed, $(I^{\otimes n})_{\Sigma_n}\simeq I$ so the diagram \ref{F2.2} for $A=I$ coincide with $...\overset{\mathrm{id}}{\rightarrow}I\overset{\mathrm{id}}{\rightarrow}I\overset{\mathrm{id}}{\rightarrow}...$ $\square$.

\subsection*{Relation with symmetric powers} Let $\mathcal C^\otimes$ be as above. Suppose that $\mathcal C^\otimes$ is stable, $\,\mathbb Q$-linear and idempotent-complete (in the sense of Section 4.4.5 of \cite{key-3}). Let $A$ be an object of $\mathcal C^\otimes$. Note that the n-th symmetric group $\Sigma_n$ acts on $A^{\otimes n}$. So we have the endomorphism $\varphi:=\frac{1}{n!}\Sigma_{\sigma_i\in\Sigma_{n}}\sigma_i$ of $A^{\otimes n}$. Now, define the n-th symmetric power $S^n(A)$ as the image $\mathrm{im}(\varphi).$ 
\begin{lem}\label{L2.7}
	The composition
	$$S^n(A)\hookrightarrow A^{\otimes n} \rightarrow (A^{\otimes n})_{h\Sigma_n}$$ is a weak equivalence. 
\end{lem}
Firstly, let us explain what do we mean by $(A^{\otimes n})_{h\Sigma_n}$. Note that the finite set $[n]$ induces the inclusion $\mathrm B \Sigma_n\hookrightarrow \mathrm{Fin_*}$. Applying $\mathrm{Fun_{N(Fin_*)}(-,\mathcal C^\otimes})$ we get the functor
$$\mathrm{CAlg(\mathcal C)\longrightarrow Rep_{\Sigma_n}(\mathcal C)=Fun(B\Sigma_n,\mathcal C)}$$
which maps $A$ to $A^{\otimes n}.$ Note that $\mathcal C$ can be considered as the full subcategory of constant functors. The corresponding inclusion $\mathcal C\rightarrow \mathrm{Rep_{\Sigma_n}(\mathcal C)}$ preserves limits so admits left-adjoint 
$$(-)_{h\Sigma_n}: \mathrm{Rep_{\Sigma_n}(\mathcal C)} \rightarrow \mathcal C$$ which is called homotopy coinvariants. 
\\Note that the element $\dfrac{1}{n!}\Sigma_{\sigma_i\in \Sigma_n}\sigma_i$ is an idempotent in the sense of Section 4.4.5 of \cite{key-3}. So we have the decomposition
$$A^{\otimes n}\simeq S^n(A)\oplus B$$
in $\mathcal C$ where $B=ker(\dfrac{1}{n!}\Sigma \sigma_i)$.
\\Note that in this case we have the fixed homotopies 
$$\frac{1}{n!}\Sigma \sigma_i|_{S^n(A)}\overset{\alpha}{\sim}id_{S^n(A)}\,\,\,\,\,\,\, , \frac{1}{n!}\Sigma \sigma_i|_B\overset{\beta}{\sim} 0.$$
Using this one can check that $S^n(A)$ and $B$ are $\Sigma_n$-equivariant objects. 
%$$\sigma|_{S^nA}\overset{\gamma}{\sim}\sigma|_{S^nA}\cdot \frac{1}{n!}\Sigma \sigma_i=\frac{1}{n!}\Sigma \sigma_i\overset{\gamma}{\sim}id_{S^nA}.$$ 
%$$\frac{1}{n!}\Sigma \sigma_i\cdot\sigma|_{B} =\frac{1}{n!}\Sigma \sigma_i|_B\overset{\gamma}{\sim}0.$$ 
\\\textit{The proof of Lemma:} Note that the action on  $S^n(A)$ is trivial. So $(S^n(A))_{h\Sigma_n}\simeq S^n(A)$. It suffice to show that $B_{h\Sigma_n}=0$. By stability of $\mathcal C$ is equivalent to 
$$\mathrm{H^0(Map_{\mathcal C}(B_{h\Sigma_n},c)\simeq H^0(Map_{Rep_{\Sigma_n}(\mathcal C)}(B,c)=0}$$ for any $\mathrm c\in \mathcal C$. By definition the element of $\rm H^0(Map_{Rep_{\Sigma_n}(\mathcal C)}(B,c)$ is a map $\phi:B\rightarrow c$ and homotopies $\phi\cdot\sigma\overset{\gamma_\sigma}{\sim}\phi$ for any $\sigma\in\Sigma_n $. So
$$0\overset{\beta}{\sim}\phi\cdot\Sigma \sigma_i\overset{\Sigma\gamma_{\sigma_i}}{\sim}n!\phi.\,\,\,\,\, \square$$

\section{Rational cohomology of a Kato-Nakayama space}

Let $(X,\alpha\!:\!\mathcal{O}_{X}\!\rightarrow\mathcal{M}_{X})$
be a fs log analytic space. For any point $x\in X$ let us define
the map $\mathrm{arg}\!:\!\mathcal{O}_{X,x}^{*}\longrightarrow S^{1}$;
$\mathrm{f}\longmapsto\mathrm{arg}(\mathrm{f}(x))$. Following \cite{'key-33},
let us construct the topological space $X^{\mathrm{log}}$. The point of $X^{\mathrm{log}}$ is the pairs $(x,\varphi)$
where $x$ is a complex point of $X$ and $\varphi:\mathcal{M}_{X,x}^{gr}\longrightarrow S^{1}$
is a map such that the following diagram is commutative
\[
\xymatrix{\mathcal{O}_{X,x}^{*}\ar[r]\ar[d]^{\mathrm{arg}} & \mathcal{M}_{X,x}^{gr}\ar[ld]^{\varphi\,\,\,\,\,\,\,.}\\
	S^{1}
}
\]
Notice that we have the projection $\pi:X^{\mathrm{log}}\longrightarrow X$.
For any open $U\subseteq X$ and $\mathbf{m}\in\mathcal{M}_{X}^{gr}(U)$
let us define the map 
\[
\mathrm{arg(\mathbf{m})}\!:\!\pi^{-1}(U)\longrightarrow S^{1}
\]
\[
(x,\varphi)\longmapsto\varphi(\mathrm{\mathbf{m}}_{x}).
\]
The topology on $X^{\mathrm{log}}$ is given the weak topology defined
by the map $\pi\!:\!X^{\mathrm{log}}\longrightarrow X$ and the family
of maps $\mathrm{arg(\mathbf{m})}$. The topological space $X^{\mathrm{log}}$
is called \textit{the Kato-Nakayama space of} $(X,\alpha)$. Note
that $X^{\mathrm{log}}$ does not depend on the map $\mathcal{O}_{X}\!\rightarrow\mathcal{M}_{X}$
, but only on the group completion $\mathcal{M}_{X}^{gr}$.

By definition, \textit{the Betti cohomology} of a fs log analytic
space $(X,\alpha)$ is a singular cohomology of $X^{\mathrm{log}}$.
So for any abelian group $\mathrm{A}$ we have 
\[
H_{Betti}^{*}((X,\alpha),\mathrm{\mathrm{A}})\overset{\mathrm{def}}{=}H^{*}(X,Rf_{*}\underline{\mathrm{\mathrm{A}}})\simeq\mathbb{H}^{*}(X,R\pi_{*}\underline{\mathrm{\mathrm{A}}})
\]
where $f$ is the map $X^{\mathrm{log}}\longrightarrow pt$.

Suppose we have a description of $R\pi_{*}\underline{\mathrm{\mathrm{A}}}$
in terms of the log structure on $X$. Then one could define the Betti
cohomology independently of the construction of $X^{\mathrm{log}}$.
The first approximation to such the description can be given by the
following lemma proved by Kato and Nakayama.
\begin{lem}
	(Lemma (1.5) of\textcolor{red}{{} }\cite{key-4}) For any abelian group
	$\mathrm{A}$ we have a canonical isomorphism
	\[
	R^{q}\pi_{*}\underline{\mathrm{\mathrm{A}}}\simeq\mathrm{A}\otimes_{\mathbb{Z}}\bigwedge^{q}\mathcal{\overline{M}}_{X}^{gr}.
	\]
\end{lem}
%\begin{rem}
%(\textit{the symmetric algebra of a pointed complex}) For a pointed
%complex of sheaves of $\mathbb{Q}$-vector spaces $\mathcal{K}^{\bullet}$
%(\textcolor{black}{i.e. the complex }$\mathcal{K}^{\bullet}$\textcolor{black}{{}
%together with the fixed arrow}\textbf{ $\mathbb{Q}^{\bullet}\!\overset{1_{\mathcal{K}^{\bullet}}}{\longrightarrow}\mathcal{K}^{\bullet}$})
%let us define the symmetric algebra $\mathrm{Sym}\mathcal{K}^{\bullet}=\underset{\longrightarrow}{\mathrm{lim}}\,S^{n}\mathcal{K}^{\bullet}$
%	as the inductive limit. Here the morphism $S^{n-1}\mathcal{K}^{\bullet}\hookrightarrow S^{n}\mathcal{K}^{\bullet}$
%\textcolor{black}{is induced by the map} $\mathcal{K}^{\bullet\otimes n-1}\otimes_{\mathbb{Q}}\mathbb{\mathbb{Q}^{\bullet}}\overset{id\otimes1_{\mathcal{K}^{\bullet}}}{\longrightarrow}\mathcal{K}^{\bullet}{}^{\otimes n}$.
%Note that the functor $\mathrm{Sym}$ is a left-adjoint to the forgetful
%	functor $\mathrm{Ch_{*}^{\bullet}(Sh}(X,\mathbb{Q}))\longrightarrow\mathbf{DG\!-\!CAlg}(\mathrm{Sh}(X,\mathbb{Q}))$
%	from the category of pointed complexes to the category of commutative
%	DG-algebras with unit.
%\end{rem}
Now, let us give the explicit description of $R\pi_{*}\underline{\mathrm{\mathrm{A}}}$
in the case $\mathrm{A=}\mathbb{Q}$. Let $(X,\alpha)$ be as above.
Let us denote by $\mathrm{exp}(\alpha)$ the complex 
\[
...\longrightarrow\underset{(-1)}{0_{\,_{\,}}}\longrightarrow\underset{(0)}{\mathcal{O}_{X_{\,}}}\longrightarrow\underset{(1)}{\mathcal{M}_{X}^{gr}}\longrightarrow\underset{(2)}{0_{\,_{\,}}}\longrightarrow...
\]
where the differential given by the composition $\mathcal{O}_{X}\!\overset{\mathrm{exp}}{\rightarrow}\mathcal{O}_{X}^{*}\!\hookrightarrow\mathcal{M}_{X}^{gr}$.
Note that the inclusion of zero cohomology $1_{\mathrm{exp}(\alpha)}:\bQ\to \mathrm{exp}(\alpha)_\bQ$ allows to consider $\mathrm{exp}(\alpha)_\bQ=\mathrm{exp}(\alpha)\otimes \bQ$ as an object of the under category $Ch^\bullet(Sh(X^{an},\bQ))_{\setminus \bQ}.$ Then the maps 
$$\mathrm{exp}(\alpha)_\bQ^{\otimes n}\otimes \bQ\sa{\mathrm{id}\otimes 1_{\mathrm{exp}(\alpha)}} \mathrm{exp}(\alpha)_\bQ^{\otimes n}$$
induce the inclusion 
$$S^{n-1}(\mathrm{exp}(\alpha)_\bQ)\to S^{n}(\mathrm{exp}(\alpha)_\bQ). $$
We put
\begin{equation}\label{Sym}
	\mathrm{Sym}^*(\mathrm{exp}(\alpha)_\bQ):=\mathrm{colim}_nS^n(\mathrm{exp}(\alpha)_\bQ).
\end{equation}
By results of the previous section (\ref{Sym}) is a free $\mathbb E_\infty$-algebra in $Ch^\bullet(Sh(X^{an},\bQ))_{\setminus \bQ}.$ 
\begin{thm}\label{R_pi_Q}
	There exist a natural quasi-isomorphism 
	of $\mathbb E_\infty$-algebras 
	$$ R\pi_{*}\mathbb{Q}\iso \mathrm{Sym}^*(\mathrm{exp}(\alpha)_\bQ) .$$
\end{thm}
Firstly, let us construct the quasi-isomorphism $\rm exp(\alpha)\iso \tau_{\leq 1}R\pi_*\bZ$. For a topological space $Y$ let us denote by $C(Y)$ the sheaf of continuous $Y$-valued functions on $X^{log}.$ Note that the sheaf $C(\bR)$ is soft. For $x\in X^{an}$ let us denote by $i$ and $\phi$ be the inclusions $x\hra X^{\rm an}$ and $\pi^{-1}(x)\hra X^{\rm log}.$ By Proposition 1.2.5. of \cite{'key-33} $X^{log}$ is a Hausdorff space. So $\pi$ is universally closed and we have the proper base change theorem.  Then 
$$i^{-1}R^1\pi_*C(\bR)\iso R^1\pi_*(\phi^{-1}C(\bR))\iso H^1(\pi^{-1}(x),C(\bR))=0$$ for any $x\in X^{\rm an}.$
\begin{lem}
	Let $$\bZ \ra A \ra B$$ be an exact sequence of sheaves and $R^1\pi_*A=0$. Then $\tau_{\leq 1}R\pi_*\bZ\iso \pi_*(\rm Cone(\it A\ra B)).$
\end{lem}
\proof Note that Cone commutes with $\pi_*$ so we have the diagram with exact rows 
$$\begin{tikzcd}
	{\pi_*B[-1]} \arrow[d] \arrow[r] & \pi_*\rm Cone(\it A\to B\rm) \arrow[d] \arrow[r]  & \pi_*A \arrow[d] \\
	{R\pi_*B[-1]} \arrow[r]          & R\pi_*\mathbb Z \arrow[r] & R\pi_*A         
\end{tikzcd}$$ 
It remains to use the long exact sequences of cohomology $\square$.
\\Let $\rm exp:\it C(\bR)\ra C(S^\mathrm{1})$ be the map induces by 
$$\rm exp:\it \bR\ra S^{\rm 1},\,\,\,\, \it a\mapsto e^{2\pi ia}.$$
Notice that the sequence of sheaves 
$$\bZ\ra C(\bR)\sa{exp} C(S^1)$$ 
is exact.
\begin{rem}
	For $\rm U\in X^{log}$ let us compute the cokernel of homomorphism ${\rm Map(U,\bR)\ra Map(U,\it S^{\rm 1}).}$ Note that a map $\phi\in \rm Map(U,\it S^{\rm 1})$ can be lifted to $\rm Map(U,\bR)$ iff $[\phi]=0\in [U,S^1]\iso H^1(U,\bZ)$. Indeed, let $\underline 1:X\ra S^1$ be a constant map with the value 1. Suppose that $\phi$ is homotopy equivalent to $\underline{1}$. But $exp(\underline{0})=\underline 1$ so we can use the homotopy lifting properties:
	$$
	% https://tikzcd.yichuanshen.de/#N4Igdg9gJgpgziAXAbVABwnAlgFyxMJZABgBpiBdUkANwEMAbAVxiRAA0QBfU9TXfIRRkAjFVqMWbdgB0ZeALbwABAEluvEBmx4CREeXH1mrRCDkK6OABYAjW8oBKGvjsH7SY6salmAygB6ItziMFAA5vBEoABmAE4QCkgATNQ4EEgAzN6SpiAwAB5oLiDxiUgGIOlZOSZscmjWWMpy2AotMkxgsHEMWGAwysE8sQlJiGRVGYgiI6VjSJPViKkSdWZyXT19A8rEJWXjlcurO3lQdHDWYSFcQA
	\begin{tikzcd}
		X \arrow[d] \arrow[r, "\underline 0"]                            & \mathbb R \arrow[d, "exp"] \\
		X\times I \arrow[r, "\phi \sim \underline 1"] \arrow[ru, dashed] & S^1.                       
	\end{tikzcd}
	$$
	Hence the sheaf $R^1\pi_*\bZ$ is a sheafification of the presheaf $W\mapsto H^1(\pi^{-1}W,\bZ).$
\end{rem}
By the Lemma we have the quasi-isomorphism $\tau_{\leq 1}R\pi_*\bZ\iso \pi_*\rm Cone(exp).$ Let us define the map of complexes
\begin{equation}\label{F3.1}
	\begin{tikzcd}
		\mathcal O \arrow[d] \arrow[r, "exp"] & \mathcal M^{\rm gr} \arrow[d] \\
		\pi_* C(\mathbb R) \arrow[r, "exp"]         & \pi_* C(S^1)                       
	\end{tikzcd}	
\end{equation}
by the rule 
$f\in \cO(U) \mapsto \pi^*\rm Re\it f  ,\,\, \,\,
h\in\cM(U)\mapsto \{\, h^*:(x,\phi_x)\mapsto \phi_x(h_x)\,\}.$ 
\begin{lem}
	The map \ref{F3.1} induces quasi-isomorphism 
	$$\rm exp(\alpha)\iso \tau_{\leq 1}R\pi_*\bZ.$$
\end{lem}
\proof It suffice to check that the induced map $\overline {\mathcal M}^{gr}_X\ra R^1\pi_*\bZ$
is an isomorphism. Let $x\in X^{\rm an}$ be a closed point. Let us consider the commutative diagram 
$$
\xymatrix@=7pt{ & colim_{x\in V}C(\pi^{-1}V,\mathbb{R})\ar@{->}[rr]^{exp}\ar@{-}[d] &  & colim_{x\in V}C(\pi^{-1}V,S^{1})\ar@{->}[dd]^{Res}\\
	colim_{x\in V}\mathcal{O}(V)\ar@{->}[rr]\ar@{->}[dd]^{Res}\ar@{->}[ru] &  & colim_{x\in V}\mathcal{M}^{gr}(V)\ar@{->}[dd]\ar@{->}[ru]\\
	& C(\pi^{-1}(x),\mathbb{R})\ar@{-}[r]\ar@{<-}[u] &  & C(\pi^{-1}(x),S^{1})\ar@{<-}[l]\\
	\mathbb{C}\ar@{->}[rr]^{exp}\ar@{->}[ru] &  & \,\,\mathbb{C}^{*}\oplus\mathcal{M}_{x}^{gr}\,\,\ar@{->}[ru]
}
$$
where $Res$ are the restriction maps. Applying $H^1$ we get the commutative diagram 
$$
\begin{tikzcd}
	colim_{x\in V}\mathcal M^{gr}(V) \arrow[rr] \arrow[d, "\wr"] &  & {colim_{x\in V}H^1(\pi^{-1}V,\mathbb Z)} \arrow[d, "Res"] \\
	\mathcal M^{gr}_x \arrow[rr]                                 &  & {H^1(\pi^{-1}(x),\mathbb Z)}                             
\end{tikzcd}
$$
Let $\tilde{\pi}$ be the map $\pi^{-1}(x)\to x.$ Then by Proper base change theorem the restriction map 
$$(R^1\pi_*\bZ)_x=\rm colim_{x\in V}H^1(\pi^{-1}V,\bZ)\sa{Res} H^1(\pi^{-1}(x),\bZ)= R^1\tilde{\pi}_*\bZ$$
is an isomorphism. So we may assume that $X^{an}$ is a point.
\\In this case the map $\overline {\mathcal M}^{gr}_X\ra R^1\pi_*\bZ$ can be expressed as the composition 
\begin{equation}\label{F3.2}
	\rm \overline {\mathcal M}^{gr}\iso \bZ^r\ra Map((\mathit{S}^1)^r,\mathit{S}^1)\ra [(\mathit{S}^1)^r,\mathit{S}^1]
\end{equation}
$$(a_1,..,a_r)\mapsto ((\phi_1,..,\phi_r)\mapsto (\phi_1^{a_1},..,\phi_r^{a_r})).$$
Note that the map \ref{F3.2} is defined by the images of (0,..,1,..,0) which correspond to the projections
$$\pi_i: S^1\times ... \times S^1\ra S^1.$$
It remains to note that $[\pi_i]$ generate $H^1((S^1)^r,\bZ).$ $\square$
\\ \textit{The proof of Theorem \ref{R_pi_Q}:} Let $A^\otimes, B^\otimes$ be a symmetric monoidal categories and $F^\otimes: A^\otimes \leftrightarrows B^\otimes :G^\otimes$ be a monoidal adjunction. Then by Remark \ref{MonAdj} we have the pairs of adjoint functors 
$$F:A_{\setminus I_A} \leftrightarrows B_{\setminus I_B}:G$$
$$F:CAlg(A_{\setminus I_A}) \leftrightarrows CAlg(B_{\setminus I_B}):G.$$ 
Suppose that $\mathcal A\in CAlg(A_{\setminus I_A})$ is an algebra and  $a\sa{\varphi}\mathcal A$, $Sym(a)\sa{\psi}\mathcal A$ is a pair of adjoint morphisms.
Then $F(\varphi),\, F(\psi)$ are also adjoint. Indeed, $F$ commutes with the forgetful functors $U_A:CAlg(A_{\setminus I_A})\ra A_{\setminus I_A}$ and $U_B:CAlg(B_{\setminus I_B})\ra B_{\setminus I_B}$ by Remark \ref{MonAdj}. So it suffice to show that $F$ preserves the unit of adjunction $\mathrm{Id} \ra U\cdot \rm Sym^*.$ Let $a\in A_{\setminus I_A}$. Then the map $a\ra \mathrm{Sym}^*(a)$ can be expressed as the canonical morphism
\begin{equation}\label{F3.3}
	a\to colim_n (a^{\otimes n})_{h\Sigma_n}.
\end{equation} 
Note that $F$ is monoidal and preserves colimits. So $F$ commutes with $(-^{\otimes n})_{h\Sigma_n}$ and maps the canonical morphism \ref{F3.3} to the canonical morphism.
\\Now, let $\tau_{\leq 1}R\pi_*\bQ\ra R\pi_*\bQ$ be the canonical inclusion. Note that $R\pi_*$ is lax-monoidal so $R\pi_*\bQ$ is an algebra with the unit $R^0\pi_*\bQ=\bQ.$ So we get the map 
\begin{equation}\label{F3.4}
	\mathrm{Sym}^*(\tau_{\leq 1}R\pi_*\bQ)\ra R\pi_*\bQ.	
\end{equation} 
To complete the proof, it suffices to show that the map is a quasi-isomorphism on any fiber. Let $i:x\hookrightarrow  X^{an}$ be a point. By the previous $\mathrm{Sym}^*(\tau_{\leq 1}R\pi_*\bQ)_x\ra (R\pi_*\bQ)_x$ is adjoint to $\tau_{\leq 1}R\pi_*\bQ_x\ra R\pi_*\bQ_x.$ So using the proper base change we may assume that $X^{an}$ is a point. In this case, $X^{log}$ is a torus, $exp(\alpha)\iso \bZ\oplus \bZ^n[-1]$ for some $n$ and $\mathrm{Sym}^*(exp(\alpha))_\bQ\iso \Lambda^*(\bQ^n)$.
So (\ref{F3.4}) is the map $\Lambda^*(\bQ^n) \to R\pi_*\bQ$ induced by the inclusion of the one-dimension cycles $\bQ^n[-1]\to R\pi_*\bQ.$ It is well known that such a map is a quasi-isomorphism. $\square.$
% One can apply the Remark to the functor $F=\oplus H^*$ which maps the complex A to the complex $\oplus H^i(A)$. Then it should be proven that $\Lambda^*(\bQ^n) \to R\pi_*\bQ$ is a quasi-iso.
\section{From log structures to motivic sheaves}

\subsection*{Extension on a smooth étale site.}

Let $X$ be a scheme of finite type over $\mathrm{k}$. The natural inclusion of site $X_{et}\hookrightarrow\mathrm{Sm/}X_{et}$
give rise to the pair of functors
\[
\mathrm{Ex}_{X}:\mathrm{Sh}(X_{et},\bZ)\rightleftarrows\mathrm{Sh_{et}(Sm/}X,\bZ):\mathrm{Res}_{X}.
\]
They are given as follows: If $\mathcal{F}\in\mathrm{Sh_{et}(Sm}/X,\bZ)$
, we have $\mathrm{Res}_{X}\mathcal{F}(X^{\prime})=\mathcal{F}(X^{\prime})$
for étale $X$-schemes $X^{\prime}$. If $\mathcal{G}$ is a sheaf on
$X_{et}$, then the sheaf $\mathrm{Ex}_{X}(\mathcal{G})$ is associated
to the presheaf 
\[
Z\mapsto\underset{\rightarrow}{\lim}\,\mathcal{G}(X^{\prime})
\]
on $\mathrm{Sm/}X$, where the limit extends over diagrams of the
form 
\[
\xymatrix{Z\ar[dr]\ar@{-->}[r] & X^{\prime}\ar@{-->}[d]^{et}\\
	& X
}
\]
By results of (\cite{'key-5}, Section VII, 4.0, 4.1) we have  
\begin{lem}
	The functor $\mathrm{Ex}_X$ is left adjoint to $\mathrm{Res}_X$ and exact. Moreover, $\mathrm{Ex}_X$ commutes with the inverse image functors $f^{-1}$ for any $f$.
\end{lem}
Now, we want to extend any virtual log structure on $X_{et}$ to the étale sheaf on $\mathrm{Sm/}X_{et}$. Namely, for the log structure $\mathcal{M}^{gr}$
on $X$ and smooth $f\!:\!Y\longrightarrow X$ let us define the sheaf
\[
\mathrm{Ex_{\mathit{X}}^{log}}\!\mathcal{M}^{gr}(\begin{array}{c}
	Y\\
	\,\downarrow^{f}\\
	X
\end{array})\overset{\mathrm{def}}{=}\Gamma(Y,f^{*}\mathcal{M}^{gr})
\]
where $f^{*}$ is a pullback of log structures. 

Observe that $\mathrm{Ex_{\mathit{X}}^{log}}$ can be defined by the
more functorial way. Let
$$\mathcal{\underline{M}}_{X}^{gr},\mathcal{\underline{O}}_{X}^{*}:\mathbf{vLog}_{X}\longrightarrow\mathrm{Sh(}X_{et},\bZ)$$
be the functors which maps the log structure $\mathcal{O}^{*}\hookrightarrow\mathcal{M}^{gr}$
to the sheaf $\mathcal{M}^{gr}$ and $\mathcal{O}^{*}$. Let us denote
by $\mathrm{Ex}_{X}(\mathcal{\underline{M}}_{X}^{gr})$ and $\mathrm{Ex}_{X}(\mathcal{\underline{O}}_{X}^{*})$
the compositions with $\mathrm{Ex}_{X}$. Observe that we have the natural
transformation $$\mathrm{Ex}_{X}(\mathcal{\underline{O}}_{X}^{*})\longrightarrow\mathrm{Ex}_{X}(\mathcal{\underline{M}}_{X}^{gr}).$$
Also, let us denote by $\mathcal{\underline{O}}_{\mathrm{Sm/}X}^{*}$
the constant functor $\mathbf{vLog}_{X}\longrightarrow\mathrm{Sh_{et}(Sm/}X,\bZ)$
with value $\mathcal{O}_{\mathrm{Sm/}X}^{*}$ (here $\mathcal{O}_{\mathrm{Sm/}X}^{*}$
is an étale sheaf of inverse function on $\mathrm{Sm/}X$). For any
morphism of scheme $f\!:\!X\longrightarrow Y$, the natural map $f^{-1}\mathcal{O}_{Y}^{*}\longrightarrow\mathcal{O}_{X}^{*}$
gives rise to the natural transformation $$\mathrm{Ex}_{X}(\mathcal{\underline{O}}_{X}^{*})\longrightarrow\mathcal{\underline{O}}_{\mathrm{Sm/}X}^{*}.$$
Now we can define $\mathrm{Ex^{log}}$ as the pushout
\[
\xymatrix{\mathrm{Ex}_{X}(\mathcal{\underline{O}}_{X}^{*})\ar[r]\ar[d] & \mathrm{Ex}_{X}(\mathcal{\underline{M}}_{X}^{gr})\ar[d]^{\,\,\,\,\,\,\,.}\\
	\mathcal{\underline{O}}_{\mathrm{Sm/}X}^{*}\ar[r] & \mathrm{Ex_{\mathit{X}}^{log}}
}
\]
Notice that the map $\mathrm{Ex}_{X}(\mathcal{\underline{O}}_{X}^{*})(\varphi)\longrightarrow\mathrm{Ex}_{X}(\mathcal{\underline{M}}_{X}^{gr})(\varphi)$
is an inclusion for any virtual log structure $\varphi$. So the map
$\mathrm{Ex}_{X}(\mathcal{\underline{O}}_{X}^{*})\longrightarrow\mathrm{Ex}_{X}(\mathcal{\underline{M}}_{X}^{gr})$
is an inclusion in the abelian category $\mathrm{Fun}(\mathbf{vLog}_{X},\mathrm{Sh_{et}(Sm/}X,\bZ))$
and we have the short exact sequence 
\[
\mathrm{Ex}_{X}(\mathcal{\underline{O}}_{X}^{*})\rightarrow\mathcal{\underline{O}}_{\mathrm{Sm/}X}^{*}\oplus\mathrm{Ex}_{X}(\mathcal{\underline{M}}_{X}^{gr})\rightarrow\mathrm{Ex_{\mathit{X}}^{log}}.
\]

\subsection*{The natural transformation $\mathrm{Ex^{log}}$.}

Let $f:Y\longrightarrow X$ be a morphism of schemes. Let us construct
the functorial map $\psi_{f}:f^{-1}\mathrm{Ex_{\mathit{X}}^{log}}\longrightarrow\mathrm{Ex_{\mathit{Y}}^{log}}f^{*}$.
Recall that we defined $f^{*}\mathcal{M}_{X}^{gr}$ as the cocartesian
square. So applying $\mathrm{Ex}$ and using commutation $\mathrm{Ex}$
with $f^{-1}$ we get the square 
\begin{equation}\label{'F5.1}
	\xymatrix{f^{-1}\mathrm{Ex}_{X}(\mathcal{O}_{X}^{*})\ar[r]\ar[d] & f^{-1}\mathrm{Ex}_{X}(\mathcal{M}_{X}^{gr})\ar[d]\\
		\mathrm{Ex}_{Y}(\mathcal{O}_{Y}^{*})\ar[r] & \mathrm{Ex}_{Y}(f^{*}\mathcal{M}_{X}^{gr})
	}
\end{equation}
On the other hand, let us consider the diagram
\[
\xymatrix{\mathcal{O}_{X}^{*}\ar[r]\ar[d]^{id} & f_{*}\mathcal{O}_{Y}^{*}\ar[d]^{id}\\
	\mathrm{Res}_{X}(\mathcal{O}_{\mathrm{Sm/}X}^{*})\ar[r] & f_{*}\mathrm{Res}_{Y}(\mathcal{O}_{\mathrm{Sm/}Y}^{*})
}
\]
where the horizontal maps induced by $f$. Using commutation $\mathrm{Res}$
with $f_{*}$ and adjunction $\mathrm{Res}-\mathrm{Ex}$, $f_{*}-f^{-1}$
we get the commutative diagram
\[
\xymatrix{f^{-1}\mathrm{Ex}_{X}(\mathcal{O}_{X}^{*})\ar[r]\ar[d] & \mathrm{Ex}_{Y}(\mathcal{O}_{Y}^{*})\ar[d]\\
	f^{-1}\mathcal{O}_{\mathrm{Sm/}X}^{*}\ar[r] & \mathcal{O}_{\mathrm{Sm/}Y}^{*}
}
\]
So we have the commutative diagram 
\begin{equation}
	\xymatrix{f^{-1}\mathrm{Ex}_{X}(\mathcal{O}_{X}^{*})\ar[r]\ar[d] & f^{-1}\mathcal{O}_{\mathrm{Sm/}X}^{*}\oplus f^{-1}\mathrm{Ex}_{X}(\mathcal{M}_{X}^{gr})\ar[d]\\
		\mathrm{Ex}_{Y}(\mathcal{O}_{Y}^{*})\ar[r] & \mathcal{O}_{\mathrm{Sm/}Y}^{*}\oplus\mathrm{Ex}_{Y}(f^{*}\mathcal{M}_{X}^{gr})
	}
\end{equation}
and we define $\psi_{f}:f^{-1}\mathrm{Ex_{\mathit{X}}^{log}}(\mathcal{M}_{X}^{gr})\longrightarrow\mathrm{Ex_{\mathit{Y}}^{log}}(f^{*}\mathcal{M}_{X}^{gr})$
as the map of cokernels. 
\begin{prop}
	The map $\psi_{f}$ is an isomorphism.
\end{prop}
\textit{\negmedspace{}\negmedspace{}\negmedspace{}\negmedspace{}\negmedspace{}Proof:}
Let us consider the diagram \ref{'F5.1}. By the exactness
of $\mathrm{Ex}$ the square \ref{'F5.1} is cocartesian
and the horizontal maps are inclusions. So the map between the cokernels
of horizontal arrows is an isomorphism. By Snake lemma we immediately
get that the canonical maps 
\[
\mathrm{ker}(f^{-1}\mathrm{Ex}_{X}(\mathcal{O}_{X}^{*})\rightarrow\mathrm{Ex}_{Y}(\mathcal{O}_{Y}^{*}))\longrightarrow\mathrm{ker}(f^{-1}\mathrm{Ex}_{X}(\mathcal{M}_{X}^{gr})\rightarrow\mathrm{Ex}_{Y}(f^{*}\mathcal{M}_{X}^{gr}))
\]
\[
\mathrm{coker}(f^{-1}\mathrm{Ex}_{X}(\mathcal{O}_{X}^{*})\rightarrow\mathrm{Ex}_{Y}(\mathcal{O}_{Y}^{*}))\longrightarrow\mathrm{coker}(f^{-1}\mathrm{Ex}_{X}(\mathcal{O}_{X}^{*})\rightarrow\mathrm{Ex}_{Y}(\mathcal{O}_{Y}^{*}))
\]
are isomorphisms. So 
\[
\mathrm{ker}(\psi_{f})\simeq\mathrm{ker}(f^{-1}\mathcal{O}_{\mathrm{Sm/}X}^{*}\longrightarrow\mathcal{O}_{\mathrm{Sm/}Y}^{*})
\]
\[
\mathrm{coker}(\psi_{f})\simeq\mathrm{coker}(f^{-1}\mathcal{O}_{\mathrm{Sm/}X}^{*}\longrightarrow\mathcal{O}_{\mathrm{Sm/}Y}^{*})
\]
and we can use Corollary \ref{CA.2}. $\Square$
\begin{cor}\label{'C5.3}
	Then the collection of functors $\mathrm{Ex_{\mathit{X}}^{log}}$
	gives rise to the natural transformation of contravariant pseudofunctors
	\[
	\mathrm{Ex^{log}}\!:\!\mathbf{vLog}\longrightarrow\mathrm{Sh_{et}(Sm/}-,\mathbb{Q}).
	\]
\end{cor}

\subsection*{From log structures to motivic sheaves. }

Let R be a commutative ring and $S$ be a scheme. Let us denote by $\mathrm{DA_{et}^{eff}}(S,\mathrm R)$ and $\mathrm{DA_{et}}(S,\mathrm R)$
be the categories of effective and spectral étale motivic sheaves
on $S$ (see \cite{'key-19} for definition, see also \cite{'key-18}
and \cite{'key-22} for more details). Notice that for any morphism
$f\!:\!S'\longrightarrow S$ we have the inverse image functor $f^{*}\!:\!\mathrm{DA_{et}}(S,\mathrm R)\longrightarrow\mathrm{DA_{et}}(S',\mathrm R)$
which can be define as the Kan extension of the composition
\[
\xymatrix{\mathrm{Sm/}S\ar[r]^{-\times_{S}S'} & \mathrm{Sm/}S'\ar[r] & \mathrm{DA_{et}}(S',\mathrm R).}
\]
We can also define the inverse image functor $f^{*}\!:\!\mathrm{DA_{et}^{eff}}(S,\mathrm R)\longrightarrow\mathrm{DA_{et}^{eff}}(S',\mathrm R)$
by the same way. Moreover, by construction (see \cite{'key-22}) $\mathrm{DA_{et,R}}$
and $\mathrm{DA_{et,R}^{eff}}$ can be considered as the
presheaves of $\infty$-categories 
\[
S\longmapsto\mathrm{DA_{et}}(S,\mathrm R);\,\,\,\,S\longmapsto\mathrm{DA_{et}^{eff}}(S,\mathrm R)
\]
\[
(f \!: \! S'\to S)\longmapsto f^{*}.
\]
\begin{rem}\label{'R5.4}
	For any $S\in\mathrm{Sch/k}$ we have the $\mathbb{A}^{1}$-localization
	functor 
	\[
	\mathrm{L_{\mathbb{A}^{1}}}\!:\!\mathrm{D(Sh_{et}(Sm}/S,\mathrm R))\longrightarrow\mathrm{DA_{et}^{eff}}(S,\mathrm R)
	\]
	and the infinity $\Gm$-suspension functor
	\[
	\mathrm{\Sigma_{\Gm}^{\infty}}\!:\!\mathrm{DA_{et}^{eff}}(S,\mathrm R)\longrightarrow\mathrm{DA_{et}}(S,\mathrm R).
	\]
	%Here $\mathbb{L}$ is the Lefschetz motive. It is the motive that
	%corresponds to the cokernel of the natural inclusion $\mathbb{Q}[\infty_{X}]\hookrightarrow\mathbb{Q}[\mathbb{P}_{X}^{1}]$. 
	By construction $\mathrm{L_{\mathbb{A}^{1}}}$ and $\mathrm{\mathrm{\Sigma_{\Gm}^{\infty}}}$
	commute with inverse image functors. So we have the composition of
	natural transformations
	\[
	\mathrm{\mathrm{\Sigma_{\Gm}^{\infty}}}\!\cdot\!\mathrm{L_{\mathbb{A}^{1}}}\!:\!\mathrm{D(Sh_{et}(Sm/}-,\mathrm R))\longrightarrow\mathrm{DA_{et,\mathrm R}}
	\]
\end{rem}

\begin{defn}\label{'D5.5}
	Let $X=(\uX,\cM_X)$ be a fs log scheme. Let us denote by $\rm M^{gr}_{\it X}$ the image of $\cM_{X}$ under the composition 
	\begin{equation}\label{'F5.3}
		\xymatrix{\mathbf{vLog}^{\mathrm{fs}}_{\uX}\ar[rr]^{\!\!\!\!\!\!\!\!\! \mathrm{Ex^{log}\otimes}\mathrm{R}} &  & \mathrm{Sh_{et}(Sm/}\underline{X},\mathrm{R})\ar@{^{(}->}[r]^{\heartsuit} & \mathrm{D}(\mathrm{Sh_{et}(Sm/}\underline{X},\mathrm{R}))\ar[r]^{\ \ \ \ \ \mathrm{\mathrm{\Sigma_{\mathbb{G}_m}^{\infty}}}\cdot\mathrm{L_{\mathbb{A}^{1}}}} & \mathrm{DA_{et}}(\uX,\mathrm R).}
	\end{equation}
\end{defn}
The motivic sheaves $\rm M^{gr}_{\mathit X}$ can be considered as a motivic analog of log structures. They play a critical role further in the text.

\section{The homological motive of a log scheme}
%Let $f:X\to Y$ be the morphism of schemes and $f^{-1}$ be the inverse functor $Sh_{et}(Sm/Y,\bQ)\to Sh_{et}(Sm/X,\bQ)$. By Corollary A.2?, the canonical map $f^{-1}\cO^*_{Sm/X}\to \cO^*_{Sm/Y}$ is an isomorphism. 

Now, let $(X,\cM)$ be a virtual log scheme. Let $R=\bQ$. By the previous, the inclusion $\cO^*\hookrightarrow \cM$ induces the map $\bQ(1)[1]\to \rm M^{gr}.$  It makes $\rm M^{gr}(-1)[-1]$ an object of the under category $DA_{et}(X,\bQ)_{/\bQ}.$ So we have a free algebra
$$\bQ^{log}_{X}:=\mathrm{Sym}^{*}({\rm M^{gr}}(-1)[-1]).$$
Let $f$ be the canonical morphism $X\to Spec(k).$ By analogy with Theorem \ref{R_pi_Q}, we define \textit{the homological motive of} $(X,\cM)$ by the formula 
$$[X]^{log}:=\underline{Hom}(Rf_*\bQ^{log}_X,\bQ).$$
\begin{rem}
	In \cite{'key-23} Deligne defined for any $\bQ$-linear, monoidal and  idempotent complete category the notions of Schur functors. In particular, we have well-defined motivic symmetric powers $S^n:DA_{et}(X,\bQ)\to DA_{et}(X,\bQ).$ By Lemma \ref{L2.7}, we get the weak equivalence
	$$\mathrm{Sym}^{*}({\rm M^{gr}}(-1)[-1])\iso \mathrm{colim}_n S^n({\rm M^{gr}}(-1)[-1]).$$
	So the motivic sheaf $\bQ^{log}_{X}$ can be defined without using of any higher algebra.
\end{rem}
Let $(g,g^\#):(X,\cM_{X})\to (Y,\cM_{Y})$ be the morphism of virtual log schemes. Let us construct the map $[X]^{log}\to [Y]^{log}.$ 
Firstly, note that (\ref{'F5.3}) maps $g^*\cM_Y$ to $g^*\rm M^{gr}_{\it Y}$ by Corollary \ref{'C5.3} and Remark \ref{'R5.4}. Hence $g^{\#}$ maps to 
$$\!\mathrm{S^{q}}(g^{*}\mathrm{M_{\mathit{Y}}^{gr}(-1)[-1])}\longrightarrow\mathrm{S^{q}(M_{\mathit{X}}^{gr}(-1)[-1])}.$$ 
On the other hand, $g^{*}$ is monoidal so $\mathrm{S^{q}}g^{*}(\mathrm{M_{\mathit{Y}}^{gr}(-1)[-1])}\simeq g^{*}\mathrm{S^{q}(M_{\mathit{Y}}^{gr}(-1)[-1])}$
and we get the map
\[
g^*\mathrm{Sym^{*}(M_{\mathit{Y}}^{gr}(-1)[-1])}\longrightarrow \mathrm{Sym^{*}(M_{\mathit{X}}^{gr}(-1)[-1])}.
\]
Then, by adjunction, we have
\[
\mathrm{Sym^{*}(M_{\mathit{Y}}^{gr}(-1)[-1])}\longrightarrow Rg_*\mathrm{Sym^{*}(M_{\mathit{X}}^{gr}(-1)[-1])}.
\]
and it remains to apply $Rf_*$ and $\underline{Hom}(-,\bQ).$
\subsection*{The motive of trivial log structure.}
Let $(X,\cO^*)$ be a log scheme with trivial log structure. Then 
$$[X]^{log}\simeq \mathrm [X]$$
where $[X]$ is the homological motive of the underlying scheme $X$. Indeed, it suffices to check that $Rf_*\mathbb Q^{\mathrm{log}}_X\simeq \underline{Hom}([X],\bQ)$. Let us note that ${Rf_*\mathbb Q^{\mathrm{log}}_X\simeq Rf_*\mathbb Q_{X}}$ and use Yoneda Lemma. For any smooth $Y$ over $\mathrm k$ we have
$$\mathrm{Hom_{DM_{gm}}(\mathit{ [Y](-q)[-p]},\underline{Hom}([X],\bQ))\simeq Hom_{DM_{gm}}([\mathit Y] \!\otimes\! [\mathit X],\mathbb Q}(q)[p])=H^{p,q}(X \!\times\!_{\rm k}  Y,\mathbb Q).$$
On the other hand, 
$$\mathrm{Hom_{DM_{gm}}(\!\mathit{[Y](\!-q \!)[\!-m ], \! Rf_*\mathbb Q_{X})} \!\simeq\! Hom_{DA_{et}(\mathit{X\times_{\rm k} Y})}([\mathit Y \times_{\rm k}  \mathit X],\mathbb Q \mathit{(q)[m]})}.$$ 
\subsection*{The motive of a fs log scheme} Let $X$ be a fs log scheme. Recall that in Section 3 we constructed  virtualization functor
$$v:\mathrm{Sch_{fs}^{Log}\!/k\longrightarrow Sch_{fs}^{vLog}\!/k}$$ 
which maps a log scheme $(X,\mathcal M\to \cO)$ to the pair $(X,\,\cO^*\hookrightarrow \cM).$
We define the Kato-Nakayama motive of $X$ as the motive of the virtualization 
$$[X]^{log}:=[v(X)]^{log}.$$
\subsection*{Functoriality}
Note that we construct the functor 
\begin{equation}\label{'F6.2}
	Sch^{vLog}_{fs}/k\ra Ho(DM(k,\bQ))	
\end{equation}
to the triangulated category of Voevodsky motives. On the other hand, $DM(k,\bQ)$ has the natural DG-enrichment \cite{BV06} and, consequently, can be considered as a stable $(\infty,1)-category$. 
\begin{prop}\label{'P6.1} The functor \ref{'F6.2} gives rise to the $\infty$-functor 
	\begin{equation}\label{'key}
		[-]^{log}:\rm Sch^{vLog}/k\ra DM(k,\bQ).
	\end{equation}
\end{prop}
It suffice to show the dual functor $(X,\cM)\mapsto Rf_*\bQ^{log}_X$ can be lifted on the level of $\infty$-category. To prove the Proposition let us reformulate the construction of (\ref{'F6.2}). 
Let us consider the 1-functors\footnote{from the category of schemes to the category of all 1-categories} 
$$\rm Ho(DA_{et,\bQ}):\mathit X \mapsto Ho(DA_{et}(\mathit X,\bQ), \ \ \it f\mapsto f^*;$$
$$\rm Ho(DA_{et,\bQ})^{op}:\mathit X \mapsto Ho(DA_{et}(\mathit X,\bQ)^{op}, \ \ \it  f\mapsto Rf_*.$$
Here $\rm Ho$ means the homotopy category. %Using the adjunction we get 
%$$\int Ho \iso \int Ho.$$
Combining Corollary \ref{'C5.3} and Remark \ref{'R5.4} we get the natural transformation of $\infty$-functors $\Psi:\mathbf{vLog}^{\mathrm{fs}}\longrightarrow\mathrm{DA_{et,\mathbb{Q}}}$
which is defined as the composition
\[
\xymatrix{\mathbf{vLog}^{\mathrm{fs}}\ar[rr]^{\!\!\!\!\!\!\!\!\! \mathrm{Ex^{log}\otimes}\mathbb{Q}} &  & \mathrm{Sh_{et}(Sm/}-,\mathbb{Q})\ar@{^{(}->}[r]^{\heartsuit} & \mathrm{D}(\mathrm{Sh_{et}(Sm/}-,\mathbb{Q}))\ar[r]^{\ \ \ \ \ \ \  \mathrm{\mathrm{\Sigma_{\mathbb{L}}^{\infty}}}\cdot\mathrm{L_{\mathbb{A}^{1}}}} & \mathrm{DA_{et,\bQ}}.}
\] 
So we have the natural transformation $\mathrm{Ho}(\Psi):\mathbf{vLog}^{\mathrm{fs}}\ra\mathrm{Ho(DA_{et,\mathbb{Q}}})$
Let us define the natural transformation 
$$R\Gamma^{mot}: \rm Ho(DA_{et,\bQ})^{op}\ra \underline{DA_{et}(k,\bQ)};$$
$$(F\in DA_{et}(X,\bQ))\mapsto Rf_*F.$$
Here $\rm \underline{DA_{et}(k,\bQ)}$ is a constant sheaf and $f$ is a canonical morphism $X\to Spec(k).$

Now, we have the commutative diagram 
\[
\xymatrix{(\int_{\mathrm{Sch/k^{op}}}\mathrm{Ho(DA_{et,\mathbb{Q}}}))^{\mathrm{op}}\ar[r]^{ \ \thicksim} & \int_{\mathrm{Sch/k}}\mathrm{Ho(DA_{et,\mathbb{Q}})^{op}}\ar[r]^{\,\,\,\,\,\,\,\, \int R\Gamma^{mot}} & \mathrm{Sch/k\times \rm DA_{et}(k,\bQ)}\ar[d]^{\rm pr_{2}}\\
	\mathrm{Sch_{fs}^{vLog}/k}\ar[u]_{\int\mathrm{Sym}^{*}\mathrm{Ho}(\Psi)}\ar[rr]^{\mathrm {the \,\, functor \,\, (6.2)}} &  & \mathrm{DA_{et}(k,\bQ)}
}
\]
To construct $[-]^{log}$ it remains to use the $\infty$-category version of Grothendieck construction
\\ \textit{Proof of Proposition \ref{'P6.1}:}  Let us consider the $\infty-$functors \footnote{here $\mathbf{ Pr}^*$ mean the category of presentable categories together with the morphisms given by continuous and cocontinuous functors}
\[
\mathrm{DA_{et,\mathbb{Q}}}\!:\!\mathrm{Sch/k^{op}}\longrightarrow\mathbf{Pr}^{\mathrm{R}};\,\,\,\,\,\,\,\,X\mapsto\mathrm{DA_{et}}(X,\mathbb{Q}),\,f\mapsto f^{*};
\]
\[
\mathrm{DA_{et,\mathbb{Q}}^{op}}\!:\!\mathrm{Sch/k}\longrightarrow\mathbf{Pr}^{\mathrm{L}};\,\,\,\,\,X\mapsto\mathrm{DA_{et}}(X,\mathbb{Q}\mathrm{)^{op}},\,f\mapsto Rf_{*}.
\]
By construction, $R\Gamma^{mot}$ can be lifted to the natural transformation of $\infty$-functors $R\Gamma^{mot}:\rm DA_{et,\bQ}^{op}\ra \underline{DA_{et}(k,\bQ)}.$ 
Then we define $[-]^{log}$ as the composition
\begin{equation}\label{'diag_funct}
	\xymatrix{(\int_{\mathrm{Sch/k^{op}}}\mathrm{DA_{et,\mathbb{Q}}})^{\mathrm{op}}\ar[r]^{\ \ \varphi} & \int_{\mathrm{Sch/k}}\mathrm{DA_{et,\mathbb{Q}}^{op}}\ar[r]^{\!\!\!\!\!\!\!\!\!\!\! \int R\Gamma^{mot}} & \mathrm{Sch/k\times \rm DA_{et}(k,\bQ)}\ar[d]^{\rm pr_{2}}\\
		\mathrm{Sch_{fs}^{vLog}/k}\ar[u]_{\int\mathrm{Sym}^{*}\Psi}\ar[rr]^{[-]^{log}} &  & \mathrm{DA_{et}(k,\bQ)}
	}	
\end{equation}

Here $\int$ mean the $(\infty,1)$-Grothendieck construction (see Appendix \ref{B}) and $\varphi$ is the equivalence of categories (see Proposition \ref{PB.4}). $\square.$
\subsection*{Constructability of $[X]^{log}$}
Let $\underline X$ be a scheme. Recall that a motivic sheaf on $\underline X$ is called constructible if it can be presented as a finite colimits of $Y_i(n_i)[m_i]$ for some $Y_i\in \mathrm{Sm/}\underline X$ and $n_i,m_i\in \mathbb Z$. For any $\underline X$ constructible motivic sheaves form full subcategory of $\mathrm{DA_{et}( \underline X,\mathbb Q)}$ which is stable under six operations.
\begin{thm}\label{T5.6}
	%Suppose that $X$ is irreducible, geometrically unibranch.  
	For any fs virtual log scheme $(X,\cM)$ the motivic sheaf $\mathbb Q^{\mathrm{log}}_X$ is constructible.
\end{thm}
\proof Firstly, suppose that $X$ is normal. The category of constructible motivic sheaves closed under tensor product and direct summands. So it suffice to check that $\mathrm{ M^{gr}}$ is constructible. The motivic sheaf $\mathrm{M^{gr}}$ is the extension of $\mathbb{Q}(1)[1]$ by $\mathrm{\overline{M}^{gr}=\Sigma^\infty_{\mathbb G_m}L_{\mathbb A^1}Ex_X(\overline{\mathcal M}^{gr})}$ (see Section 6). So it suffice to check that $\mathrm{\overline{M}^{gr}}$ is constructible. The six operations preserve constructible objects. Hence we may assume that $\mathcal{\overline{M}}^{gr}$ is a local system. For normal schemes the fiber functor 
$$\mathrm{Fib}_x: LocSys(\mathit X)_{\mathbb{Q}}\sa{\sim} Rep_{cont}(\pi_1^{et}(\mathit{X,x)})_{\mathbb{Q}}$$
is an equivalence of categories between finite type local systems and finite continuous representations. 
%[Stack??? \url{https://stacks.math.columbia.edu/tag/0DV4}]
So the action of $\pi_{1}^{et}(X,\overline{x})$ on $\overline{\mathcal{M}}^{gr}_x$ is continuous and the stabilizer
of each element is an open subgroup. Hence the kernel $\pi_{1}^{et}(X,\overline{x})\longrightarrow\mathrm{Aut}(\overline{\mathcal{M}}^{gr}_x)$
is open and the image is a finite group. Consequently, there is the Galois cover $f: Y\rightarrow X$ such that the action of $\pi_1^{et}(X,x)$ on the stack $\overline{\mathcal{M}}^{gr}_x$ factors through $Aut(Y/X)$.

Now, using the equivalence of categories  we conclude that  $\overline{\mathcal{M}}^{gr}=\bigoplus_i \mathrm{Fib_x^{-1}}(\rho_i) $ where $\rho_i$ are some irreducible representations of $Aut(Y/X)$. So $\mathrm{\overline M^{gr}}=\bigoplus_i [\rho_i]$ where $[\rho_i]:=\rm \Sigma^\infty_{\mathbb G_m}L_{\mathbb A^1}Ex_X(Fib_x^{-1}(\rho_i)).$ But for any irreducible representation $\rho_i$ the corresponding motivic sheaf $[\rho_i]$ is constructible. Indeed, for the representable local system  $\mathbb{Q}[Y]$ we have 
$$\mathrm{Fib_x} \bQ[Y]=\bQ[Aut(Y/X)]\iso\bigoplus_j \rho_j^{\rm dim(\rho_j)}$$
there $j$ runs of all irreducible representations of $Aut(Y/X)$. On the other hand, the motivic sheaf $[Y]$ is constructible.

For arbitrary $X$ let us use induction by $\dim(X).$ Let $\dim(X)=0$. By results of \cite{'key-22}, any motivic sheaf admits h-descent. So $DA_{et}(X,\bQ)\iso DA_{et}(X_{red},\bQ)$ and we may assume that $X$ is reduced.
Then $X=\sqcup_mSpec(L_m)$ where $L_m/k$ are algebraic extensions. So ${\bQ^{log}_X=\oplus_m i_m*\bQ^{log}_{Spec(L_m)}.}$
%Spec(L_m) is normal$

Suppose we check the statement for all $Y$ with $\dim(Y)<n.$ Let $X^\nu$ be a normalization of $X$. The canonical morphism $\varphi:X^\nu\to X$ is birational. So we have an abstract blow-up square 
$$\begin{tikzcd}
	Z' \arrow[d, "\phi"] \arrow[r, "i'"] & X^\nu \arrow[d, "\varphi"] \\
	Z \arrow[r, "i"]                     & X                         
\end{tikzcd}.$$
Note that $\bQ^{log}_X$ is a motivic sheaf so it admits h-descent. Moreover, the functor $\cM_{X}\mapsto \bQ^{log}_X$ commutes with inverse images. Hence we get the exact triangle 
$$R(i\phi)_*\bQ^{log}_{Z'}\to i_*\bQ^{log}_{Z}\oplus R\varphi_*\bQ^{log}_{X^\nu}\to \bQ^{log}_{X}$$ where $\dim(Z')$ and $\dim(Z)$ less than n. $\square$

Note that for $X=Spec(k)$ the subcategory of constructible motivic sheaves is equivalent to the category of geometrical Voevodsky motives $DM_{gm}(k,\bQ).$ As a consequence we get 
\begin{cor}
	For any log scheme $(X,\cM)$ the motive $[X]^{log}$ belongs to $DM_{gm}(k,\bQ).$
\end{cor} 
\proof The six operations preserve constructability 
%so for normal $X$ the statement follows from Theorem \ref{T5.6}. For arbitrary $X$ let us use induction by $\dim(X).$ If $\dim(X)=0$ we may assume that $X=Spec(L)$ for some algebraic field extension $L/k$ so $X$ is normal. Suppose we check the statement for all $Y$ with $\dim(Y)<n.$ Let $X^\nu$ be a normalization. Then the canonical morphism $X^\nu\to X$ is birational. Hence, we get the log abstract blow-up sequence $$[Z']^{log}\to [Z]^{log}\oplus [X^\nu]^{log}\to [X]$$ where $\dim(Z')$ and $\dim(Z)$ less than n
. $\square$
\subsection*{The structure of $\mathbb E_\infty$-coalgebra on $[X]^{log}$}

%\subsection*{Monoidality of $\bQ^{log}_X.$}

\begin{comment}

Let $X$ be a fs log scheme. Then the inclusion $\mathcal O^*\hookrightarrow \mathcal M^{\mathrm{gr}}$ induces the map  $\bQ\to\mathrm{\Omega_{\mathbb G_m}M^{gr}}$. So $\mathrm{\Omega_{\mathbb G_m}M^{gr}}$ can be considered as an object of $\rm DA_{et}(\mathit{\underline{X}},\bQ)_{\bQ/}$. By Theorem \ref{T2.6} and Lemma \ref{L2.7} we have
$$\mathrm{Sym^*(\mathrm{\Omega_{\mathbb G_m}M^{gr}})\simeq colim_kS^k(\mathrm{\Omega_{\mathbb G_m}M^{gr}})}.$$
Using Lemma \ref{L5.5} and the isomorphisms \ref{F5.3} we conclude that $\mathrm{colim_kS^k(\mathrm{\Omega_{\mathbb G_m}M^{gr}})\simeq S^n(\mathrm{\Omega_{\mathbb G_m}M^{gr}})}$ where $n$ is the maximal rank of $\mathcal{ \overline M}^{\mathrm{gr}}_{X,x}$. So 
\begin{equation}\label{F5.4}
\mathbb Q^{\mathrm{log}}_X\simeq \mathrm{Sym^*(\mathrm{\Omega_{\mathbb G_m}M^{gr}})}
\end{equation} 
and $\mathbb Q^{\mathrm{log}}_X$ has the structure of free $E_\infty$-algebra. 
\begin{rem}
Using the isomorphism \ref{F5.4} one can equipped the groups $H^{p,q}(X,\mathbb Q)$ with the structure of bigraded supercommutative algebra with the multiplication induces by the maps 
$$\mathbb Q^{\mathrm{log}}_X(q)\otimes\mathbb Q^{\mathrm{log}}_X(q')\simeq(\mathbb Q^{\mathrm{log}}_X\otimes\mathbb Q^{\mathrm{log}}_X)(q+q')\overset{m\otimes\rm id}{\longrightarrow}\mathbb Q^{\mathrm{log}}_X(q+q').$$
\end{rem}
\end{comment}
Note that the construction of $\bQ^{log}_X$ is functorial so we have the functor
\begin{equation}\label{F5.5}
	\begin{tikzcd}
		\mathrm{vLog/}\underline X \arrow[rr, "\mathcal M^{\mathrm{gr}}\mapsto \Omega_{\mathbb G_m}\mathrm{M^{gr}}"] &  & {\mathrm{DA_{et}(}\underline X,\mathbb Q)_{\mathbb Q/}} \arrow[rr, "\mathrm{Sym^*}"] &  & {\mathrm{CAlg(DA_{et}}(\underline X,\mathbb Q)).}
	\end{tikzcd}
\end{equation}
Now, let us consider $\mathrm{vLog}/\underline X$ as the monoidal category together with cocartesian monoidal structure and $\mathrm{CAlg(DA_{et}}(\underline X,\mathbb Q))$ as the monoidal category relatively to the tensor product of algebras.
\begin{prop}\label{P5.8}
	The functor (\ref{F5.5})  is monoidal.
\end{prop}
\proof Note that %by Prop.3.7. nlab "over-(infinity,1)-category"
the coproduct of two objects $(\mathbb Q\rightarrow E)$, $(\mathbb Q\rightarrow F)$ of $\mathrm{DA_{et}(}\underline X,\mathbb Q)_{\mathbb Q/}$ coincide with the colimit of the diagram 
$$\begin{tikzcd}
	\mathbb Q \arrow[d] \arrow[r] & E \\
	F                             &  
\end{tikzcd}$$
in $\mathrm{DA_{et}(}\underline X,\mathbb Q)$. So let us equip $\mathrm{DA_{et}(}\underline X,\mathbb Q)_{\mathbb Q/}$ with cocartesian monoidal structure. Then the first functor in the composition (\ref{F5.5}) is monoidal. 
On the other hand, by Proposition 3.2.4.7. of \cite{key-1} the monoidal structure on $\mathrm{CAlg(DA_{et}(}\underline X,\mathbb Q))$ induced by $\otimes$ is cocartesian and $\mathrm{Sym^*}$ commutes with colimits. So the second functor in the composition (\ref{F5.5}) is also monoidal. $\square$

%\subsection*{Monoidality of $[-]^{log}$.} 

Let us endow the category of virtual log schemes $\rm Sch^{vLog}_{fs}/k$ with Cartesian monoidal structure.
\begin{thm}\label{T7.5}
	The $\infty$-functor $[-]^{log}:\rm Sch^{vLog}_{fs}/k \ra DM_{gm}(k,\bQ)$ is monoidal. 
\end{thm}
\proof By Remark \ref{MonAdj} for any morphism of schemes $f:X\ra Y$ we have the adjoint pair 
$\rm \mathit{Rf_*}: CAlg(DA_c(X,\bQ))\leftrightarrows CAlg(DA_c(Y,\bQ)):\mathit{f^*}$. Note that the anti-equivalence 
$$\rm Hom(-,\bQ): DM_{gm}(k,\bQ)\to DM_{gm}(k,\bQ)$$
maps algebras to coalgebras. So we can define the new contravariant functor 
\begin{equation}\label{F7.6}
	\rm Sch^{vLog}/k\ra coCAlg(DM_{gm}(k,\bQ))
\end{equation}  
replaced in the diagram (\ref{'diag_funct}) $\rm DA_{et}$ with $\rm CAlg(DA_{et})$.  
\\By Example 3.2.4.5 of \cite{key-1} the motivic tensor product $\otimes$ induces on $\rm coCAlg(DM_{gm}(k,\bQ))$ the structure of symmetric monoidal category. The functor $[-]^{log}$ can be recovered as the composition of \ref{F7.6} with the forgetful functor $\rm coCAlg(DM_{gm}(k,\bQ))\ra DM_{gm}(k,\bQ)).$ So it suffice to check that \ref{F7.6} is monoidal.

Note that the monoidal structure on $\rm coCAlg(DM_{gm}(k,\bQ))=Alg(DM_{gm}(k,\bQ)^{op})^{op}$ is Cartesian by Proposition 3.2.4.7. of \cite{key-1}. So it remains to show that \ref{F7.6} commutes with products or, equivalently, check that the functor 
$$X\mapsto ([X]^{log})^*$$
maps a product of log schemes to the tensor product of algebras. Let $X,Y$ be virtual log schemes. 
Observe that the functors 
\[
\begin{tikzcd}
	\mathrm{vLog/}\underline X \arrow[rr, "\mathcal M^{\mathrm{gr}}\mapsto \Omega_{\mathbb G_m}\mathrm{M^{gr}}"] &  & {\mathrm{DA_{et}(}\underline X,\mathbb Q)_{\mathbb Q/}} \arrow[rr, "\mathrm{Sym^*}"] &  & {\mathrm{CAlg(DA_{et}}(\underline X,\mathbb Q))}
\end{tikzcd}
\]
is monoidal for each $X$ (Proposition \ref{P5.8}) and commute with inverse image functors. 
%By monaidality of the functor of ? $vLog_X\ra CAlg(DA(X,\bQ))$
So we have the isomorphism of algebras 
$$\rm\bQ^{log}_{X\times Y}\iso \bQ^{log}_X\boxtimes \bQ^{log}_Y.$$
On the other hand, Cisinski proved the Künneth formula for motivic sheaves (see \cite{key-11}). So we have
\begin{equation}\label{F7.7}
	\rm Rf_*(\bQ^{log}_X\boxtimes \bQ^{log}_Y)\iso R\alpha_{X*}\bQ^{log}_X\otimes R\alpha_{Y*}\bQ^{log}_Y 
\end{equation}
where $f,\alpha_X,\alpha_Y$ are canonical morphisms to $\rm Spec(k).$ It suffice to show that \ref{F7.7} is an isomorphism of $\bE_\infty$-algebras. 
\begin{lem}\label{L7.5}
	Let $A$ and $B$ be $\bE_\infty$-algebras in $\mathrm{DA_{et}}(X\times Y,\bQ).$ Then the canonical map 
	$$\psi:Rf_*A\otimes Rf_*B\ra Rf_*(A\otimes B)$$
	can be extended to the homomorphism of $\mathbb E_\infty$-algebras.
\end{lem}
\proof Note that $\psi$ is adjoint to $\rm\varepsilon_A\otimes\varepsilon_B\in Hom(\mathit{f^*Rf_*A\otimes f^*Rf_*B,A\otimes B})$ and apply the Remark \ref{MonAdj}.$\square$
\\Using Remark \ref{MonAdj} and the Lemma we can define the homomorphism of algebras $\mu$ as the composition:
\begin{equation}\label{F7.8}
	\begin{tikzcd}
		\rm R\alpha_{X*}(R\pi_{X*}\pi^*_X\mathbb Q^{log}_X)\otimes R\alpha_{Y*}( R\pi_{Y*}\pi^*_Y\mathbb Q^{log}_Y)  \arrow[r, "\sim"]                            & \rm Rf_*\pi^*_X\mathbb Q^{log}_X\otimes Rf_*\pi^*_Y\mathbb Q^{log}_Y \arrow[d, "\psi"] \\
		\rm R\alpha_{X*}\mathbb Q^{log}_X\otimes R\alpha_{Y*}\mathbb Q^{log}_Y \arrow[u, "\rm R\alpha_{X*}(\eta)\otimes R\alpha_{Y*}(\eta)"] \arrow[r, "\mu"] & \rm Rf_*(\mathbb Q^{log}_X\boxtimes\mathbb Q^{log}_Y).                                 
	\end{tikzcd}
\end{equation}
Finally, comparing the diagram \ref{F7.8} with the construction of the isomorphism \ref{F7.7} (see \cite{key-11}) we conclude that $\mu$ is an isomorphism of algebras. $\square$

\section{Betti realization of $[X]^{log}$}
Let k be a field of characteristic
zero together with the embedding $\mathrm{\sigma : k \hookrightarrow \mathbb C}$. Let us fix the ring of coefficients R. In the paper \cite{Ayo10}, Ayoub constructed for any $X$ the Betti realization of motivic sheaves
$$\mathrm{Betti_X}: \mathrm{DA_{et}(}X,\mathrm R)\rightarrow \mathrm{D(Sh(} X^{\mathrm{an}},\mathrm R))$$
%such that the collection of functors $\{ \mathrm{Betti_X} \}$ commutes with six operations. 
Let us fix $X$. The construction of $\mathrm{Betti_X}$ is the following. Let us consider the diagram 
\begin{equation}\label{Betti_diag}
	\begin{tikzcd}
		\mathrm{D(Sh_{et}(Sm/} X,\mathrm{R})) \arrow[rr, "\mathrm{An^*}"] \arrow[d, "\mathrm L_{\mathbb{A}^1}"]      &  & \mathrm{D(Sh_{top}(Sm/} X^{\mathrm{an}},\mathrm{R})) \arrow[d, "\mathrm L_{\mathbb{D}}", shift left=4] \arrow[rrd, "\mathrm{Res}_X^{an}"]                                            &  &                                                                                    \\
		\mathrm{DA_{et}^{eff}(} X,\mathrm{R}) \arrow[d, "\Sigma_{\mathbb{G}_m}^{\infty} "] \arrow[rr, "\mathrm{An^*}"] &  & \mathrm{DA_{an}^{eff}(} X,\mathrm{R}) \arrow[d, "\Sigma_{\mathring{\mathbb D}}^{\infty} "', shift right=1] \arrow[u, "i", shift left=1] \arrow[rr, "\mathrm{Res}_X^{an} \cdot i"] &  & \mathrm{D(Sh(} X^{\mathrm{an}},\mathrm{R})) \arrow[llu, "\mathrm{Ex}_X^{an}"', bend right, shift right=0] \\
		\mathrm{DA_{et}(} X,\mathrm{R}) \arrow[rr, "\mathrm{An^*}"]                                                    &  & \mathrm{DA_{an}(} X,\mathrm{R}) \arrow[u, "\Omega_{\mathring{\mathbb D}}^{\infty} "', shift right=4]                                                                     &  &                                                                                   
	\end{tikzcd}
\end{equation}
Here we denote by $\mathrm{An^*}$ the Kan extensions of the analytification functor
$$\mathrm{Sm}/X\to \mathrm{Sm}/X^{\rm an},\,\, X\mapsto (X\times_{\mathrm k} \mathbb C)^{an}.$$
By results of \cite{Ayo10} the functors $\Sigma_{\mathring{\mathbb D}}^{\infty}$ and  $\mathrm{Res}_X^{an} \cdot i$ are equivalences of categories. So we can define $\mathrm{Betti}_X$ as the  composition $\mathrm{Res}_X^{an} \cdot i \cdot \Omega_{\mathring{\mathbb D}}^{\infty} \cdot \mathrm{An^*}$. 
\\Let us denote by $\mathrm{DA_{et,c}(X,R)}$  the subcategory of compact objects. The main result of \cite{Ayo10} is the following.
\begin{thm}
	(Theorem 3.19 of \cite{Ayo10}). The restriction of $\{ \mathrm{Betti}_X \}$ on the subcategory of constructible objects
	$$\mathrm{Betti}:\mathrm{DA_{et,c}}\rightarrow \mathrm{D(Sh(-^{an},R)}$$
	preserves six operations. 
\end{thm} 
\begin{rem}
	The original theorem consists the categories $\mathbf{SH}_{\mathfrak M,c}(X)$ rather than $\mathrm{DA_{et,c}}(X,\mathrm R)$. By construction  $\mathbf{SH}_{\mathfrak M,c}(X)=\mathrm{lim}_{\Omega_{\mathbb{G}_{m}}}\mathrm{Ho_{\mathbb A^1, Nis}(PSh(Sm/X,\mathfrak M)}$. One can put $\mathfrak M=\mathrm{R-Mod}$. On the other hand, the proof of the theorem in the étale case tautologically coincides with the case of Nisnevich topology. 
	%Note that the objects $ Y(n)\otimes R, Y\in\mathrm{Sm}/X, \, n\in \mathbb Z $ form a generating family of compact generators of $\mathrm{DA_{et}(}X,\mathrm R)$ by Theorem 5.2.4 [EtaleMotivesCisinsky]????.
\end{rem}
\subsection*{The Betti realization of $\bQ^{log}_X$ and $[X]^{log}$.} 
Let $\mathrm R= \mathbb Z$. For a virtual log structure ${\alpha :\mathcal{O}^* \rightarrow \mathcal{M}^{\mathrm{gr}}}$  on $X$ we will denote by $\mathrm{M^{gr}_\mathbb Z}$ the motivic sheaf $\Sigma^{\infty}_\mathbb L \mathrm{L_{\mathbb A^1}Ex^{log}_X(\mathcal M^{gr})}$. Note that $ \mathrm{M^{gr}=M^{gr}_\mathbb Z\otimes \mathbb Q}$.

Let $\alpha :\mathcal{O}^*_{an} \rightarrow \mathcal{M}^{\mathrm{gr}}_{an}$ be the analytification. Recall that we denoted by $\mathrm{exp}(\alpha)$ the complex 
\[
...\longrightarrow\underset{(-1)}{0_{\,_{\,}}}\longrightarrow\underset{(0)}{\mathcal{O}_{an_{\,}}}\longrightarrow\underset{(1)}{\mathcal{M}_{an}^{gr}}\longrightarrow\underset{(2)}{0_{\,_{\,}}}\longrightarrow...
\]
where the differential given by the composition $\mathcal{O}_{an}\!\overset{\mathrm{exp}}{\rightarrow}\mathcal{O}_{an}^{*}\!\hookrightarrow\mathcal{M}_{an}^{gr}$.
\begin{thm}\label{T6.3}
	There is the natural quasi-isomorphism $\mathrm{Betti}_X(\mathrm{M^{gr}_\mathbb Z})[-1]\simeq \mathrm{exp}(\alpha )$.
\end{thm}
Let $$ \mathrm{An^*: Sh_{et}(Sm/}X)\longrightarrow \mathrm{Sh_{top}(Sm/}X^{\mathrm{an}})$$ be the Kan extension of the analytification functor.  For any $ W\in \mathrm{Sm/}X^{\mathrm{an}} $ we have the map 
\begin{equation}\label{F6.1}
	\underset{W \overset{g}{\longrightarrow} Y^{\mathrm{an}}}{colim} \mathcal{O}^*(Y)\longrightarrow \mathcal{O}^*_{\mathrm{an}};
\end{equation}
$$ (f\in \mathcal O^*(Y)) \mapsto g^*f. $$
We define the map
\begin{equation}\label{F6.2}
	\psi: \mathrm{An^*}(\mathcal{O}^*_{\mathrm{Sm/}X})\longrightarrow \mathcal{O}^*_{\mathrm{Sm/}X^{\mathrm{an}}}
\end{equation} 
as the sheafification of \ref{F6.1}. 
\begin{lem}
	The map $\psi$ is an isomorphism.
\end{lem}
\proof  The analytification functor commutes with fiber product. %(see SGA) 
Consequently, $\mathrm{An^*}$ is exact. %(see "morphisms of sites"?)
In Appendix \ref{A}, for any abelian algebraic group $\mathcal G$, we defined the complex of sheaves $\boldsymbol Q_{\geq{-1}}(\mathcal G)$ with $^0(\boldsymbol Q_{\geq{-1}}(\mathcal G))=\mathcal G$. Note that $\psi$ can be extended to the map of complexes 
\begin{equation}\label{key}
	\boldsymbol Q_{\geq{-1}}(\mathbb G_m)\longrightarrow \boldsymbol Q_{\geq{-1}}(\mathbb G_m^{\mathrm{an}})	
\end{equation}
which induced by the map 
\begin{equation}\label{key}
	\underset{W \overset{g}{\longrightarrow} Y^{\mathrm{an}}}{colim} \mathbb Z[\mathbb G_m](Y) \overset{\psi^{'}}{\longrightarrow} \mathbb Z[\mathbb G_m^{\mathrm{an}}](W);
\end{equation}
$$ \Sigma_i a_i[f_i] \mapsto \Sigma_i a_i[g^*f_i]. $$
Note that $\mathrm{An^*}$ is exact and consequently, commutes with cohomology. So it suffice to check that $\psi^{'}$ is an isomorphism. But $\psi^{'}$ is the canonical isomorphism between $ \mathrm{An^*}(h_{\mathbb{G}_m})  $ and $ h_{\mathbb{G}_m^{\mathrm{an}}} $ which can be define by the same way for any Kan extension and any presentable object. $ \square $

Let us consider the small analytic and the smooth analytic sites of $X^{\mathrm{an}}$. As in the case of the small and the smooth étale sites there exists the pair of adjoint functors 
$$\mathrm{Ex_\mathit{X}^{an}:Sh(\mathit{X}^{an},\mathrm R)\rightleftarrows Sh_{top}(Sm/\mathit X^{an},\mathrm R):Res_\mathit{X}^{an}}. $$
So we can define the functor 
$$ \mathrm{Ex^{log}_{an}:vLog/\mathit X \longrightarrow Sh_{top}(Sm/X^{an},\mathbb Z)} $$ in the same way as in Section 6. Namely, we define $ \mathrm{Ex^{log}_{an}\mathcal M^{gr}_{an}} $ as the colimit of the diagram
$$ \begin{tikzcd}
	\mathrm{Ex_{\mathit X^{an}}\mathcal O^*_{an}} \arrow[d] \arrow[r] & \mathrm{Ex_{\mathit X^{an}}\mathcal M^{gr}_{an}} \\
	\mathrm{\mathcal O^*_{Sm/\mathit X^{an}}}                         &   \,\,\,\,\,\,\,\,\,\,\,\,\,\,\,\,\,\,\,\,\,\,\,\,\,\,\,\,\,\,\,\,\,\,\,.                                              
\end{tikzcd} $$
Now, let us prove that the map \ref{F6.2} induces the canonical isomorphism 
$$ \phi : \mathrm{An^*(Ex^{log}\mathcal M^{gr})\longrightarrow Ex^{log}_{an}\mathcal M^{gr}_{an}}. $$
Indeed, by construction of $\mathrm{Ex^{log}}$ there is the pushout 
$$\begin{tikzcd}
	\mathrm{ \mathrm{An^*} Ex_{\mathit X}(\mathcal O^*)} \arrow[r] \arrow[d] & \mathrm{ \mathrm{An^*} Ex_{\mathit X}(\mathcal M^{gr})} \arrow[d] \\
	\mathrm{\mathrm{An^*}(\mathcal O^*_{Sm/X})} \arrow[r]                    & \mathrm{\mathrm{An^*}(Ex^{log}\mathcal M^{gr})}                  
\end{tikzcd}$$
On the other hand, let us consider the diagram 
$$ \begin{tikzcd}
	\mathrm{Ex_{\mathit X^{an}} \mathrm{An^*} (\mathcal O^*)} \arrow[r] \arrow[d] & \mathrm{Ex_{\mathit X^{an}} \mathrm{An^*} (\mathcal M^{gr})} \arrow[d] \\
	\mathrm{Ex_{\mathit X^{an}}\mathcal O^*_{an}} \arrow[d] \arrow[r]           & \mathrm{Ex_{\mathit X^{an}}\mathcal M^{gr}_{an}} \arrow[d]           \\
	\mathrm{\mathcal O^*_{Sm/\mathit X^{an}}} \arrow[r]                         & \mathrm{Ex^{log}_{an} \mathcal M^{gr}_{an}}                         
\end{tikzcd} $$
Here the top square is cocartesian by construction of $\mathrm{\mathcal M^{gr}_{an}}$ and the bottom square is cocartesian by construction of $\mathrm{Ex^{log}_{an}}$. So the large square is also cocartesian and the existence of $\phi$ follows from the Lemma.
\begin{lem}
	The morphisms
	\begin{equation}\label{F6.5}
		\mathrm{An^*Ex}_X(\mathcal{O}^*_{small})\rightarrow \mathrm{An^*}(\mathcal{O}^*_{Sm/X})
	\end{equation}
	and 
	\begin{equation}\label{F6.6}
		\mathrm{Ex}_X^{an}\mathrm{An^*}(\mathcal{O}^*_{small})\rightarrow \mathrm{Ex}_X^{an} (\mathcal{O}^*_{an.small})\rightarrow \mathcal{O}^*_{Sm/X^{an}}
	\end{equation}
	are the same up to the canonical isomorphism \ref{F6.2}.
	
\end{lem}

\proof Note that $\mathrm{Ex}_X$, $\mathrm{Ex}_{X^{an}}$ and $\mathrm{An^*}$ commute with sheafification. So it suffice to check the statement on the level of preshaves.

Let $W$ belongs to $\mathrm{Sm}/X^{\rm an}$. Then the map \ref{F6.5} induces by the composition
$$\underset {W\rightarrow U\times_XW}{colim}(\underset{U\rightarrow V^{am}}{colim} \mathcal{O}^*(V\sa{et}X))\ra \underset {W\rightarrow U\times_XW}{colim} \mathcal{O}^*_{an}(U)\ra \mathcal{O}^*_{an}(W) $$ 
where the both maps induced by pullbacks of functions.
Note that $W$ is smooth over $X$. So we can drop the external colimit and assume that $U=im(W)$:
$$\underset{im(W)\rightarrow V^{am}}{colim} \mathcal{O}^*(V\sa{et}X)\ra \mathcal{O}^*_{an}(im(W))\ra \mathcal{O}^*_{an}(W). $$
On the other hand, the arrow \ref{F6.6} is the sheafification of the map
$$\underset{W\rightarrow Y^{an}}{colim}(colim \, \mathcal{O}^*(V\sa{et}X))\ra
\underset{W\rightarrow Y^{an}}{colim} \mathcal{O}^*(Y\sa{sm}X))$$
where the internal colimit taken by the collection of the diagrams 
\begin{center}
	\begin{tikzcd}
		Y \arrow[rd, "id"] \arrow[r] & V\times_XY \arrow[d] \\
		& Y                   
	\end{tikzcd}
	.
\end{center}
Note that this collection coincides with
\begin{center}
	\begin{tikzcd}
		Y \arrow[rd] \arrow[r] & V \arrow[d, "et \, \, \, \, ."] \\
		& X                
	\end{tikzcd}
\end{center} 
So we can rewrite the arrow \ref{F6.6} as the map
$$\underset{W\rightarrow V^{an}}{colim}\mathcal{O}^*(V\sa{et}X)\ra \underset{W\rightarrow Y^{an}}{colim}\mathcal{O}^*(Y\sa{sm}X).$$
Finally, let us consider the diagram 
\begin{center}
	\begin{tikzcd}
		{\,\,\,\,\,\,\,\, \underset{W\rightarrow V^{an}}{colim}\mathcal{O}^*(V\sa{et}X)} \arrow[d, "id", no head, Rightarrow] \arrow[rrr] &  &                                & {\underset{W\rightarrow Y^{an}}{colim}\mathcal{O}^*(Y\sa{sm}X)} \arrow[d] \\
		{\underset{im(W)\rightarrow V^{an}}{colim}\mathcal{O}^*(V\sa{et}X)} \arrow[rr]                                   &  & \mathcal{O}_{an}^*(im(W)) \arrow[r] & \mathcal{O}_{an}^*(W)                                                                                            
	\end{tikzcd}
\end{center}
where the upper horizontal arrow is \ref{F6.5}, the lower horizontal arrows is \ref{F6.6} and the right vertical arrow is \ref{F6.1}. Now, one can check that the diagram is commutative.$\square$

Let $X$ and $\alpha$ be as above. Let us define the complex $ \mathcal Exp \mathrm{(\alpha)\in Ch^{\bullet}(Sh_{top}(Sm/\mathit X^{an}})$: 
\[
...\longrightarrow\underset{(-1)}{0_{\,_{\,}}}\longrightarrow\underset{(0)}{\mathcal{O}_{\mathrm{Sm/\mathit X^{an}}}}\overset{exp}{\longrightarrow}\underset{(1)}{\mathrm{Ex^{log}_{an} \mathcal M_{an}^{gr}}}\longrightarrow\underset{(2)}{0_{\,_{\,}}}\longrightarrow...
\]
Let us define the analytic motivic sheaf $[\mathcal Exp\mathrm{(\alpha)]}$ using the functors $\mathrm{L_{\mathbb D}}$ and $\Sigma^\infty_{\mathring{\mathbb D}}$. 

Note that the functor $$ \mathrm{Res\cdot \mathit i \cdot \Omega^{\infty}_{\mathring{\mathbb D}}:DA_{an}(\mathit X,\bZ)\longrightarrow D(Sh(\mathit X^{an},\bZ))} $$ maps $[\mathcal Exp\mathrm{(\alpha)]}$ to $\rm exp(\alpha)$. Indeed, by previous there exists the canonical map $$\mathrm{Ex_{\mathit X^{an}}(\mathit{exp(\alpha)})\longrightarrow} \mathcal Exp\mathrm{(\alpha)}$$ which is a quasi-isomorphism. But the functor $\Sigma^\infty_{\mathring{\mathbb D}} \cdot \mathrm{L_{\mathbb D}} \cdot \mathrm{Ex_{X^{an}}}$ is adjoint to $\mathrm{Res\cdot \mathit i \cdot \Omega^{\infty}_{\mathring{\mathbb D}}}$ and hence inverse.
\\\textit{Proof of Theorem \ref{T6.3}:} Observe that the commutative diagram 
$$ \begin{tikzcd}
	0 \arrow[d] \arrow[r]                                               & \mathcal O_{\mathrm{Sm/\mathit X^{an}}} \arrow[d, "exp"] \\
	\mathrm{An^*(Ex^{log}\mathcal M^{gr}}) \arrow[r, "\phi"] & \mathrm{Ex^{log}_{an}\mathcal M^{gr}_{an}}             
\end{tikzcd} $$ 
define the map of complexes $$\mathrm{An^*(Ex^{log}\mathcal M^{gr})[-1] \overset{\phi^{'}}{\longrightarrow} \mathcal E\mathit{xp}(\alpha)}. $$
Note that $\mathrm{An^*}$ commutes with $\mathrm{L_{-}}$ and $\Sigma^{\infty}_{-}$. So it suffice to show that $\phi^{'}$ is an $\mathbb A^1$-equivalence or equivalently check that $$\mathrm{L_{\mathbb D}(Cone(\phi^{'}))}\simeq \mathrm{L_{\mathbb D}(\mathcal O_{Sm/\mathit X^{an}})}=0.$$
Indeed, let us consider the commutative diagram  
$$ \begin{tikzcd}
	{\mathbb Z[pt_X]\otimes \mathcal O_{\mathrm{Sm/\mathit X^{an}}}} \arrow[d, "{"1"}"] \arrow[rrd, "\mathrm{id}_{\mathcal O}", bend left]  &                                                                                                       &                                         \\
	{\mathbb Z[\mathbb A^1_X]\otimes \mathcal O_{\mathrm{Sm/\mathit X^{an}}}} \arrow[r, "pr \times id"]               & \mathcal O_{\mathrm{Sm/\mathit X^{an}}}\otimes \mathcal O_{\mathrm{Sm/\mathit X^{an}}} \arrow[r, "m"] & \mathcal O_{\mathrm{Sm/\mathit X^{an}}} \\
	{\mathbb Z[pt_X]\otimes \mathcal O_{\mathrm{Sm/\mathit X^{an}}}} \arrow[u, "{"0"}"] \arrow[rru, "0"', bend right] &                                                                                                       &                                        
\end{tikzcd}
$$ 
where $m$ is the multiplication. This diagram defines the $\mathbb A^1$-homotopy between $\mathrm{id}_{\mathcal{O}}$ and the zero endomorphism. $\square$ 

Now, let $X$ and $\alpha$ be as above, $X^{log}$ be the corresponding Kato-Nakayama space and $\pi:X^{log}\rightarrow X^{an}$  be the canonical map. 
\begin{thm}\label{T6.6}
	There is the canonical quasi-isomorphisms 
	$$\mathrm{Betti}_X(\mathbb Q^{\mathrm{log}}_X)\simeq R\pi_*\mathbb Q; \,\,\,\,\,\,\,\,\,\,\,\,\, [X]^{log}\iso Sing_*(X^{log},\bQ).$$
\end{thm}
\proof Each $\bQ^{log}_X$ is constructible. So the second quasi-isomorphism can be derived from the first using the six operations. Note that the functor $\mathrm{Betti}_X$ is monoidal. So we have
$$\mathrm{Betti_X \mathbb Q^{\mathrm{log}}_X\simeq Betti_XS^n(\Omega_{\mathbb G_m}M^{gr})\simeq S^n(Betti_X\Omega_{\mathbb G_m}M^{gr})\simeq  S^n(Betti_XM^{gr}[-1])}$$
and 
$$\mathrm{S^n(Betti_XM^{gr}[-1])\simeq  S^n(Betti_X(M^{gr}_\mathbb Z)[-1]\otimes \mathbb Q)\simeq S^n(}\mathrm{exp}(\alpha)\otimes \mathbb Q).  $$
Hence $\mathrm{Sym}^{*}(\mathrm{M^{gr}}(-1)[-1])\iso\mathrm{Sym}^*(\mathrm{exp}(\alpha)\otimes \mathbb Q)$ by Theorem \ref{T2.6}.
It remains to use Theorem \ref{R_pi_Q}. $\square$ 

\section{Logarithmic $h$-descent}

In analogy to Kato's logarithmic geometry,
we will say that a morphism of virtual log schemes $(f,f^{\#})$
is strict if $f^{\#}$ is an isomorphism.
\begin{defn}
	We will say that the morphism $f\!:\!Y\longrightarrow X$ of virtual
	log schemes is a\textit{ strict h-cover} if the underlying morphism
	$\underline{Y}\longrightarrow\underline{X}$ is a h-cover and $f$
	is strict. Let $\mathcal{F}$ be a presheaf of complexes on $\mathrm{Sch_{fs}^{vLog}\!/k}$.
	We will say that $\mathcal{F}$ \textit{admits strict h-descent} if
	for any strict h-cover $f\!:\!Y\longrightarrow X$ the canonical map
	of complexes 
	\[
	\mathcal{F}(X)\rightarrow\mathrm{Tot^{\bullet}}(\mathcal{F}(Y)\rightarrow\mathcal{F}(Y\times_{X}Y)\rightarrow\mathcal{F}(Y\times_{X}Y\times_{X}Y)\rightarrow...)
	\]
	is a quasi-isomorphism.
\end{defn}
\begin{defn} 
	A \textit{logarithmic abstract blow-up} is a Cartesian square
	\begin{equation}\label{'F4.1}
		\xymatrix{Z^{\prime}\ar[r]^{i^{\prime}}\ar[d]^{p_{Z}} & X^{\prime}\ar[d]^{p}\\
			Z\ar[r]^{i} & X
		}
	\end{equation}
	in $\mathrm{Sch_{fs}^{vLog}\!/k}$ such that the maps $i$ and $i^{\prime}$
	are strict, the map $p|_{X^{\prime}\setminus Z^{\prime}}\!:\!X^{\prime}\setminus Z^{\prime}\longrightarrow X\setminus Z$
	is an isomorphism in $\mathrm{Sch_{fs}^{vLog}\!/k}$ and the underlying
	square of ordinary schemes is an abstract blow-up square.
\end{defn}
Let us suppose that all morphism in square (4.1) are strict. In this
case we will call such a square the \textit{strict log abstract blow-up}.
In the case then underlying maps $p\!:\!\underline{X}^{\prime}\longrightarrow\underline{X}$
and $p_{Z}\!:\!\underline{Z}^{\prime}\longrightarrow\underline{Z}$
are identity maps we will say that the square (4.1) is \textit{geometrically
	trivial}. 
\begin{rem}
	Notice that any morphism of virtual log schemes $(f,f^{\#})\!:\!(X,\mathcal{M}_{X}^{gr})\longrightarrow(Y,\mathcal{M}_{Y}^{gr})$
	can be decomposed as $(X,\mathcal{M}_{X}^{gr})\overset{(id,f^{\#})}{\longrightarrow}(X,f^{*}\mathcal{M}_{Y}^{gr})\overset{(f,id)}{\longrightarrow}(Y,\mathcal{M}_{Y}^{gr})$.
	So the log abstract blow-up (4.1) can be decomposed into strict and
	geometrically trivial parts:
	\[
	\xymatrix{(Z^{\prime},\mathcal{M}_{Z^{\prime}}^{gr})\ar[d]^{(i^{\prime},id)}\ar[r]^{(id,p_{Z}^{\#})} & (Z^{\prime},p_{Z}^{*}\mathcal{M}_{Z}^{gr})\ar[r]^{(p_{Z},id)}\ar[d]^{(i^{\prime},id)} & (Z,\mathcal{M}_{Z}^{gr})\ar[d]^{(i,id)}\\
		(X^{\prime},\mathcal{M}_{X^{\prime}}^{gr})\ar[r]^{(id,p^{\#})} & (X^{\prime},p^{*}\mathcal{M}_{X}^{gr})\ar[r]^{(p,id)} & (X,\mathcal{M}_{X}^{gr})
	}
	\]
\end{rem}
\begin{example}\label{'E4.4}
	Let us denote by $\mathbb{A}_{log}^{n}$ the virtualization of canonical
	log structure on $Spec(\mathrm{k}[\mathbb{N}^{n}])$. Notice that
	$\mathbb{A}_{log}^{n}=\mathbb{A}_{log}^{1}\times...\times\mathbb{A}_{log}^{1}$
	and $\mathcal{M}_{\mathbb{A}_{log}^{1}}^{gr}=j_{*}\mathcal{O}_{\mathbb{G}_{m}}^{*}$,
	$\overline{\mathcal{M}}_{\mathbb{A}_{log}^{1}}^{gr}=i_{*}\mathbb{Z}$
	where $pt\overset{i}{\hookrightarrow}\mathbb{A}^{1}\overset{j}{\hookleftarrow}\mathbb{G}_{m}$.
	Let $pt_{log}$ be a virtualization of a standard log point. Notice
	that $\mathcal{M}_{pt_{log}}^{gr}=\mathrm{k}^{*}\oplus\bZ$
	and the virtualization of an inclusion $pt_{log}\overset{i}{\hookrightarrow}\mathbb{A}_{log}^{1}$
	is the strict morphism. 
	
	Contrasting with ordinary logarithmic geometry the morphism $pt_{log}\overset{i}{\hookrightarrow}\mathbb{A}_{log}^{1}$
	is a section of the projection $\mathbb{A}_{log}^{1}\overset{p}{\longrightarrow}pt_{log}$
	which arises from the generator of the group\linebreak{}
	$H_{et}^{0}(\overline{\mathcal M}_{\mathbb{A}_{log}^{1}}^{gr})=\mathbb{Z}$. 
	
	It is easy to see that for any $X\in\mathrm{Sch_{fs}^{vLog}\!/k}$
	the square 
	\[
	\xymatrix{X\times pt_{log}\ar[r]^{i}\ar[d] & X\times\mathbb{A}_{log}^{1}\ar[d]\\
		X\ar[r]^{i} & X\times\mathbb{A}^{1}
	}
	\]
	is a geometrically trivial log abstract blow-up.
\end{example}
\begin{defn}
	Let $\mathcal{F}$ be a presheaf of complexes on $\mathrm{Sch_{fs}^{vLog}\!/k}$.
	We will say that $\mathcal{F}$ \textit{admits logarithmic h-descent}
	if $\mathcal{F}$ admits strict h-descent and for any log abstract
	blow-up of the form \ref{'F4.1} there is the exact triangle of complexes
	\[
	\mathcal{F}(X)\longrightarrow\mathcal{F}(X^{\prime})\oplus\mathcal{F}(Z)\longrightarrow\mathcal{F}(Z^{\prime}).
	\]
\end{defn}
The main result of this section is the following:
\begin{thm}\label{'T4.6}
	Let $\mathcal{F}$ and $\mathcal{F}^{\prime}$ are the presheaves
	of complexes on $\mathrm{Sch_{fs}^{vLog}\!/k}$ and $\varphi\!:\!\mathcal{F}\longrightarrow\mathcal{F}^{\prime}$
	be a map. Suppose that both $\mathcal{F}$ and $\mathcal{F}^{\prime}$
	admit logarithmic h-descent and for any $\underline{X}\times pt_{log}^{n}$ with $\underline X\in \rm Sch/k$
	the map 
	\[
	\varphi_{\underline{X}\times pt_{log}^{n}}:\mathcal{F}(\underline{X}\times pt_{log}^{n})\longrightarrow\mathcal{F}^{\prime}(\underline{X}\times pt_{log}^{n})
	\]
	is a quasi-isomorphism. Then $\varphi_{Y}$ is a quasi-isomorphism
	for any virtual fs logarithmic scheme $Y$.
\end{thm}
\textit{\negmedspace{}\negmedspace{}\negmedspace{}\negmedspace{}\negmedspace{}Proof:}
Notice that by additivity we may assume that $\mathcal{F}^{\prime}=0$.
Let us denote by $\mathscr{C}$ the full subcategory of $\mathrm{Sch_{fs}^{vLog}\!/k}$
generated by all virtual log schemes of the form $\underline{X}\times pt_{log}^{n}$.
The proof of the theorem contain several steps. At the n-th step we
define the category $\mathscr{C}^{(n)}$ such that $\mathscr{C}\subset\mathscr{C}^{(n)}\subset\mathscr{C}^{(n-1)}$
and check that acyclicity of $\mathcal{F}(X)$ for any $X\in\mathscr{C}^{(n)}$
entails  acyclicity of $\mathcal{F}$. 

\textit{\negmedspace{}\negmedspace{}\negmedspace{}\negmedspace{}\negmedspace{}Step
	1:} For an ordinary scheme $X$ let us fix open immersion $j:U\hookrightarrow X$
, local system $\Lambda$ on $U$ and an element $\xi$ of $Ext^{1}(j_{\#}\Lambda,\mathcal{O}_{X}^{*})$.
We will denote by $j_{\#}X_{\xi}$ the virtual log scheme corresponding
to $\xi$. Let $\mathscr{C}^{(1)}$ be the category containing all
$j_{\#}X_{\xi}$. Now, let us suppose that $\mathcal{F}|_{\mathscr{C}^{(1)}}$
is acyclic. Then $\mathcal{F}(Y)$ is acyclic for any $Y\in\mathrm{Sch_{fs}^{vLog}\!/k}$.
Indeed, we can check that statement using induction by dimension of
$\underline{Y}$. Let $\mathrm{dim}(\underline{Y})=0$. Then any constructible
sheaf on $\underline{Y}$ is a local system and $\mathcal{F}(Y)\simeq0.$
Now, suppose that for all $Y\in\mathrm{Sch_{fs}^{vLog}\!/k}$ with
$\mathrm{dim}(\underline{Y})<n$ we have $\mathcal{F}(Y)\simeq0$.
Let $X$ be an fs virtual log scheme with $\mathrm{dim}(\underline{X})=n$
. Recall that the ghost sheaf $\overline{\mathcal{M}}_{X}^{gr}$ is
constructible. Let $...\subset \underline{Z}_{2} \subset \underline{Z}_{1} \subset \underline{X}$ 
be the corresponding stratification and let $U=X\diagdown Z_{1}$.
So $\Lambda=\overline{M}_{X}^{gr}\mid_{U}$ is a local system and
we can define the log scheme $j_{\#}X_{\xi}$ by the Cartesian diagram
\[
\xymatrix{\mathcal{M}_{j_{\#}X}^{gr}\ar[r]\ar[d] & j_{\#}\Lambda\ar[d]\\
	\mathcal{M}_{X}^{gr}\ar[r] & \overline{\mathcal{M}}_{X}^{gr}
}
\]
By definition we have the log abstract blow-up
\[
\xymatrix{Z_{1}\ar[r]\ar[d] & X\ar[d]\\
	\underline{Z}_{1}\ar[r] & j_{\#}X_{\xi}
}
\]
But $\mathcal{F}(\underline{Z}_{1}),$$\mathcal{F}(Z_{1})$ and $\mathcal{F}(j_{\#}X)$
are acyclic so $\mathcal{F}(X)\simeq0.$

\textit{\negmedspace{}\negmedspace{}\negmedspace{}\negmedspace{}\negmedspace{}Step
	2:} Now, let $\mathscr{C}^{(2)}$ be the category generated by all
$j_{\#}Y_{\xi}$ with constant $\Lambda$ and $\xi=0$. We should
check that if $\mathcal{F}|_{\mathscr{C}^{(2)}}$ is acyclic then
$\mathcal{F}$ is acyclic. By previous it suffice to show that $\mathcal{F}|_{\mathscr{C}^{(1)}}$
is acyclic. 

Let $j_{\#}X_{\xi}$ be a log scheme corresponding to some $\Lambda$
and $\xi$. Using strict $h$-descent we may assume that $X$ is normal.
So, by definition, for any geometric point $\overline{x}$ the fiber
$\Lambda_{\overline{x}}$ is a finitely generated abelian group. The
action of $\pi_{1}^{et}(U,\overline{x})$ is continuous so the stabilizer
of each element is an open subgroup. Hence the kernel $\pi_{1}^{et}(U,\overline{x})\longrightarrow\mathrm{Aut}(\Lambda_{\overline{x}})$
is open and the image is a finite group. So there exist the Galois
cover $\widetilde{U}\overset{\pi}{\longrightarrow}U$ such that $\pi^{*}\Lambda$
is constant. 

Now, let $X^{\prime}\overset{\nu}{\longrightarrow}X$ be the normalization
of $\widetilde{U}$ in $X$. Observe that $\nu$ is a $h$-cover.
Indeed, by Lemma 29.53.15 of \cite{'key-1} $\nu$ is finite and $\mathrm{Im}(\nu)\supset\mathrm{Im}(\widetilde{U})=U$.
So $v$ is a finite surjective morphism. Notice that a normalization
commute with smooth base change (see Lemma 37.17.2 of \cite{'key-1}).
Moreover, the normalization of $\widetilde{U}$ in $U$ coincide with
$U$ because $\pi$ is finite. So $X^{\prime}\times_{X}U\simeq\widetilde{U}$
and $v^{-1}(j_{\#}\Lambda)\simeq\widetilde{j}_{\#}(\pi^{-1}\Lambda)=\widetilde{j}_{\#}(\mathbb{Z}^{n})$.
Here $\widetilde{j}$ is an open inclusion of $\widetilde{U}$ into
$X^{\prime}$. 

The log structure $v^{*}\mathcal{M}_{j_{\#}X}^{gr}$ corresponding
to some element $\xi$ of 
\[
\mathrm{Ext^{1}}(\tilde{j}_{\#}\mathbb{Z}^{n},\mathcal{O}_{X^{\prime}}^{*})\simeq\mathrm{Ext^{1}}(\mathbb{Z}^{n},\mathcal{O}_{\widetilde{U}}^{*})=H_{Zar}^{1}(\widetilde{U},\mathcal{O}^{*})^{n}
\]
So there exists a Zariski cover $\{W_{i}\longrightarrow X^{\prime}\}$
such that $\xi|_{W_{i}}=0$. Now, let us compose the normalization
with h-cover $X^{\prime}\overset{\nu}{\longrightarrow}X$ and the
Zariski cover $\{W_{i}\longrightarrow X^{\prime}\}$. We get the h-cover
$Y\longrightarrow X$ such that for any product $Y\times_{X}...\times_{X}Y$
the pullback of the log structure has the form $j_{\#}Y_{\xi}$ with
constant $\Lambda$ and $\xi=0$. So using strict h-descent we get
$\mathcal{F}(j_{\#}X_{\xi})=0$.

\textit{\negmedspace{}\negmedspace{}\negmedspace{}\negmedspace{}\negmedspace{}Step
	3:} Let us prove the main statement. For all $j_{\#}X_{\xi}$ with
constant $\Lambda=\mathbb{Z}^{n}$ and $\xi=0$ there exist the log
abstract blow-up 
\[
\xymatrix{\underline{Z}_{1}\times pt_{log}^{n}\ar[r]\ar[d] & \underline{X}\times pt_{log}^{n}\ar[d]\\
	\underline{Z}_{1}\ar[r] & j_{\#}X_{\xi}
}
\]
So we can use induction by dimension of $\underline{X}$. $\square$

\section{Logarithmic motivic cohomology}
Let us use the notation of section 6. Let R be a commutative ring. 
\begin{defn}
	For $q=0,1$ we define the complexes $R\Gamma(X,\mathrm{R}^{\mathrm{log}}(q))$ by the formulas
	$$R\Gamma(X,\mathrm{R}^{\mathrm{log}}):=\rm Hom^\bullet_{\mathrm{DA_{et}}(X,\mathrm{R})}([X],[X]); \ \ \ \mathit{ R}\Gamma(X,\mathrm{R}^{\mathrm{log}}(1)):= Hom^\bullet_{\mathrm{DA_{et}}(X,\mathrm{R})}([X],M^{gr}_{\it X}).$$
	and the \textit{logarithmic motivic cohomology} as
	$$H^{p,q}_{log}(X,\mathrm{R}):=H^p(R\Gamma(X,\mathrm{R}^{\mathrm{log}}(q))).$$ 
\end{defn}
Now, let us assume that $\mathrm{R}=\bQ$. In this case we can define $H^{p,q}_{log}$ for $q>1$ using motivic symmetric powers.
\begin{rem}
	( \textit{motivic symmetric powers} ) Let $S$ be a scheme. For $n>1$
	and $F\in\mathrm{DA_{et}}(S,\mathbb{Q})$ let us consider the endomorphism
	$\varphi:=\frac{1}{n!}\Sigma_{\sigma\in\Sigma_{n}}\sigma$
	of $F^{\otimes n}$. Notice that $\varphi^{2}=\varphi$ and, consequently,
	the motivic sheaf $F^{\otimes n}$ can be canonically factorized as
	the direct sum 
	\begin{equation}\label{'F5.4}
		F^{\otimes n}\simeq\mathrm{im}(\varphi)\oplus\mathrm{ker}(\varphi).
	\end{equation}
	The direct term $\mathrm{im}(\varphi)$ is usually denoted by $S^{n}F$
	and called the n-th symmetric power of $F$ (see \cite{'key-28} for
	details). The factorization \ref{'F5.4} gives rise to the
	endofunctor 
	\[
	S^{n}\!:\!\mathrm{DA_{et}}(S,\mathbb{Q})\longrightarrow\mathrm{DA_{et}}(S,\mathbb{Q})
	\]
	of triangulated categories. Notice that this functor is not triangulated. 
	
	On the other hand, the n-th symmetric power can be considered as the DG-functor.
	Indeed, we can define the map of complexes $\mathrm{Hom}^{\bullet}(E,F)\longrightarrow\mathrm{Hom}^{\bullet}(S^{n}E,S^{n}F)$
	as the composition of the map $\mathrm{Hom}^{\bullet}(E,F)\longrightarrow\mathrm{Hom}^{\bullet}(E^{\otimes n},F^{\otimes n})$
	with the canonical projection %\linebreak{}
	$\mathrm{Hom}^{\bullet}(E^{\otimes n},F^{\otimes n})\longrightarrow\mathrm{Hom}^{\bullet}(S^{n}E,S^{n}F)$.
	It is easy to see that such a collection of maps preserve compositions
	and associativity. 
	\\We can also define the motivic exterior powers $\Lambda^n$ replace $\varphi$ with $\frac{1}{n!}\Sigma_{\sigma\in\Sigma_{n}}(-1)^{|\sigma|}\sigma$. Note that $ S^n(F[-1])\simeq \Lambda^nF[-n]$ for any $F$ and $n$. 
\end{rem}
\begin{defn}
	Using motivic symmetric powers we put
	\begin{equation}
		R\Gamma(X,\mathrm{\mathbb{Q}^{log}}(q)\!)\!\overset{\mathrm{def}}{=}\!\mathrm{Hom}_{\mathrm{DA_{et}}(X,\mathbb{Q})}^{\bullet}([X],S^{q}(\mathrm{M_{\mathit{X}}^{gr}}[-1])).
	\end{equation}
	%We also define the motivic cohomology $H^{p,q}(X,\mathbb{Q})$ of
	%a virtual log scheme $X$ as the cohomology of the complex of $\mathbb{Q}$-vector
	%spaces $R\Gamma(\mathrm{\mathbb{Q}_{\mathit{X}}^{log}}(q)\!)$
	and
	\[
	H^{p,q}_{log}(X,\mathbb{Q})\!\overset{\mathrm{def}}{=}\!H^{p}(R\Gamma(X,\mathrm{\mathbb{Q}^{log}}(q)\!)\!).
	\]
	Further, for simplicity, we will use the sheaves $\mathrm{\Lambda^{q}M_{\mathit{X}}^{gr}}$
	instead $S^{q}(\mathrm{M_{\mathit{X}}^{gr}}[-1])$.
\end{defn}

\subsection*{Functoriality.} 
Let us fix $q$. Let $(f,f^{\#}):(X,\mathcal M^{gr}_X)\to (Y,\mathcal M^{gr}_Y)$ be the
morphism of virtual log schemes. Now, we explain how to construct the map 
\begin{equation}\label{'key}
	R\Gamma(f,f^{\#}):R\Gamma(Y,\bQ^{\rm log}\it (q))\ra R\rm \Gamma(\it X,\bQ^{\rm log}\it (q)).
\end{equation}
Firstly, note that (\ref{'F5.3}) maps $f^*\cM_Y$ to $f^*\rm M^{gr}_{\it Y}$ by Corollary \ref{'C5.3} and Remark \ref{'R5.4}. Hence $f^{\#}$ maps to 
$$\!\mathrm{\Lambda^{q}}f^{*}\mathrm{M_{\mathit{Y}}^{gr}}\longrightarrow\mathrm{\Lambda^{q}M_{\mathit{X}}^{gr}}.$$ 
On the other hand, $f^{*}$ is monoidal so $\mathrm{\Lambda^{q}}f^{*}\mathrm{M_{\mathit{Y}}^{gr}}\simeq f^{*}\mathrm{\Lambda^{q}M_{\mathit{Y}}^{gr}}$
and there exists the map 
\[
\mathrm{\Lambda^{q}M_{\mathit{Y}}^{gr}}\longrightarrow Rf_{*}\mathrm{\Lambda^{q}M_{\mathit{X}}^{gr}}.
\]
It remains to apply the functor $\mathrm{Hom}_{\mathrm{DA_{et}( Y},\mathbb{Q})}^{\bullet}([Y],-)$
and use the natural quasi-isomorphism 
$$
\mathrm{Hom}_{\mathrm{DA_{et}(}Y,\mathbb{Q})}^{\bullet}(\mathbf{M}(Y),Rf_{*}\mathrm{\Lambda^{q}M_{\mathit{X}}^{gr}})\simeq\mathrm{Hom}_{\mathrm{DA_{et}(}X,\mathbb{Q})}^{\bullet}(f^{*}\mathbf{M}(Y),\mathrm{\Lambda^{q}M_{\mathit{X}}^{gr}}).
$$
Note that we construct the set of functors 
\begin{equation}\label{'F6.2}
	[R\Gamma(-,\bQ^{\mathrm{log}}(q))]:\rm Sch^{vLog}/k\ra D(Vect_\bQ);
\end{equation}
$$H^{p,q}_{log}:\rm Sch^{vLog}_{fs}/k\ra Vect_\bQ.$$
Moreover, after choosing a suitable model structure, (\ref{'F6.2}) can be lifted to the presheaf of complexes  
$$R\Gamma(-,\bQ^{\mathrm{log}}(q)):\rm Sch^{vLog}/k\ra Ch^\bullet(Vect_\bQ).$$ Indeed, it is equivalent to show that \ref{'F6.2} can be lifted to a $(\infty,1)$-functor ${\rm Sch^{vLog}/k\ra Ch^\bullet(Vect_\bQ).}$ By analogy with the construction of $[-]^{log}$ , it can be done using the following diagram 
\[
\xymatrix{(\int_{\mathrm{Sch/k^{op}}}\mathrm{DA_{et,\mathbb{Q}}})^{\mathrm{op}}\ar[r]^{\ \ \varphi} & \int_{\mathrm{Sch/k}}\mathrm{DA_{et,\mathbb{Q}}^{op}}\ar[r]^{\!\!\!\!\!\!\!\!\!\!\! \int R\Gamma} & \mathrm{Sch/k\times \rm D(Vect_\bQ)}\ar[d]^{\rm pr_{2}}\\
	\mathrm{Sch_{fs}^{vLog}/k}\ar[u]_{\int\Lambda^{q}\Psi}\ar[rr]^{R\Gamma(-,\mathrm{\mathbb{Q}^{log}}(q))} &  & \mathrm{D(Vect_\bQ)}
}
.
\]
Again, here $\int$ mean the $(\infty,1)$-Grothendieck construction (see Appendix \ref{B}) and $\varphi$ is the equivalence of categories from Proposition \ref{PB.4}.

\subsection*{$\Al$-invariance and log h-descent.} 
Note that the presheaves $R\Gamma(-,\bQ^{\mathrm{log}}(q))$ is $\Al$-invariant by construction. It should also be noted that $R\Gamma(-,\bQ^{\mathrm{log}}(q))$ admit log $h$-descent. Indeed, the presheaves
$R\Gamma(-,\bQ^{\mathrm{log}}(q))$ admit strict $h-$descent iff for
any $X\in\mathrm{Sch_{fs}^{vLog}/k}$ the $\Gm$ - spectrum
$\mathrm{\Lambda^{q}M_{\mathit{X}}^{gr}}$ admits $h-$descent. But
any $\Gm$ - spectrum $E\in\mathrm{DA_{et}(}\underline{X},\mathbb{Q})$
admits $h-$descent by results of \cite{'key-22}. Let 
\begin{equation}\label{'F6.4}
	\xymatrix{Z^{\prime}\ar[r]\ar[d] & X^{\prime}\ar[d]\\
		Z\ar[r] & X
	}
\end{equation}
be a log abstract blow-up square. We should check that the sequence
\[
R\Gamma(X,\bQ^{\mathrm{log}}(q))\longrightarrow R\Gamma(X',\bQ^{\mathrm{log}}(q))\oplus R\Gamma(Z,\bQ^{\mathrm{log}}(q))\longrightarrow R\Gamma(Z',\bQ^{\mathrm{log}}(q))
\]
is an exact triangle. The presheaf $R\Gamma(-,\bQ^{\mathrm{log}}(q))$
admits strict $h-$descent so we may assume that the
square \ref{'F6.4} is geometrically trivial (i.e $\underline{X}^{\prime}=\underline{X}$
and $\underline{Z}^{\prime}=\underline{Z}$). Let us denote by $i$
the inclusion $\underline{Z}\hookrightarrow\underline{X}$ and let
$j:\underline{U}\hookrightarrow\underline{X}$ be the inclusion of
the complement. It suffice to show that the sequence 
\begin{equation}\label{'F6.5}
	\mathrm{\Lambda^{q}M_{\mathit{X}}^{gr}}\longrightarrow\mathrm{\Lambda^{q}M_{\mathit{X^{\prime}}}^{gr}}\oplus i_{*}\mathrm{\Lambda^{q}M_{\mathit{Z}}^{gr}}\longrightarrow i_{*}\mathrm{\Lambda^{q}M_{\mathit{Z^{\prime}}}^{gr}}
\end{equation}
is an exact triangle. Applying the localization sequence we get the
map of exact triangle
\begin{equation}\label{'F6.6}
	\xymatrix{j_{\#}\mathrm{\Lambda^{q}M_{\mathit{U}}^{gr}}\ar[r]\ar[d]^{\wr} & \mathrm{\Lambda^{q}M_{\mathit{X}}^{gr}}\ar[d]\ar[r] & i_{*}\mathrm{\Lambda^{q}M_{\mathit{Z}}^{gr}}\ar[d]^{\,\,\,\,\,\,\,\,\,\,\,\,\,.}\\
		j_{\#}\mathrm{\Lambda^{q}M_{\mathit{U^{\prime}}}^{gr}}\ar[r] & \mathrm{\Lambda^{q}M_{\mathit{X^{\prime}}}^{gr}}\ar[r] & i_{*}\mathrm{\Lambda^{q}M_{\mathit{Z^{\prime}}}^{gr}}
	}
\end{equation}
By definition of log abstract blow-up the left horizontal map is a
quasi-isomorphism. Hence the exactness of (\ref{'F6.5})
follows formally from (\ref{'F6.6}). %
\begin{comment}
Let $\xymatrix{C[-1]\ar[r]^{\varphi}\ar[d] & A\ar[r]\ar[d]^{\wr\psi} & B\ar[d]^{\,\,\,\,\,\,\,\,\,\,\,\,\,.}\\
C^{\prime}[-1]\ar[r] & A^{\prime}\ar[r] & B^{\prime}
}
$be the map of exact triangles. The the exactness of $B\rightarrow B^{\prime}\oplus C\rightarrow C^{\prime}$
follows from the map of exact triangles $\xymatrix{C[-1]\ar[r]^{id}\ar[d]^{\psi\varphi} & C[-1]\ar[r]\ar[d]^{0} & 0\ar[d]^{\,\,\,\,\,\,\,\,\,\,\,\,\,.}\\
A^{\prime}\ar[r] & B^{\prime}\ar[r] & C^{\prime}
}
$
\end{comment} 

\subsection*{Weight one motivic cohomology}

In ordinary theory of motivic cohomology we have $H^{p,1}(X,\mathbb{Q})=H^{p-1}(X,\mathcal{O}^{*}\otimes_{\mathbb{Z}}\mathbb{Q})$
for a smooth scheme $X$. One would expect a similar result for the
motivic cohomology of a log schemes. At least for smooth and log
smooth $X$ we have the following result.
\begin{prop}\label{'P6.2}
	Let $X=(\underline{X},\cO^*\overset{\alpha}{\hookrightarrow} \mathcal{M}_{X}^{gr})$ be a smooth and log smooth scheme.
	Let us apply
	the natural transformation 
	$
	\mathrm{Ex_{\mathit{\underline{X}}}^{log}}\longrightarrow\mathrm{L}_{\mathbb{A}^{1}}\mathrm{Ex_{\mathit{\underline{X}}}^{log}}\longrightarrow\mathrm{\Omega_{\mathbb{G}_m}^{\infty}}\mathrm{\Sigma_{\mathbb{G}_m}^{\infty}}\mathrm{L}_{\mathbb{A}^{1}}\mathrm{Ex_{\mathit{\underline{X}}}^{log}}
	$
	to $\alpha$. Then the map
	\[
	\mathrm{Ex_{\mathit{\underline{X}}}^{log}}\mathcal{M}_{X}^{gr}\otimes \bQ \longrightarrow\mathrm{\Omega_{\mathbb{G}_m}^{\infty}}\mathrm{M_{\mathit{X}}^{gr}}
	\]
	is a quasi-isomorphism in $\mathrm{D(Sh_{et}(Sm/}\underline X,\mathbb{Q})$.
	In particular, 
	\[
	H^{p,1}(X,\mathbb{Q})\iso H_{et}^{p-1}(\underline{X},\mathcal{M}_{X}^{gr}\otimes \mathbb{Q}).
	\]
\end{prop}
\textit{Proof:}  Notice that $\mathrm{Ex}_{Y}(i_{*}\mathbb{Q})=i_{*}\mathbb{Q}$
by construction and the functor 
\[
i_{*}:\mathrm{Sh_{et}(Sm}/D,\mathbb{Q})\longrightarrow\mathrm{Sh_{et}(Sm/}Y,\mathbb{Q})
\]
is exact (because for any $W$ over $Y$ the functor $\mathrm{Res}_{W}(i_{*})$
is exact). So 
\[
H_{et}^{p}(Y,\mathrm{Ex}_{Y}(i_{*}\mathbb{Q}))=H_{et}^{p}(D,\mathbb{Q})=\begin{cases}
	\mathbb{Q} & p=0\\
	0 & p>0
\end{cases}
\]
and $\mathrm{Ex}_{Y}(i_{*}\mathbb{Q})$ is $\mathbb{A}^{1}$-invariant. Using étale descent we may assume that the log structure $\alpha$ admits a global chart. So it
suffice to check that the units of adjunction
$$\eta_{\cO^*}:\cO^*\to \Omega_{\mathbb{G}_m}^\infty \Sigma_{\mathbb{G}_m}^\infty \cO^*$$
$$\eta_{i_*\bQ}:i_*\bQ\to \Omega_{\mathbb{G}_m}^\infty \Sigma_{\mathbb{G}_m}^\infty i_*\bQ$$ 
are isomorphisms. 
\\Let $F$ be $\cO^*$ or $i_*\bQ$. Using the localization and Voevodsky Cancellation theorem we see that, in both cases, $F\iso \Omega_{\mathbb{G}_m}^\infty \Sigma_{\mathbb{G}_m}^\infty F.$ Moreover, in both cases $\rm End_{DA^{eff}_{et}(\it X,\bQ)}(\it F)\iso \bQ.$ So $\eta_F$ is an isomorphism or 0. The second case is impossible
because $\eta_F$ corresponds to $\mathrm{id}_{\Sigma_{\Gm}^\infty F}$ which generates $\rm End_{DA^{eff}_{et}(\it X,\bQ)}(\Sigma_{\Gm}^\infty \it F).$
$\square$

\section{Motivic cohomology of a log projective line}

In this section we prove $\mathbb{P}_{log}^{1}$-invariance of log
motivic cohomology. By definition $\mathbb{P}_{log}^{1}$ is the gluing
of $\mathbb{A}_{log}^{1}$ and $\mathbb{A}^{1}$ along $\mathbb{G}_{m}$.
We have the following commutative diagram 
\[
\xymatrix{ & \mathbb{P}_{log}^{1}\ar[d] & pt_{log}\ar[l]\ar[d]^{\,\,\,\,\,\,.}\\
	\mathbb{A}^{1}\ar[r]^{j}\ar[ur] & \mathbb{P}^{1} & \infty_{\mathrm{k}}\ar[l]_{i}
}
\]

Let $\mathcal{O}_{\mathbb{P}^{1}}^{*}\hookrightarrow\mathcal{M}^{gr}$
be the virtual log structure on $\mathbb{P}^{1}$ associated with
$\mathbb{P}_{log}^{1}$. Notice that $\mathcal{M}^{gr}=j_{*}\mathcal{O}_{\mathbb{A}^{1}}^{*}$.
Let us give the explicit description of the motivic sheaves $\mathrm{\Lambda^{q}M^{gr}}$.
\begin{lem}\label{'L7.1}
	Let $Rj_{*}:\mathrm{DA_{et}^{eff}(}\mathbb{A}^{1},\mathbb{Q})\longrightarrow\mathrm{DA_{et}^{eff}(}\mathbb{P}^{1},\mathbb{Q})$
	be the direct image functor. Then
	\[
	\Omega_{\mathbb{L}}^{\infty}\mathrm{M^{gr}}\simeq Rj_{*}\mathcal{O}_{\mathrm{Sm/}\mathbb{A}^{1}}^{*}.
	\]
\end{lem}
\textit{Proof:} By Yoneda Lemma we should check that 
\[
\mathrm{Hom}_{\mathrm{DA_{et}(}\mathbb{P}^{1})}^{\bullet}(Y,\mathrm{M^{gr}})\simeq\mathrm{Hom}_{\mathrm{DA_{et}^{eff}(}\mathbb{A}^{1})}^{\bullet}(Y\times_{\mathbb P^1}{\mathbb{A}^{1}},\mathcal{O}_{\mathrm{Sm/}\mathbb{A}^{1}}^{*})
\]
for any $Y\in\mathrm{Sm/}\mathbb{P}^{1}$. Then by Proposition \ref{'P6.2} it
suffice to show that $R^{k}j_{*}\mathcal{O}_{\mathbb{A}^{1}}^{*}\in\mathrm{Sh(}\mathbb{P}_{et}^{1},\mathbb{Q})$
vanish for $k>0$. Moreover, by smoothness of $\mathbb{P}^{1}$ it
enough to show $R^{1}j_{*}\mathcal{O}_{\mathbb{A}^{1}}^{*}=0$.
%because for smooth scheme H^k(\mathcal O^*)=0 for k>1
The sheaf $R^{1}j_{*}\mathcal{O}_{\mathbb{A}^{1}}^{*}$ is an étale
sheafification of the presheaf 
\[
\widetilde{\mathbf{Pic}}:X\longmapsto\mathrm{Pic}(X\!\times_{\mathbb{P}^{1}}\!\mathbb{A}^{1}).
\]
Let us denote by $\mathbf{Pic}$ the presheaf which maps $X$ to $\mathrm{Pic}(X)$.
Let us consider the short exact sequence $\mathcal{O}_{\mathbb{P}^{1}}^{*}\longrightarrow j_{*}\mathcal{O}_{\mathbb{A}^{1}}^{*}\longrightarrow i_{*}\mathbb{Q}$.
Applying long exact sequence of cohomology we get the surjection of
presheaves $\mathbf{Pic}\longrightarrow\widetilde{\mathbf{Pic}}$
\begin{comment}
because $H_{et}^{1}(X,i_{*}\mathbb{Q})=H_{et}^{1}(X_{0},\mathbb{Q})=0$
\end{comment}
. But the Zariski sheafification of $\mathbf{Pic}$ is zero. Indeed, by
definition of plus construction we have 
\[
\mathrm{Pic^{+}}(X)=\underrightarrow{\mathrm{holim}}\mathrm{(ker}(\Pi_{i}\mathrm{Pic}(U_{i})\overset{\partial_{U_{i}}}{\longrightarrow}\Pi_{i,j}\mathrm{Pic}(U_{i}\cap U_{j}))
\]
where colimit taken by all Zariski covers of $X$. Let $\varphi\in\mathrm{ker(}\partial_{U_{i}})$
for some cover. Suppose $\varphi$ is represented by collection of vector
bundles $E_{i}$ on $U_{i}$. For each $U_{i}$ there exist a Zariski
cover $\{W_{ki}\longrightarrow U_{i}\}$ such that $E_{i}|_{W_{ki}}$ are trivial.
So $\varphi$ maps to zero by the canonical map %\linebreak{}
$\mathrm{ker(}\partial_{U_{i}})\longrightarrow\mathrm{ker(}\partial_{W_{ki}})$,$[\varphi]=0$
and $\mathrm{Pic^{+}}(X)=0$. $\square$ 
\begin{cor}
	The motivic spectrum $\mathrm{M^{gr}}$ coincide with $Rj_{*}\mathbb{Q}(1)[1]$.
\end{cor}
\textit{Proof:} This is a formal consequence of Verdier duality and
the localization property. Let us denote by $\pi:\mathbb{P}^{1}\longrightarrow Spec(\mathrm{k})$
the map to the point. By Verdier duality
\[
\mathbb{Q}(1)\simeq\pi^{*}\mathbb{Q}(1)\simeq\pi^{!}\mathbb{Q}[-2].
\]
So $i^{!}\mathbb{Q}(1)[1]\simeq i^{!}\pi^{!}\mathbb{Q}[-1]=\mathbb{Q}[-1]$
and 
\[
\mathrm{Hom}_{\mathrm{DA_{et}(}\mathbb{P}^{1})}^{\bullet}(i_{*}\mathbb{Q}[-1],\mathbb{Q}(1)[1])\simeq\mathrm{Hom}_{\mathrm{DA_{et}(k)}}^{\bullet}(\mathbb{Q}[-1],\mathbb{Q}[-1])=\mathbb{Q}.
\]
Hence there exists only one non-split extension $\mathbb{Q}(1)[1]\rightarrow \mathcal F \rightarrow i_{*}\mathbb{Q}$. 

On the other hand, applying the localization to the motivic sheaf
$\mathbb{Q}(1)[1]\in\mathrm{DA_{et}(}\mathbb{P}^{1})$ we get the
triangle $i_{*}\mathbb{Q}[-1]\rightarrow\mathbb{Q}(1)[1]\rightarrow Rj_{*}\mathbb{Q}(1)[1]$
and the non-split extension 
\[
\mathbb{Q}(1)[1]\rightarrow Rj_{*}\mathbb{Q}(1)[1]\rightarrow i_{*}\mathbb{Q}.
\]
So if $\mathrm{M^{gr}}\neq Rj_{*}\mathbb{Q}(1)[1]$ then $\mathrm{M^{gr}}$
and, consequently, $\Omega_{\mathbb{L}}^{\infty}\mathrm{M^{gr}}$
are split. We get the contradiction with Lemma \ref{'L7.1}. $\square$

Now, let $A$ be an abelian category and
\[
K^{0}\longrightarrow K^{1}\longrightarrow...\longrightarrow K^{n-1}\longrightarrow K^{n}
\]
be a complex concentrated in the finite numbers of components. Let
\[
F^{n}K^{\bullet}\subset F^{n-1}K^{\bullet}\subset....\subset F^{0}K^{\bullet}=K^{\bullet}
\]
be the silly filtration with 
\[
F^{j}K^{\bullet}=\begin{cases}
	K^{i} & \mathrm{if}\,\,j\leq i\leq n\\
	0 & \mathrm{otherwise}
\end{cases}
\]

Let $\mathfrak{E}:D(A)\longrightarrow\mathcal{T}$ be a triangulated
functor to some triangulated category. Suppose that for $m<n$ we
have $\mathfrak{E}(K^{0})=\mathfrak{E}(K^{1})=...=\mathfrak{E}(K^{m})=0$.
Then $\mathfrak{E}(K^{\bullet})\simeq\mathfrak{E}(F^{m+1}K^{\bullet})$.
This assertion allow us to construct the Koszul complexes and compute the
external powers $\text{\ensuremath{\mathrm{\Lambda^{q}M^{gr}}}}$.

Let us consider the complex $\mathcal{O}_{\mathbb{Q}}^{*}\longrightarrow\mathcal{M}_{\mathbb{Q}}^{gr}\in\mathrm{Ch^{\bullet}(Sh_{et}(Sm/}\mathbb{P}^{1},\mathbb{Q}))$.
The q-th exterior power of $\mathcal{O}_{\mathbb{Q}}^{*}\longrightarrow\mathcal{M}_{\mathbb{Q}}^{gr}$
is the complex
\[
\underset{-q}{S^{q}\mathcal{O}_{\mathbb{Q}}^{*}}\longrightarrow\underset{-q+1}{S^{q-1}\mathcal{O}_{\mathbb{Q}}^{*}\otimes\mathcal{M}_{\mathbb{Q}}^{gr}}\longrightarrow...\longrightarrow\underset{-1}{\mathcal{O}_{\mathbb{Q}}^{*}\otimes\Lambda^{q-1}\mathcal{M}_{\mathbb{Q}}^{gr}}\longrightarrow\underset{0}{\Lambda^{q}\mathcal{M}_{\mathbb{Q}}^{gr}}
\]
Notice that any sheaf of vector spaces are flat. It follows that the functor %\linebreak{}
$\mathrm{Ch^{\bullet}(Sh_{et}(Sm/}\mathbb{P}^{1},\mathbb{Q}))\longrightarrow\mathrm{D(Sh_{et}(Sm/}\mathbb{P}^{1},\mathbb{Q})$
is monoidal and commute with $\Lambda^{q}$ and $S^{q}$. 
Using this and monoidality of $\mathrm{\Sigma_{\mathbb{L}}^{\infty}}\mathrm{L}_{\mathbb{A}^{1}}$
we get

\[
\mathrm{\Sigma_{\mathbb{L}}^{\infty}}\mathrm{L}_{\mathbb{A}^{1}}(S^{q-k}\mathcal{O}_{\mathbb{Q}}^{*}\otimes\Lambda^{k}\mathcal{M}_{\mathbb{Q}}^{gr})\simeq\mathrm{\Lambda^{q-k}\mathbb{Q}}(1)[n]\otimes\mathrm{\Lambda^{k}M^{gr}=0}
\]
for $k<q-1$. So by the previous 
\[
\mathrm{\Lambda^{q}(Cone(\mathbb{Q}(1)[1]\longrightarrow M^{gr}))\!\simeq\!Cone(}\mathrm{\Lambda^{q-1}M^{gr}}(1)\longrightarrow\mathrm{\Lambda^{q}M^{gr}}).
\]
But $\mathrm{\Lambda^{q}(Cone(\mathbb{Q}(1)[1]\rightarrow M^{gr}}))\!\simeq\!\mathrm{\Lambda^{q}}(i_{*}\mathbb{Q})$
and $\mathrm{\Lambda^{q}}(i_{*}\mathbb{Q})\!=\!0$. Indeed, $j^{*}\mathrm{\Lambda^{q}}(i_{*}\mathbb{Q})\!=\!\mathrm{\Lambda^{q}}(j^{*}i_{*}\mathbb{Q})\!=\!0$
and $i^{*}\mathrm{\Lambda^{q}}(i_{*}\mathbb{Q})\!=\!\mathrm{\Lambda^{q}}(\mathbb{Q})\!=\!0$.
So we proved the following 
\begin{cor}\label{'C7.3}
	For $\mathrm{M^{gr}}$ the same as before we have $\mathrm{\Lambda^{q}M^{gr}}\simeq Rj_{*}\mathbb{Q}(q)[q]$. In particular, 
	$$H^{p,q}_{log}(\mathbb{P}_{log}^{1},\mathbb{Q})\iso H^{p,q}(\mathbb{A}^{1},\mathbb{Q})\iso H^{p,q}(Spec(\mathrm{k}),\mathbb{Q}). $$
\end{cor}
The last part of the Corollary can be generalized as $\mathbb{P}_{log}^{1}$-invariance
of log motivic cohomology.
\begin{thm}\label{'T7.4}
	Let $\pi:\mathbb{P}^{1}\longrightarrow Spec(\mathrm{k})$ be the canonical map. Then for any virtual fs log scheme $X$  the maps
	\[
	H^{p,q}_{log}(X,\mathbb{Q})\stackrel{\pi_{X}^{*}}{\longrightarrow}H^{p,q}_{log}(X\times\mathbb{P}_{log}^{1},\mathbb{Q})
	\]
	are isomorphisms.
\end{thm}
\textit{Proof:} It enough to show that the morphism $(j_{X},j_{X}^{\#}):X\times\mathbb{A}^{1}\longrightarrow X\times\mathbb{P}_{log}^{1}$
induced an isomorphism on the motivic cohomology. Using resolution of singularities and Theorem \ref{'T4.6} we may assume that $X=\underline{X}\times pt_{log}^{n}$
for some $n$ and $\underline{X}$ smooth over $\mathrm{k}$.
Let us denote by $\mathrm{M^{gr}}$ the motivic spectrum associated with $\mathcal{M}_{\mathit{X\times\mathbb{P}_{log}^{\mathrm{1}}}}^{gr}.$
By construction, the map $(j_{X},j_{X}^{\#})^{*}$
is adjoint to the map 
\[
j_{X}^{*}\mathrm{\Lambda^{q}M^{gr}}\longrightarrow\mathbb{Q}(q)[q]
\]

So we have the commutative diagram 
\[
\xymatrix{ & Rj_{X*}(j_{X}^{*}\mathrm{\Lambda^{q}M^{gr}})\ar[d]^{\,\,\,\,\,\,\,\,\;\,\,\,\,\,\;\,\,\,\,\,\;\;\;\,.}\\
	\mathrm{\Lambda^{q}M^{gr}}\ar[ur]^{\eta}\ar[r]^{(j_{X},j_{X}^{\#})^{*}} & Rj_{X*}(\mathbb{Q}(q)[q])
}
\]

The vertical map is an isomorphism. So it enough to show that
for any $p$ the map 
\[
\eta_*:\mathrm{Hom}_{\mathrm{DA_{et}}}(\mathbb{Q}[-p],\mathrm{\Lambda^{q}M^{gr}})\longrightarrow\mathrm{Hom}_{\mathrm{DA_{et}}}(\mathbb{Q}[-p],Rj_{X*}(j_{X}^{*}\mathrm{\Lambda^{q}M^{gr}}))
\]
is an isomorphism.

Suppose that $n=0$ and $X=\underline{X}$. Then applying smooth
base change
(combine A.5.1.4 and A.5.1.5 of \cite{'key-22})
to the diagram 
\[
\xymatrix{X\times\mathbb{A}^{1}\ar[r]^{\pi_{\mathbb{A}^{1}}}\ar[d]^{j_{X}} & \mathbb{A}^{1}\ar[d]^{j}\\
	X\times\mathbb{P}^{1}\ar[r]^{\pi_{\mathbb{P}^{1}}} & \mathbb{P}^{1}
}
\]
we get $\mathrm{\Lambda^{q}M^{gr}}\simeq Rj_{X*}\mathbb{Q}(q)[q]$.
So $\eta$ is a quasi-isomorphism.

Now, let $n>0$ and $X=\underline{X}\times pt_{log}^{n}$. In this
case $\mathrm{M^{gr}}\simeq Rj_{X*}\mathbb{Q}(1)[1]\oplus\mathbb{Q}^{n}$.
For any direct sum $A\oplus B$ in any tensor $\mathbb{Q}$-linear
category $C$ we have the isomorphism $\mathrm{\Lambda^{q}}(A\oplus B)\simeq\bigoplus_{\mathrm{j}}\mathrm{\Lambda^{j}}A\otimes\mathrm{\Lambda^{q-j}}B$
(see \cite{'key-23}). So it suffice to show that $\eta_*$ is
an isomorphism for each summand. If $\mathrm{j}\neq0$ then 
\[
\mathrm{\Lambda^{j}}(Rj_{X*}\mathbb{Q}(1)[1])\otimes\mathrm{\Lambda^{q-j}}(\mathbb{Q}^{n})=Rj_{X*}\mathbb{Q}(j)[j]\otimes\mathbb{\mathbb{Q}}^{\binom{n}{q-j}}=Rj_{X*}\mathbb{Q}(j)[j]^{\binom{n}{q-j}}.
\]
If $j=0$ then $\mathrm{\Lambda^{q}}(\mathbb{Q}^{n})\simeq\mathbb{\mathbb{Q}}^{\binom{n}{q}}$.
We should check that 
\[
\eta_*:\mathrm{Hom}_{\mathrm{DA_{et}}}(\mathbb{Q}[-p],\mathbb{Q})\longrightarrow\mathrm{Hom}_{\mathrm{DA_{et}}}(\mathbb{Q}[-p],Rj_{X*}j_{X}^{*}\mathbb{Q})
\]
is an isomorphism for any $p$. Note that the groups $\mathrm{Hom}_{\mathrm{DA_{et}}}(\mathbb{Q}[-p],i_{X}^{!}\mathbb{Q})$
are zero for any $p$. Indeed, by Verdier duality $i_{X}^{!}\mathbb{Q}=\mathbb{Q}(-1)$.
So $\mathrm{Hom}_{\mathrm{DA_{et}}}(\mathbb{Q}[-p],i_{X}^{!}\mathbb{Q})=H^{p}(X,\mathbb{Q}(-1))$
and 
\[
H^{p}(X,\mathbb{Q}(-1))=0
\]
for any smooth $X$ and any $p$%
\begin{comment}
https://mathoverflow.net/questions/309818/a-question-about-the-vanishing-of-motivic-cohomology-in-negative-tate-twist/309971\#309971
\end{comment}
. Now, use the localization sequence $Ri_{X*}i_{X}^{!}\mathbb{Q}\rightarrow\mathbb{Q}\rightarrow Rj_{X*}j_{X}^{*}\mathbb{Q}$.
$\square$ 

\section{Log motivic cohomology as a motivic sheaf}
\subsection*{Logarithmic motivic sheaves}
We will say that a presheaf  ${E:\mathrm{Sch_{fs}^{vLog}/k\longrightarrow Ch^{\bullet}(\mathbb{Q})}}$ is an \textit{effective logarithmic motivic sheaf} if it is $\mathbb{A}^{1}$-invariant, $\mathbb{P}_{log}^{1}$-invariant
and admits logarithmic $h$-descent. 
%\begin{defn}
%An \textit{effective logarithmic motivic sheaf} is a presheaf  ${E:\mathrm{Sch_{fs}^{vLog}/k\longrightarrow Ch^{\bullet}(\mathbb{Q})}}$ which is 
%\end{defn}

By result of Section 10 and Theorem \ref{'T7.4} we get
\begin{prop}
	The presheaves $R\Gamma(-,\mathbb{Q}^{\mathrm{log}}(q))$ are effective logarithmic motivic
	sheaves.
\end{prop}
Let us denote by $\mathrm{PSh(Sch_{fs}^{vLog}/k,Ch^{\bullet}(\mathbb{Q}))}$
the category of presheaves of complexes. As an ordinary theory of
motivic sheaves we can define the functor 
\[
\Omega_{\mathbb{G}_{m}}:\mathrm{PSh(Sch_{fs}^{vLog}/k,Ch^{\bullet}(\mathbb{Q}))\longrightarrow PSh(Sch_{fs}^{vLog}/k,Ch^{\bullet}(\mathbb{Q}))}
\]
by the rule 
\[
\Omega_{\mathbb{G}_{m}}E(X)=\mathrm{Cone(}E(X\times\mathbb{G}_{m})\overset{E(1)}{\longrightarrow}E(X))[-1].
\]
\begin{comment}
Note that $H^{*}(\mathbb{A}^{1}\backslash0)=\mathbb{Q}\oplus\mathbb{Q}(1)[1]$.
So we have exact triangle $pt\overset{1}{\longrightarrow}\mathbb{A}^{1}\backslash0\longrightarrow\mathbb{Q}(1)[1].$
By def $(P^{1},\infty$)=Cone($\infty\longrightarrow P^{1})$. But
$(P^{1},\infty)=G_{m}[1]$. So $G_{m}=\mathbb{Q}(1)[1].$ Finally,
$\Omega_{\mathbb{G}_{m}}E=E(Cone(1))=fib(E(1))=Cone(E(1))[-1]$.
\end{comment}

\begin{comment}
Lemma. The functor $\Omega_{\mathbb{G}_{m}}$ preserve the subcategory
of logarithmic motivic sheaves?
\end{comment}

\begin{defn}
	A \textit{logarithmic motivic sheaf} is the following data:
	
	1) For each $n\in\mathbb{Z}$ an effective log motivic sheaf $E_{n}$. 
	
	2) For each $n\in\mathbb{Z}$ a map $\varphi_{n}:E_{n}\longrightarrow\Omega_{\mathbb{G}_{m}}E_{n+1}$
	which is a quasi-isomorphism. 
\end{defn}
Let us explain how to construct the log motivic sheaf with 
\[
E_{n}=\begin{cases}
	R\Gamma(-,\mathbb{Q}^{\mathrm{log}}(n))[n], & n\geq0;\\
	0, & \mathrm{otherwise.}
\end{cases}
\]
Let us denote by $pt_{log}$ the virtualization of standard log point. Namely, this is the pair $\rm(Spec(k),k^*\to k^*\oplus \bZ \to \bZ).$ 
One can defined another endofunctor of $\mathrm{PSh(Sch_{fs}^{vLog}/k,Ch^{\bullet}(}\mathbb{Q}))$:
\[
\Omega_{pt_{log}}:\mathrm{PSh(Sch_{fs}^{vLog}/k,Ch^{\bullet}(\mathbb{Q}))\longrightarrow PSh(Sch_{fs}^{vLog}/k,Ch^{\bullet}(\mathbb{Q})})
\]
\[
\Omega_{pt_{log}}E(X)=\mathrm{Cone(}E(X\times pt_{log})\overset{E(1)}{\longrightarrow}E(X))[-1]
\]
Here the morphism $Spec(\mathrm{k})\overset{1}{\longrightarrow}pt_{log}$
correspond to the projection $pr:\mathrm{k^{*}}\oplus\mathbb{Z}\longrightarrow\mathrm{k^{*}}$.
\begin{rem}\label{R4.4}
	Let $A$ be an abelian category. Let 
	\[
	A\longrightarrow A\oplus B\longrightarrow B
	\]
	be a split exact sequence of complexes. Then $\mathrm{Cone(}\pi_{B})[-1]\equiv A\oplus\mathrm{Cone}(id_{B})[-1]$.
	So we have the inclusion $i_{A}\!:\!A\hookrightarrow\mathrm{Cone}(\pi_{B})[-1]$
	which is a quasi-isomorphism.
\end{rem}
For fix $n\geq0$ let us construct the map $$\phi_{n}\!:\!R\Gamma(-,\mathbb{Q}^{\mathrm{log}}(n))[n]\longrightarrow\Omega_{pt_{log}}R\Gamma(-,\mathbb{Q}^{\mathrm{log}}(n+1))[n+1].$$
It suffice to construct the map 
\[
\Lambda^{n}\mathrm{M_{\mathit{X}}^{gr}}\longrightarrow\mathrm{Cone(}\Lambda^{n+1}(\mathrm{M_{\mathit{X}}^{gr}}\oplus\mathbb{Q})\overset{\Lambda^{n+1}(pr)}{\longrightarrow}\Lambda^{n+1}\mathrm{M_{\mathit{X}}^{gr}})[-1]
\]
natural by $X$. %(see Definition \ref{'D5.5}). 
Note that 
\[
\Lambda^{n+1}(\mathrm{M_{\mathit{X}}^{gr}}\oplus\mathbb{Q})\simeq\underset{i+j=n+1}{\bigoplus}\Lambda^{i}(\mathrm{M_{\mathit{X}}^{gr}})\otimes\Lambda^{i}(\mathbb{Q})=\Lambda^{n+1}(\mathrm{M_{\mathit{X}}^{gr}})\oplus\Lambda^{n}(\mathrm{M_{\mathit{X}}^{gr}})
\]
and the map $\Lambda^{n+1}(pr)$ is the projection on the first term.
So $\mathrm{ker}(\Lambda^{n+1}(pr))=\Lambda^{n}\mathrm{M_{\mathit{X}}^{gr}}$
. Hence by Remark \ref{R4.4} we have the canonical inclusion
$\phi_{n}(X):\Lambda^{n}\mathrm{M_{\mathit{X}}^{gr}}\longrightarrow\mathrm{Cone}(\Lambda^{n+1}(pr))[-1]$. 
\begin{rem}
	Let $\mathrm{DA_{et}(}X,\bQ)$ be the category of rational étale motivic sheaves
	on $X$ (see \cite{'key-19} and \cite{'key-22} for the definition).
	Note that $\mathrm{DA_{et}(}X,\bQ)$ can be construct as a Bousfield localization of the category of complexes $\mathrm{Ch^{\bullet}(Sp_{\mathbb{G}_{m}}(Sh_{et}(Sm/}X,\bQ))).$
	Here $\mathrm{Sp_{\mathbb{G}_{m}}(Sh_{et}(Sm}/X,\bQ))$ is the abelian
	category of sheaves of $\mathbb{G}_{m}$-spectra (see \cite{'key-22} or \cite{'key-19}
	for details). By construction of Bousfield localization we have the
	inclusion ${\mathrm{DA_{et}(}X,\bQ)\hookrightarrow\mathrm{Ch^{\bullet}(Sp_{\mathbb{G}_{m}}(Sh_{et}(Sm}/X,\bQ)))}$
	which commute with homotopy limits. This allow us to define the functor
	\[
	\mathrm{Mor(DA_{et}(}X,\bQ))\longrightarrow\mathrm{DA_{et}(}X,\bQ);\,\,\,\,\,\,(E\overset{f}{\longrightarrow}E^{\prime})\longmapsto\mathrm{Cone(}f)[-1]
	\]
	using the standard construction of mapping cone. 
\end{rem}

Let $\mathbb{A}_{log}^{1}$ be a virtual log line - the virtualization of a unique smooth log scheme associated with the divisor $0\in \mathbb A^1.$
Now, let us consider the composition 
\[
\mathbb{G}_{m}\hookrightarrow\mathbb{A}_{log}^{1}\longrightarrow pt_{log}
\]
of the origin complement inclusion with the projection from Example
\ref{'E4.4}. For any\linebreak{}
$X\!\in\!\mathrm{Sch_{fs}^{vLog}/k}$ we have the commutative diagram
\[
\xymatrix{X\times\mathbb{G}_{m}\ar[r] & X\times\mathbb{A}_{log}^{1}\ar[d]\\
	X\ar[u]^{1}\ar[r]^{1} & X\times pt_{log}
}
\]
So the morphism gives rise to the natural transformation
\begin{equation}\label{F4.2}
	\Omega_{pt_{log}}\longrightarrow\Omega_{\mathbb{G}_{m}}.
\end{equation}
Composing this transformation with $\{\phi_{n}\}$ we get the set
of maps 
\[
\varphi_{n}:R\Gamma(-,\mathbb{Q}^{\mathrm{log}}(n))[n]\longrightarrow\Omega_{\mathbb{G}_{m}}R\Gamma(-,\mathbb{Q}^{\mathrm{log}}(n+1))[n+1].
\]

\begin{thm}\label{T4.6}
	The data $\{R\Gamma(-,\mathbb{Q}^{\mathrm{log}}(n))[n],\varphi_{n}\}$ is a logarithmic
	motivic sheaf. 
\end{thm}
Each $\varphi_{n}$ is a quasi-isomorphism by construction. So it suffice
to prove the following 
\begin{prop}
	Let $E$ be a log motivic sheaf. Then applying natural transformation (\ref{F4.2}) to $E$ we get the quasi-isomorphism
	\[
	\Omega_{pt_{log}}E\overset{\sim}{\longrightarrow}\Omega_{\mathbb{G}_{m}}E.
	\]
\end{prop}
Proof: Note that the projection $\mathbb{A}_{log}^{1}\longrightarrow pt_{log}$
admits a section $i\!:\!pt_{log}\longrightarrow\mathbb{A}_{log}^{1}$.
Applying $E$ to the log abstract blow-up square 
\[
\xymatrix{X\times pt_{log}\ar[r]^{i_{X}}\ar[d] & X\times\mathbb{A}_{log}^{1}\ar[d]\\
	X\ar[r] & X\times\mathbb{A}^{1}
}
\]
we get that $E(i_{X})$ is a quasi-isomorphism.  

On the other hand, let us consider the strict Zariski cover of $\mathbb{P}_{log}^{1}$:
\[
\mathbb{A}_{log}^{1}\amalg\mathbb{A}^{1}\longrightarrow\mathbb{P}_{log}^{1}
\]
with $\mathbb{A}_{log}^{1}\cap\mathbb{A}^{1}=\mathbb{G}_{m}$. By
strict $h-$descent for any $X$ we have the exact sequence 
\[
E(X\times\mathbb{P}_{log}^{1})\longrightarrow E(X\times\mathbb{A}_{log}^{1})\oplus E(X\times\mathbb{A}^{1})\longrightarrow E(X\times\mathbb{G}_{m}).
\]
The map $E(X\times\mathbb{P}_{log}^{1})\longrightarrow E(X\times\mathbb{A}^{1})$
is a quasi-isomorphism because $E$ is a log motivic sheaf. This implies
that $E(X\times\mathbb{A}_{log}^{1})\longrightarrow E(X\times\mathbb{G}_{m})$
is a quasi-isomorphism. $\square$

\subsection*{Restriction on a classical scheme}
Let us fix a virtual log scheme $X\in \mathrm{Sch^{vLog}_{fs}/k}$. Note that the morphism  $ X\overset{f}{\rightarrow}\underline{\mathrm{Spec(k)}} $ induces the inclusion of categories 
\begin{equation}\label{F5.2}
	\mathrm{Sm/\underline X \hookrightarrow Sch_{fs}^{vLog}/k} 	
\end{equation}
which maps
$(\underline Y \overset{g}{\rightarrow} \underline X)$ to $(\underline Y,g^*\mathcal M_X^{\mathrm{gr}})$ and $\phi: \underline Y_1 \rightarrow \underline Y_2$ to the strict morphism. Moreover, étale covers map to strict étale covers. Hence (\ref{F5.2}) is a map of sites and we have the pair of adjoint functors 
$$f^{log}_*:\mathrm{\mathcal Sh_{(\infty,1),st.h}(Sch^{vLoh}_{fs}/k,\mathbb Q})\rightleftarrows \mathrm{\mathcal Sh_{(\infty,1),et}(Sm/\underline X,\mathbb Q}):f^{log*} $$
Here $f^{log}_*F(\underline Y\overset{g}{\rightarrow} \underline X)=F(\underline Y,g^*\mathcal M^{\mathrm{gr}}_X)$ so we can naturally consider $f^{log}_*$ as the restriction. 

Now, let $E$ be an effective log motivic sheaf. By construction $f_*^{log}E$ is $\Al$-invariant. Note that effective motivic sheaf on $X$ admits étale descent (see \cite{'key-22}). So $f_*^{log}E$ belongs to $DA_{et}^{eff}(X,\bQ).$

Finally, let $\{E_n,\varphi_n\}$ be a log motivic sheaf. Note that 
$f^{log}_*$ commutes with $\Omega_{\mathbb G_m}$. So $\{f_*^{log}E_n,f_*^{log}(\varphi_n)\}$ is an object of $DA_{et}(X,\bQ).$

\subsection*{The motivic sheaf representing log motivic cohomology} Let $\{R\Gamma(-, \mathbb Q^{\mathrm{log}}(n))[n],\varphi_n\} $ be a log motivic sheaf from Theorem \ref{T4.6}. For $(X\overset{f}{\rightarrow}\mathrm k)\in\mathrm{Sch^{vLog}_{fs}/k}$ let us denote by $\bQ_X(log) $ the restriction $\{f_*^{log}R\Gamma(-, \mathbb Q^{\mathrm{log}}(n))[n],f_*^{log}(\varphi_n)\}.$ Note that by construction we have
$$ Hom_{DA_{et}(X, \mathbb Q)}(X,\bQ_X(log)(q)[p])\simeq H^{p,q}_{log}(X,\mathbb Q).$$
Now, let us give the explicit description of the motivic sheaf $\mathbb Q_X(log)$.
\begin{thm}
	Let n be the maximal rank of $\mathcal{\overline M}^{\mathrm{gr}}_{X,x}$ where $x$ runs over all closed points of $X$. Then we have the natural weak equivalence $\mathbb Q_X(log)\simeq \mathrm{S^n(\Omega_{\mathbb G_m} M^{gr}})$.
\end{thm}
\proof Firstly, let us check that $\mathrm{S^n(\Omega_{\mathbb G_m} M^{gr})\simeq \Omega_{\mathbb G_m}^n \Lambda^n(M^{gr})}$. Note that 
$\mathrm{S^n([-1])=(\Lambda^n)}[-n]$. Then it suffice to prove the following 
\begin{lem}\label{L5.5}
	$\mathrm{S^n(A(-1))\simeq (S^nA)}(-n)$ for any $A\in \mathrm{DA_{et}(\underline X,\mathbb Q)}$.
\end{lem} 
\proof By the results of Deligne (see also \cite{'key-28}, Prop. 1.4) for any $A,B\in \rm DA_{et}(\underline X,\bQ)$ we have
$$\mathrm{S^n}(A\otimes B)\simeq \bigoplus_{|\lambda |=n} \mathrm{S^{\lambda}}(A)\otimes \mathrm{S^{\lambda}}(B).$$
Using this and the isomorphisms
\begin{center}
	\begin{tikzcd}
		\mathbb Q^{\otimes n}\simeq \mathrm{S}^n \mathbb Q & \mathbb Q(1)^{\otimes n}\simeq \mathrm{S}^n (\mathbb Q(1)) & \mathbb Q(1)\otimes \mathbb Q(-1)\simeq \mathbb Q
	\end{tikzcd}
\end{center}
one can check that $\mathbb Q(-1)^{\otimes n}\simeq \mathrm{S}^n (\mathbb Q(-1))$. So we get 
$$\mathrm S^n(A\otimes \mathbb Q(-1))\simeq \bigoplus_{|\lambda |=n} \mathrm S^\lambda(A)\otimes \mathrm S^\lambda (\mathbb Q(-1))\simeq \mathrm S^nA\otimes \mathbb Q(-n).\,\,\, \square$$
Note that $\mathbb Q_X(log)$ is an omega $\Omega_{\mathbb G_m}$-spectrum with $(\mathbb Q_X(log))_{\mathrm{k}}=\Omega^\infty_{\mathbb G_m}\mathrm{\Lambda^k M^{gr}}$. For any $\Omega_{\mathbb G_m}$-spectrum $E$ we have $E_k=\mathrm{\Omega^\infty_{\mathbb G_m} \Sigma_{\mathbb G_m}^k}E$. So it suffice to check that 
$$\mathrm{\Omega^\infty_{\mathbb G_m} \Lambda^k M^{gr}\simeq \Omega^\infty_{\mathbb G_m} \Sigma_{\mathbb G_m}^k \Omega^n_{\mathbb G_m}\Lambda^n M^{gr}}.$$
For $k\leq n$ the isomorphism induce by the set of maps $\{ \varphi_l \}_{l\in \mathbb N}$ (see Theorem \ref{T4.6}). For $k>n$ it suffice to construct the set of isomorphism 
$$\mathrm{\Sigma_{\mathbb G_m}\Lambda^m M^{gr}\simeq \Lambda^{m+1} M^{gr}}$$
for any $m\geq n$. Let 
\begin{equation}\label{F5.3}
	\mathrm{\Lambda^m M^{gr}(1)[1]\overset{\tau_m}{\longrightarrow} \Lambda^{m+1} M^{gr}}
\end{equation}
be the map induce by the $\Sigma_m$-equivariant map 
$$ \mathrm{(M^{gr})^{\otimes m}\otimes \mathbb Q(1)[1]\overset{id\otimes \alpha}{\longrightarrow}(M^{gr})^{\otimes(m+1)}}_.$$
Observe that $\tau_m$ is exactly the rightmost differential in the Koszul complex $\mathrm{\Lambda^{m+1}(Cone(\mathbb Q(1)[1] \rightarrow M^{gr}))}$. So we have 
$$\mathrm{Cone(\tau_m)\simeq \Sigma^{\infty}_{\mathbb G_m} Ex_{\mathbb Q,\underline X}\Lambda^{m+1} (\mathcal{\overline M}^{gr}_X)}$$
(see the proof of Corollary \ref{'C7.3} for more details). 
But $\mathrm{\Lambda^{m+1} (\mathcal{\overline M}^{gr}_X)}=0$ for $m+1>n$.$\square$

%Recall that in [Sh22? was constructed the functor $[-]^{log}:Sch^{vLog}_{fs}/k\ra DM_{gm}(k,\bQ)$ which maps a virtual log scheme $(X,\mathcal M^{gr}_X)$ to the motive  $[X]^{log}.$ The following Corollary explain how this notion relates with the above definition of log motivic cohomology.
\begin{cor}\label{''C14.3}
	We have the natural weak equivalence 
	$\bQ_X(log)\iso \mathrm{Sym}^*(\mathrm{M^{gr}}(-1)[-1]).$
	In particular, 
	$$H^{p,q}_{log}(X,\bQ)\iso Hom_{DM(k,\bQ)}([X]^{log},\bQ(q)[p]).$$
\end{cor}
\proof By Theorem \ref{T2.6} and Lemma \ref{L2.7} we have
$\mathrm{Sym^*(\mathrm{\Omega_{\mathbb G_m}M^{gr}})\simeq colim_kS^k(\mathrm{\Omega_{\mathbb G_m}M^{gr}})}.$
Then, using Lemma \ref{L5.5} and the isomorphisms \ref{F5.3}, we conclude that $\mathrm{colim_kS^k(\mathrm{\Omega_{\mathbb G_m}M^{gr}})\simeq S^n(\mathrm{\Omega_{\mathbb G_m}M^{gr}})}$ where $n$ is the maximal rank of $\mathcal{ \overline M}^{\mathrm{gr}}_{X,x}$. $\square$ 

\section{Log smooth schemes and motivic tubular neighborhoods}

\subsection*{Motivic cohomology of log smooth schemes}
Let $X$ be a log smooth scheme over $k\subset \bC$. Denote by $X^*$ the maximal open subset with trivial log structure. Suppose that the underlying scheme $\underline X$ is smooth. Then the Kato-Nakayama space $X^{log}$ can be constructed as so-called \textit{oriented real blow-up} (see \cite{'key-17} for more details) .  In particular, $X^{log}$ is a smooth manifold with boundary and the interior of $X^{log}$ coincides with $X^*$.  Consequently the natural inclusion
$$X^*\hookrightarrow X^{log}$$
is a homotopy equivalence.
\\Hence it would be natural to expect that the logarithmic motivic cohomology of $X^*$ and $X$ are isomorphic. In fact, the stronger statement is true.
\begin{thm}\label{'T8.1}
	Let k be an arbitrary field. Let $\mathcal{F}$ be an effective log motivic sheaf. Then the induced map 
	$$\mathcal F(X)\to \mathcal F(X^*)$$
	is a quasi-isomorphism for any smooth and log smooth $X$ over k.  	
\end{thm} 
Further, we will use the following description of log smooth
schemes.
\begin{lem}
	(Corollary 4.8 of \cite{'key-26}) Let k and $X$ be as above. Then, locally in
	étale topology, there exist charts $f_{i}:X\rightarrow\mathbb{A}^{n_{i}}\times\mathbb{A}_{log}^{m_{i}}$
	with étale $f_{i}$ . So the log structure $\mathcal{M}_{X}$
	on $\underline{X}$ is the canonical log structure associated with
	some normal crossing divisor $D\subset\underline{X}$.
\end{lem}
\begin{rem}\label{'R8.3}
	Let $\mathcal F$ be as above. Let
	\begin{equation}\label{'F8.1}
		\mathbb{G}_{m}^n\hookrightarrow\mathbb{A}_{log}^{n}
	\end{equation} be
	the strict open immersion. Then for any $Y\in\mathrm{Sch_{fs}^{vLog}/k}$
	the map 
	\[
	\mathcal F(Y\times\mathbb{A}_{log}^{n})\longrightarrow \mathcal F(Y\times\mathbb{G}_{m}^{n})
	\]
	is a quasi-isomorphism. Indeed, let us consider the strict Zariski cover of $\mathbb{P}_{log}^{1}$:
	\[
	\mathbb{A}_{log}^{1}\amalg\mathbb{A}^{1}\longrightarrow\mathbb{P}_{log}^{1}
	\]
	with $\mathbb{A}_{log}^{1}\cap\mathbb{A}^{1}=\mathbb{G}_{m}$. By
	strict $h-$descent for any $Y$ we have the exact sequence 
	\[
	\mathcal F(Y\times\mathbb{P}_{log}^{1})\longrightarrow \mathcal F(Y\times\mathbb{A}_{log}^{1})\oplus \mathcal F(Y\times\mathbb{A}^{1})\longrightarrow \mathcal F(Y\times\mathbb{G}_{m}).
	\]
	The map $\mathcal F(Y\times\mathbb{P}_{log}^{1})\longrightarrow \mathcal F(Y\times\mathbb{A}^{1})$ is a quasi-isomorphism. This implies
	that ${\mathcal F(Y\times\mathbb{A}_{log}^{1})\longrightarrow \mathcal F(Y\times\mathbb{G}_{m})}$
	is a quasi-isomorphism. 
	
	So we checked the statement for n=1. But (\ref{'F8.1}) is the product of $\Gm\hookrightarrow \Al_{log}$ and we have the chain of quasi-isomorphisms
	$$\mathcal F(Y\times \Gm^n)\to \mathcal F (Y\times {\Gm}^{n-1}\times \Al_{log})\to...\to \mathcal F(Y\times {\Gm}^{n-k}\times \mathbb A^k_{log})\to...\to \mathcal F(Y\times \mathbb A^n_{log}).$$
\end{rem}
Now, let $\mathbb{A}_{log}^{n}\setminus0\hookrightarrow\mathbb{A}_{log}^{n}$
be the strict open immersion. 
\begin{lem}\label{'L8.4}
	The map $\mathcal F(Y\times\mathbb{A}_{log}^{n})\longrightarrow \mathcal F(Y\times \mathbb{A}_{log}^{n}\setminus0)$ is a quasi-isomorphism.
\end{lem}
\textit{Proof:} By Remark \ref{'R8.3} the composition $\mathbb{G}_{m}^{n}\hookrightarrow\mathbb{A}_{log}^{n}\setminus0\hookrightarrow\mathbb{A}_{log}^{n}$
becomes a quasi-isomorphism after applying $\mathcal F(Y\times -)$. So it enough to check
that  
\begin{equation}\label{'F8.2}
	\mathcal F(Y\times\mathbb{A}_{log}^{n}\setminus0)\rightarrow \mathcal F(Y\times \mathbb{G}_{m}^{n})
\end{equation}
is a quasi-isomorphism. For $n=1$, (\ref{'F8.2}) is identity. Let $n > 1$. It suffice to show that there is Zariski covers $\bigsqcup^n_{i=1} U_i\to \Gm^n$ and $\bigsqcup^n_{i=1} V_i\to \mathbb A^n_{log}\setminus 0$ such that for any index $i_1,...,i_k$ with $k\leq n$ we have 
$$\mathcal F(Y\times(V_{i_1}\cap ... \cap V_{i_k}))\iso  \mathcal F(Y\times(U_{i_1}\cap ... \cap U_{i_k})).$$
Now, let $V_i=\mathbb A^n\setminus \{ x_i=0\}$, $U_i=\Gm^n$ and use induction by $n$. $\square$
\begin{lem}\label{'L8.5}
	Let $X$ be a smooth and log smooth scheme with canonical divisor ${D=D_1\cup...\cup D_r}$. Let us denote by $Z$ the intersection $D_1\cap...\cap D_r.$ Then the induced map
	\begin{equation}\label{'F8.3}
		\mathcal F(X)\to \mathcal F(X\setminus Z)
	\end{equation} 
	is a quasi-isomorphism.  
\end{lem}
\textit{Proof:} The question is local so we may assume that $X$ admits an étale  chart ${f:X\to \mathbb A^{n-r}\times \mathbb A^r_{log}.}$ So $Z=f^{-1}(0)$. In this case, Morel and Voevodsky \cite{'key-30} construct the common Nisnevich neighborhood of $Z\subset X$ and the zero section $0_Z$ of $N_XZ.$  Let us denote such a neighborhood by $U$ and define the log structure on $U$ as the pullback along the map $U\to X \overset{f}{\to} \mathbb A^{n-r}\times \mathbb A^r_{log}. $ Then we get the commutative diagram 
\[
\xymatrix{ & U\ar@/_{0.3pc}/[ld]_{pr_{2}}\ar@/^{0.3pc}/[rd]^{pr_{1}}\\
	Z\!\times\!\mathbb{A}_{log}^{r}\ar@/_{0.3pc}/[rd]^{f|_{Z}\!\times\!\mathrm{id}} &  &  X \,\,\,\,\,\ar@/^{0.3pc}/[ld]_{f}\\
	& \mathbb{A}^{n-r}\!\times\!\mathbb{A}_{log}^{r}
}
\]
with strict arrows, $pr_{1}^{-1}(Z)\simeq Z$ and $pr_{2}^{-1}(0_Z)\simeq Z$.
Consequently, we have two log Nisnevich squares
\[
\xymatrix{Z\!\times\!\mathbb{A}_{log}^{r}\setminus0\ar[d] & U\!\setminus\!Z\ar[r]\ar[l]\ar[d] & X\!\setminus\!Z\ar[d]\\
	Z\!\times\!\mathbb{A}_{log}^{r} & U\ar[r]\ar[l] & X
}
.
\]
Applying $\mathcal F$ and using the Lemma \ref{'L8.4} we see that the map $\mathcal F(U)\to \mathcal F(U\setminus Z)$ is a quasi-isomorphism. So the map (\ref{'F8.3}) is also a quasi-isomorphism. $\square$ 

\textit{Proof of Theorem \ref{'T8.1}: }
\\Let $D=D_1\cup...\cup D_r$ be a canonical divisor of $X$ and $r$ be a number of smooth components. For $r=1$ the statement follows from Lemma  \ref{'L8.5}. Let $r>1$. By Lemma \ref{'L8.5} it enough to prove that the map 
$$\mathcal F(X\setminus Z)\to \mathcal F(X^*)$$
is a quasi-isomorphism. So it suffice to find Zariski covers $\bigsqcup^n_{i=1} U_i\to X^*$ and $\bigsqcup^n_{i=1} V_i\to X\setminus Z$ such that for any index $i_1,...,i_k$ with $k\leq n$ we have 
$$\mathcal F(V_{i_1}\cap ... \cap V_{i_k})\iso  \mathcal F(U_{i_1}\cap ... \cap U_{i_k}).$$
Now, let $V_i=X\setminus D_i$, $U_i=X^*$ and use induction by $r.$ $\square$
\begin{cor}\label{''C15.6}
	Let $X$ be a log smooth scheme over $k.$ Suppose that the underlying scheme $\underline X$ is smooth. Then the natural map
	$$[X^*]\to [X]^{log}$$
	is a weak equivalence.
\end{cor}
\proof Let us use Yoneda Lemma. Note that by construction for any smooth $Y\in \rm Sm/k$ and any log scheme $(Z,\mathcal M^{gr}_Z)$ we have 
\begin{equation}\label{F7.3}
	\mathrm{Hom(\mathit{Y(-q)[-p]},\mathbb D([Z]^{log}))\simeq}H^{p,q}_{log}(Y\times Z,\mathbb Q).
\end{equation}
But $H^{p,q}_{log}(Y\times X,\bQ)\iso H^{p,q}(Y\times X^*,\bQ)$ by Theorem \ref{'T8.1} (let $\mathcal F:=H^{p,q}_{log}(Y\times -,\bQ)$). $\square$ 

\begin{rem}\label{''R15.7}
	Using Yoneda lemma and Corollary \ref{''C14.3} one can check that the functor $[-]^{log}$ admits descent under log abstract blow-up. That is, for any log abstract blow-up square 
	$$\begin{tikzcd}
		Z' \arrow[d] \arrow[r] & X' \arrow[d] \\
		Z \arrow[r]            & X           
	\end{tikzcd}$$
	we have the exact triangle 
	$$[Z']^{log}\to [Z]^{log}\oplus [X']^{log}\to [X]^{log}.$$
\end{rem}

\subsection*{Motivic tubular neighborhood}

Let $\underline X$ be a smooth scheme, $\underline D$ be a normal crossing divisor and $\underline U=\underline X\setminus \underline D$. Let us denote by $X$ the associated log scheme. Let $D:=(\underline D,\mathcal M_X|_{\underline D}) $. Note that the square
\begin{equation}\label{F8.1}
	\begin{tikzcd}
		D \arrow[r] \arrow[d]  & X \arrow[d]  \\
		\underline D \arrow[r] & \underline X
	\end{tikzcd}
\end{equation}
is a log abstract blow-up square. So by Remark \ref{''R15.7} and Corollary \ref{''C15.6} we have the exact triangle
$$[D]^{log}\longrightarrow [\underline D]\oplus[\underline U]\ra [\underline X].$$
This allows us to interpret the motive $\mathrm{M^{log}}(D)$  as the \textit{motivic punctured tubular neighborhoods of $\underline D$ in $\underline X$}. Let us denote the punctured tubular neighborhoods by $\mathrm{PTN}_{\underline X} \underline D$. Note that we also have the canonical 
\textit{exponential  map} 
\begin{equation}\label{key}
	exp:\mathrm{PTN}_{\underline X}\underline D\longrightarrow \underline U
\end{equation}
which induce by the strict closed immersion of log schemes $D\hookrightarrow X.$
%\\Remark. Note that in the paper [Levine]  a map
%$$\mathrm{PTN}_{\underline X}\underline D\longrightarrow \underline U$$
%was also constructed. We do not know how this is exactly related to our construction. Nevertheless, one can note that the map differ by an isomorphism in triangulated category $\mathrm{DM_{gm}(k,\mathbb Q)}$.
\begin{rem}
	Suppose $\mathrm k$ admits resolution of singularities. Let $\underline Y$ be a smooth scheme and $\underline Z$ be an arbitrary closed subscheme. Then by Hironaka  principalization theorem (\cite{H64}, see also Theorem 3.35 of \cite{K09}) there exists the abstract blow-up
	\begin{equation}\label{F8.3}
		\begin{tikzcd}
			\underline D \arrow[r] \arrow[d] & \underline X \arrow[d, "p"] \\
			\underline Z \arrow[r]                      & \underline Y                     
		\end{tikzcd}
	\end{equation}
	where $\underline X$ and $\underline D$ as above and $p^{-1}(\underline Y\setminus \underline Z)\iso\underline X\setminus \underline D$.
	Then combining \ref{F8.1} and \ref{F8.3} we get the $log$ abstract blow-up
	$$\begin{tikzcd}
		D \arrow[r] \arrow[d] & X \arrow[d] \\
		\underline Z \arrow[r]           & \underline Y          
	\end{tikzcd}$$
	and, consequently, the exact triangle
	$$[D]^{log}\longrightarrow [\underline Z]\oplus [\underline Y\setminus \underline Z]\longrightarrow [\underline Y].$$
	This allows us to interpret $[D]^{log}$  as the motivic punctured tubular neighborhoods of $\underline Z$ in $\underline Y$. 
	%Note that in this case we also have the canonical exponential map 
	%$$[D]^{log}\longrightarrow [\underline Y\setminus \underline Z].$$
\end{rem}

\section{Limit motives via log geometry}

\subsection*{Unipotent local systems}
Let $\mathring \Delta:=\{z\in\bC\,|\,0<|z|\leq 1 \}$ be a punctured disk. Let us denote by $DLocSys(\mathring \Delta)_\bQ$ the full subcategory of $D(Sh(\disk,\bQ))$ containing the complexes $K^\bullet$ such that $\mathcal H^i(K^\bullet)$ are local systems. Let $\uni$ be the full subcategory of $\dloc$ generated by the trivial local system $\underline \bQ$. By abstract nonsense we have 
\begin{equation}\label{F8.4}
	\uni \iso \rm A-Mod
\end{equation}
where $\rm A=\it REnd(\underline \bQ)\iso C^*(S^{\rm 1},\bQ)$ is the $\mathbb E_\infty$-algebra of cochains.
Indeed, for any derived local system $M^\bullet$ we can constructed the $C^*(S^1,\bQ)$-module $RHom(\underline Q,M^\bullet)$. On the other hand, let ${\bf B}C^*(S^1,\bQ)$ be the category with one object and the ring of endomorphisms $C^*(S^1,\bQ)$. Then the inclusion ${\bf B}C^*(S^1,\bQ)$ gives rise to the Yoneda extension 
$${\rm A-Mod\ra \uni.}$$ 
One can check that this functors form the equivalence of categories.
\\Let $i:pt \hookrightarrow \disk$ be the inclusion of a point. Then we have the inverse image functor 
$$i^*:\dloc \ra D(Sh(pt,\bQ))\iso D(\bQ-Vect)$$ 
and we can restrict $i^*$ on $\uni.$ On the other hand, the map $i$ induces the map of $\mathbb E_\infty$-algebras 
\begin{equation}\label{F8.5}
	C^*(S^1,\bQ)\ra C^*(pt,\bQ)\iso \bQ.
\end{equation}
The map \ref{F8.5} defines the functor 
$$-\otimes_{C^*(S^1,\bQ)}C^*(pt,\bQ):{\rm A-Mod}\ra D(\bQ-Vect)$$
Note that $-\otimes_{C^*(S^1,\bQ)}C^*(pt,\bQ)$ and $i^*$ commute with colimits. The both maps 
$$REnd_{\dloc}(\underline Q)\ra End_\bQ(\bQ)$$
$$REnd_\mathrm A(\mathrm A)\ra End_\bQ(\bQ)$$
coincide with \ref{F8.5}. Moreover, this maps differ only by the isomorphism 
$$REnd_{\dloc}(\underline Q)\iso REnd_\mathrm A(\mathrm A)$$
induced by the equivalence \ref{F8.4}. So $i^*$ and $-\otimes_{C^*(S^1,\bQ)}C^*(pt,\bQ)$ are isomorphic.
\subsection*{Nearby cycles as a bar-construction}
Let $C$ be a smooth complex curve with marked point 0. Let $f:X\ra C$ be a semi-stable degeneration with a special fiber $X_s$. Suppose that $f$ is proper and let $\mathring f:\mathring X\ra \mathring C$ be the restriction of $f$ on $\mathring X:=X\setminus X_s.$ Then by Ehresmann's theorem $\mathring f$ is a fiber bundle. Let $t$ be a point of $\mathring C$. Using proper base change we conclude that $R^n\mathring f_*\bQ$ is a local system with the fiber $H^n(X_t,\bQ).$
Now, let us chose a disk $\Delta\subset C$ centered at $s$. For $X':=f^{-1}(\Delta)$ the map $f:X'\ra \Delta$ is also proper semi-stable degeneration. Note the $X'$ can be considered as a tubular neighborhood of $X_s$ in $X$. By the previous the complex $R\mathring f_*\bQ$ belongs to $\dloc.$
\begin{lem}
	$R\mathring f_*\bQ$ is a unipotent derived local system.
\end{lem}
\proof Note that $R^m\mathring f_*\bQ=H^n(X_t,\bQ)$ is zero for sufficiently large $m$. So it enough to check that $R^b\mathring f_*\bQ$ belongs to $\uni$. Let us denote by $\gamma$ the monodromy matrix acts on $H^n(X_t,\bQ)$. Note that $\gamma$ is a unipotent operator. Now, let $\rm M$ be an arbitrary local system with unipotent monodromy $\gamma$. Then $\rm M$ belongs to $\uni$. Indeed, $(\gamma -1)^k=0$ so all eigenvalues are 1. Representing $\rm M$ as a direct sum we can assume that $\gamma$ is a Jordan block. Then we can use the induction by dimension of $\gamma$. $\square$
%if dim=1 then $M=\underline \bQ$. Suppose we done for dim=k and dim M=k+1. Then \Gamma e_1=e_1. So we have the equivariant inclusion $\underline \bQ \hookrightarrow M$ and \Gamma acts on $M/\underline \bQ$ as the Jordan block of dim=k.

Note that the equivalence $\uni \ra \rm A-Mod$ maps $R\mathring f_*\bQ$ to the complex 
$$RHom_{D(Sh(\disk,\bQ))}(\bQ,R\mathring f_*\bQ)\iso RHom_{D(Sh(\mathring X',\bQ))}(\bQ,\bQ)\iso C^*(\mathring X',\bQ)$$
such that the structure of $C^*(S^1,\bQ)$-module on $R\mathring f_*\bQ$ induces by the homomorphism of $\mathbb E_\infty$-algebras 
$$C^*(\disk,\bQ)\sa{\mathring f^*} C^*(\mathring X',\bQ).$$
So for $f:X\ra C$ and $t:pt\hookrightarrow \disk \hookrightarrow \mathring C$ we have 
\begin{equation}\label{F8.6}
	t^*(R\mathring f_*\bQ)\iso C^*(\mathring X',\bQ)\otimes_{C^*(S^1,\bQ)}C^*(t,\bQ)
\end{equation}
where $\mathring X'$ is a punctured tubular neighborhood of $X_0\subset X.$
\\Now, let $\psi_f: D(Sh(\mathring X,\bQ))\ra D(Sh(X_s,\bQ))$ be the nearby cycles functor. For semi-stable degeneration $f$ we have 
$$H^*(X_s,\psi_f\bQ)\iso H^*(X_t,\bQ).$$
So, using \ref{F8.6}, we get 
\begin{equation}\label{F8.7}
	R\Gamma(X_s,\psi_f\bQ)\iso C^*(\mathring X',\bQ)\otimes_{C^*(S^1,\bQ)}C^*(t,\bQ).
\end{equation}  

\subsection*{The limit motive of a semi-stable degeneration}
Let $k$ be an arbitrary field. Again, let $C$ be a smooth curve over $k$ and $f:X\ra C$ be a semi-stable degeneration. Let us associated with $(X,X_s)$ the canonical log structure. Let $\mathcal M^{gr}_{X_s}$ be the restriction of this structure on $X_s.$  Note that the map $f$ induces the morphism of virtual logarithmic schemes
\begin{equation}\label{F8.8}
	f:(X_s,\mathcal M^{gr}_{X_s})\ra pt_{log}.
\end{equation} 
By \ref{T7.5}, the map $[X_s]^{log}\ra [pt_{log}]^{log}=[\Gm]$ is a homomorphism of $\mathbb E_\infty$-coalgebras. Let $t:\rm Spec(k)\hookrightarrow \Gm$ be a closed point of $\Gm$. Then, replicating \ref{F8.6}, we define the geometrical motive 
$${\rm LM}_f:=([X_0]^{log})^*\otimes_{[\Gm]^*}[Spec(k)]. $$
We will call the motive ${\rm LM}_f$ \textit{the Voevodsky limit motive of f}.

Now, let us fix an inclusion $\tau:k\to \bC.$ Let $\rm Betti: DM_{gm}(k,\bQ)\to D(\bQ-Vect)$ be the associated Betti realization.  
\begin{prop}
	There exists the canonical quasi-isomorphism 
	$${\rm Betti( LM}_f)\iso R\Gamma(X_s,\psi_f\bQ)$$
	induced by \ref{F8.6}.
\end{prop}
\proof Let $X'$ and $\mathring X'$ be as above. Note that the inclusion $X_s^{an}\hookrightarrow X'$ induces the canonical log analytic structure on $X'$. So we have the commutative diagram 
$$
\begin{tikzcd}
	{C^*(X_s^{log},\bQ)} \arrow[r]                  & {C^*(X'^{log},\bQ)}               & {C^*(\mathring X',\bQ)} \arrow[l]               \\
	{C^*((pt_{log})^{log},\bQ)} \arrow[r] \arrow[u] & {C^*(\Delta^{log},\bQ)} \arrow[u] & {C^*(\mathring \Delta,\bQ)} \arrow[l] \arrow[u]
\end{tikzcd}
$$
where each map is a homomorphism of $\mathbb E_\infty$-algebras. Combining Corollary 4.17. of \cite{PS08} with Theorem \ref{R_pi_Q} we conclude that the right-upper and the right-lower horizontal arrows are quasi-isomorphisms.

On the other hand, one can check that the functor $X\mapsto C^*(X^{log},\bQ)$  admits descent under virtual logarithmic blow-ups. Indeed, using Theorem \ref{R_pi_Q}, we can check this by the same way as in the Section 10. So we get the exact triangle 
$$C^*(X',\bQ)\to C ^*(X'^{log},\bQ)\oplus C^*(X_0^{an},\bQ)\to C^*(X_0^{log},\bQ).$$
But the inclusion $X_0^{an}\to X'$ is a homotopy equivalence. So the left-upper and the left-lower horizontal arrows are quasi-isomorphisms.

It remains to note that, by Proposition \ref{T6.6}, $C^*((pt_{lpg})^{log},\bQ)\ra C^*(X_0^{log},\bQ)$ is the (covariant) Betti realization of the map $[\Gm]^*\to [X_0]^{log*}$ % which induced by \ref{F8.8}
.$\square$

\section{Motivic monodromy filtration}

\subsection*{Limit Hodge structure}

Let $\Delta$ be a complex unit disk and $f:X\to \Delta$ be a proper semi-stable degeneration. By Ehresmann's theorem the restriction of $f$ on $X-X_0$ is a fiber bundle over $\mathring \Delta$. So $f$ defines the variation of Hodge structures over $\mathring \Delta$ Namely, let us fix $t\in \mathring \Delta$ and $n\in \mathbb N$. We will denote by $\check{D}$ the \textit{completed local period domain}. Recall that $\check{D}$ is the space of flags $...\subset F^p \subset ... \subset F^1\subset F^0=H^n(X_t,\bC)$ with the property ${\rm dim}(F^p/F^{p-1})=h^{p,q}(X_t).$ Note that $f$ is proper so $Rf_*\bQ$ defines the local system over $\mathring\Delta$ with the fibers $H^n(X_t,\bQ).$ Then we can define the \textit{local period map} $\Phi: \mathring \Delta \rightarrow \check D$
$$t'\in \mathring \Delta \mapsto \{
{\rm filtration \, on \, \mathit{H^n(X_t,\bC)} \, induced \, by \,}
{\rm the \, Hodge \, filtration \, on \,} H^n(X_{t'},\bC). \} $$  
Note that an isomorphism $H^n(X_t,\bC)\iso H^n(X_{t'},\bC)$ depends on the choosing of a path so the map $\Phi$ is multi-valued. However, one can ask is it possible to continue $\Phi$ holomorphically on all $\Delta$?

Let $\gamma$ ba a monodromy matrix which act on  $H^n(X_t,\bC)$ and $N=-{\rm log}(\gamma)$. Then we can define $\tilde \Psi(z):=exp (zN)\Phi(z)$ where $z$ is the coordinate on the universal cover of $\mathring \Delta.$ Note that $\tilde\Psi(z)=\tilde \Psi(z+1)$ so we have well-defined map $\Psi: \mathring \Delta\to \check D$. Then $\Psi$ can be holomorphically continue on all $\Delta.$

Observe that the extension of $\Psi$ defines the new filtration $F_{\rm lim}^\bullet$ on $H^n(X_t,\bC)$ given by $\Psi(0).$ On the other hand, by Borel results, the monodromy matrix $\gamma$ is unipotent. So there is some $k$ such that $N^k=0$ and $N$ defines one more filtration 
$$...\subset 0\subset W_0\subset ...\subset W_{2k}=H^n(X_t,\bQ)$$ 
with the properties:
\vspace{-\topsep} 
\begin{itemize}
	\item[-] for each $r$ we have $N(W_r)\subset W_{r-2}$;
	\item[-] the operator $N$ defines the isomorphisms ${\rm Gr}^W_r(H^n(X_t,\bQ))\iso {\rm Gr}^W_{r-2}(H^n(X_t,\bQ)). $
\end{itemize}
\vspace{-\topsep}
Moreover, Schmid proved \cite{Sch73} that the data 
$(H^n(X_t,\bZ),W_\bullet, F^\bullet_{\rm lim})$
defines the mixed Hodge structure. This Hodge structure called \textit{the limit Hodge structure}. 

\subsection*{The Steenbrink construction}
There is another, pure algebraic construction of the limit Hodge structure proposed by Steenbrink (see \cite{St76}, \cite{St95} and \cite{PS08}). Namely, he constructed two filtered bicomplexes $(C^{\bullet,\bullet},W(M)^{\bullet,\bullet})_r)$  and $(A^{\bullet,\bullet},W(M)^{\bullet,\bullet}_{\bC, r})$ with $C^{\bullet,\bullet}\in {\rm Bicomp}(Sh(X_0,\bQ))$ and $A^{\bullet,\bullet}\in {\rm Bicomp}(Sh(X_0,\bC))$, the filtered quasi-isomorphism 
$$\mathrm{Tot}(C^{\bullet,\bullet}\otimes \bC)\iso \mathrm{Tot}(A^{\bullet,\bullet})$$ 
and the Hodge filtration on $\mathrm{Tot}(A^{\bullet,\bullet})$ given by the partial totalizations 
$$F^pA={\rm Tot}(s_{\geq p}A^{\bullet,\bullet});\,\,\,s_{\geq r}A^{i,j}=\begin{cases}
	A^{i,j} & ,\, j\geq p\\
	0 & \rm ,\, otherwise.
\end{cases} $$
Steenbrink proved that the data 
\begin{equation}\label{F9.1}
	\psi_f\bZ^{\rm Hdg}:=(\psi_f\bZ,({\rm Tot}(C^{\bullet,\bullet},{\rm Tot}(W(M)^{\bullet,\bullet}),({\rm Tot}(A^{\bullet,\bullet}),{\rm Tot}(W(M)_\bC^{\bullet,\bullet}),FA))
\end{equation}
forms the mixed Hodge complex of sheaves over $X_0$. %Here $\psi_f:D(Sh(X_t,\bZ))\to D(Sh(X_0,\bZ))$ is the nearby cycles functor. 
Then the mixed Hodge structure on $\mathbb H^n(X_0,\psi_f\bZ^{\rm Hdg})$ is isomorphic to the Schmid limit Hodge structure on $H^n(X_t,\bQ)$.

Let us describe the bicomplexes $C^{\bullet,\bullet}$ and $W(M)^{\bullet,\bullet}_r$ explicitly. Note that we have the canonical log structure on $X$. Let $\mathcal M\to \cO_{X_0}$ be the restriction of this structure on $X_0.$ Then we have the associated exponential complex $exp$ given as the cone of the map 
$$\cO_{X_0^{an}}\sa{exp}\mathcal M^{gr}_{X_0^{an}}.$$
Moreover, $f$ induced the morphism of smooth log analytic spaces $(X^{an},X_0^{an})\to (\Delta,0).$ The restriction of this morphism on $X_0$ gives the morphism of virtual log analytic spaces 
\\${(X_0^{an},\cO^*_{an}\hookrightarrow \cM_{an})\to pt_{log}}$. So we get the element $t\in \Gamma(X_0^{an},\cM_{an}).$ 
\\Let us define the following bicomplexes: 
\\the bicomplex $W^{\bullet,\bullet}_r$
$$...\ra \underset{-1}{0} \ra \underset{0}{S^r(exp)[1]} \sa{\wedge t} \underset{1}{S^{r+1}(exp)[2]}\sa{\wedge t} \underset{2}{S^{r+2}(exp)[3]} \sa{\wedge t} ... $$
where the horizontal differential given by
\begin{equation}\label{F9.2}
	S^i(\cO_{an})\otimes \Lambda^j(\cM_{an})\sa{{\rm id}\otimes \wedge t} S^i(\cO_{an})\otimes \Lambda^{j+1}(\cM_{an})
\end{equation} 
and $W_r^{p,q}=0$ for $p<0;$ 
\\the bicomplex $\tilde W(M)^{\bullet,\bullet}_r$
$$...\ra \underset{-1}{0} \ra \underset{0}{S^{r+1}(exp)[1]} \sa{d} \underset{1}{S^{r+3}(exp)[2]}\sa{d} \underset{2}{S^{r+5}(exp)[3]} \sa{d} ... $$
where the differential $d$ given by the composition of \ref{F9.2} with the inclusion
\begin{equation}\label{F9.3}
	\bQ\otimes S^i(exp)\hookrightarrow S^{i+1}(exp)
\end{equation}
$$ a\otimes w \mapsto a\cdot w.$$
Here $\tilde W(M)^{p,\bullet}_r=S^{r+2p+1}(exp)[p]$ and $\tilde W(M)^{p,q}_r=0$ for for $p<0.$
\\Note that we can also define $W(M)^{\bullet,\bullet}_r$ for $r<-1$ if we put
$$S^k(-)=0\,\,\,{\rm for}\,\,\,k<0.$$ Let $\tilde C^{\bullet,\bullet}$ be the bicomplex
$$...\ra \underset{-1}{0} \ra \underset{0}{\mathrm{Sym}^*(exp)[1]} \sa{\wedge t} \underset{1}{\mathrm{Sym}^*(exp)[2]}\sa{\wedge t} \underset{2}{\mathrm{Sym}^*(exp)[3]} \sa{\wedge t} ...\,\,.$$
Note that we have the canonical inclusions 
$$W^{\bullet,\bullet}_r\hookrightarrow \tilde W(M)^{\bullet,\bullet}_{r+j}$$
$$\tilde W(M)^{\bullet,\bullet}_r\hookrightarrow \tilde W(M)^{\bullet,\bullet}_{r+j}$$
$$W^{\bullet,\bullet}_r\hookrightarrow \tilde C^{\bullet,\bullet}$$
given by compositions of \ref{F9.3}. Let us define $$C^{\bullet,\bullet}:=\tilde C/W_0$$
$$W(M)_r^{\bullet,\bullet}:={\rm im}(\tilde W(M)_{r}^{\bullet,\bullet}\to C^{\bullet,\bullet}).$$ 
%Observe that up to renaming this definition coincides with the definition from Section ? of [?]. 
%Steenbrink compute (see ?) that the graded pieces ${\rm Gr}^{W(N)}_r$ is isomorphic to $$  

Now, we should  make a few remarks. Let $N$ be the number of irreducible components of $X_0.$ Note that $\tilde W(M)^n_r\subseteq \ W_0^n$ if $r+2n+1\leq n.$ So $\tilde W(M)^n_r=0$ for $n\leq -r-1.$ Let $n>-r$ and $r<-(N-1)$. Then $r+2n+1>N$ and 
$$({\rm Gr}^{W(M)}_r(C^{\bullet,\bullet}))^n=\Lambda^{r+2n+1}\bar{\mathcal M}^{gr}[-r-n]=0.$$ 
On the other hand, let $n=-r$ and $r<-(N-1).$ Then $-r+1>N$ and $({\rm Gr}^{W(M)}_r(C^{\bullet,\bullet}))^n=\Lambda^{-r+1}\bar{\mathcal M}^{gr}[-r-n]=0.$ Finally, assume that $r>N-1$. Then again $({\rm Gr}^{W(M)}_r(C^{\bullet,\bullet}))^n=0.$ So we get
\begin{lem}\label{L9.1}
	The inclusions $W(M)_{r-1}\subset W(M)_r$ are quasi-isomorphism for $r\leq -(N-1)$ and $r\geq N-1.$ Moreover, for $r\geq N-1$ the canonical inclusions $W(M)^{\bullet,\bullet}_r\hookrightarrow C^{\bullet,\bullet}$ are quasi-isomorphisms.
\end{lem}
\proof It remains to check that the canonical inclusion 
$\tilde W(M)_{N-1}\hookrightarrow \tilde C$ is a quasi-isomorphism. Using the spectral sequence of a bicomplex, it can be reduce to the fact that the maps 
$S^{N-1+2n+1}(exp)\ra \mathrm{Sym}^*(exp)$ 
are quasi-isomorphisms. $\square$
\vspace{3mm}

Let $s_{\geq N-1}C^{\bullet,\bullet}$ be the complex with $s_{\geq N-1}C^{n,*}=C^{n,*}$ for $n\geq N-1$  and $s_{\geq N-1}C^{n,*}=0$ otherwise. Using the same arguments as in the proof of the lemma one can check that the canonical inclusion $s_{\geq N-1}C^{\bullet,\bullet}\ra C^{\bullet,\bullet}$ is a quasi-isomorphism. Indeed, it suffice th check that the complexes $W(M)_{N-1}^n$ us quasi-isomorphic to zero for $n>N-1$. But 
$$W(M)_{N-1}^n=\mathrm{Cone}(S^n(exp)\ra S^{r+2n+1}(exp))\sim 0 $$
for $n\leq N.$
\begin{rem}\label{R9.2}
	Using the same arguments we can also replace the complex $A^{\bullet,\bullet}$ with $s_{\geq N-1}A^{\bullet,\bullet}$. Then, analyzing the graded pieces, we conclude that the Steenbrink limit Hodge complex \ref{F9.1} can be replaced with the data 
	$$(({\rm Tot}(s_{\geq N-1}C),{\rm Tot}(s_{\geq N-1}W(M)),({\rm Tot}(s_{\geq N-1}A),{\rm Tot}(s_{\geq N-1}W(M)_\bC),FA)).$$ 
	%$$0\subset s_{\geq N-1}W(M)_{-N+1}\subset...\subset s_{\geq N-1}W(M)_{N-1}$$
	\tikzset{every picture/.style={line width=0.75pt}} %set default line width to 0.75pt        
	$$\begin{tikzpicture}[x=0.75pt,y=0.75pt,yscale=-1,xscale=1]
		%uncomment if require: \path (0,189); %set diagram left start at 0, and has height of 189
		
		%Straight Lines [id:da43580424001456675] 
		\draw    (374.8,8) -- (374.8,158) ;
		%Straight Lines [id:da5013852791996187] 
		\draw    (228.8,29) -- (374.8,158) ;
		%Straight Lines [id:da5263683019308434] 
		\draw    (225.8,158) -- (442.8,158) ;
		\draw [shift={(444.8,158)}, rotate = 180] [color={rgb, 255:red, 0; green, 0; blue, 0 }  ][line width=0.75]    (10.93,-3.29) .. controls (6.95,-1.4) and (3.31,-0.3) .. (0,0) .. controls (3.31,0.3) and (6.95,1.4) .. (10.93,3.29)   ;
		%Straight Lines [id:da1609641637694499] 
		\draw  [dash pattern={on 0.84pt off 2.51pt}]  (374.8,158) -- (374.8,186) ;
		%Shape: Rectangle [id:dp9169669402562557] 
		\draw  [color={rgb, 255:red, 255; green, 255; blue, 255 }  ,draw opacity=0 ][fill={rgb, 255:red, 221; green, 223; blue, 224 }  ,fill opacity=1 ][dash pattern={on 0.84pt off 2.51pt}] (228.8,7) -- (373.8,7) -- (373.8,29) -- (228.8,29) -- cycle ;
		%Shape: Right Triangle [id:dp9771911776672471] 
		\draw  [color={rgb, 255:red, 0; green, 0; blue, 0 }  ,draw opacity=0 ][fill={rgb, 255:red, 221; green, 223; blue, 224 }  ,fill opacity=1 ] (374.27,156.92) -- (228.8,29) -- (373.35,27.97) -- cycle ;
		%Straight Lines [id:da41306401119952985] 
		\draw  [dash pattern={on 0.84pt off 2.51pt}]  (228.8,29) -- (440.8,29) ;
		%Straight Lines [id:da6464978127007233] 
		\draw    (229.8,184) -- (228.81,9) ;
		\draw [shift={(228.8,7)}, rotate = 89.68] [color={rgb, 255:red, 0; green, 0; blue, 0 }  ][line width=0.75]    (10.93,-3.29) .. controls (6.95,-1.4) and (3.31,-0.3) .. (0,0) .. controls (3.31,0.3) and (6.95,1.4) .. (10.93,3.29)   ;
		
		% Text Node
		\draw (212,-4.6) node [anchor=north west][inner sep=0.75pt]    {$r$};
		% Text Node
		\draw (448,154.4) node [anchor=north west][inner sep=0.75pt]    {$n$};
		% Text Node
		\draw (376.8,161.4) node [anchor=north west][inner sep=0.75pt]  [font=\small]  {$N-1$};
		% Text Node
		\draw (199.8,152.4) node [anchor=north west][inner sep=0.75pt]  [font=\small]  {$-N$};
		% Text Node
		\draw (211,24.4) node [anchor=north west][inner sep=0.75pt]  [font=\small]  {$0$};
		% Text Node
		\draw (403,84.4) node [anchor=north west][inner sep=0.75pt]    {$0$};
		% Text Node
		\draw (273,109.4) node [anchor=north west][inner sep=0.75pt]    {$0$};
		% Text Node
		\draw (344,164.4) node [anchor=north west][inner sep=0.75pt]    {$0$};
		% Text Node
		\draw (280,50.4) node [anchor=north west][inner sep=0.75pt]    {$s_{\geq N-1}W( M)_{r}$};

	\end{tikzpicture}$$
	\begin{center}
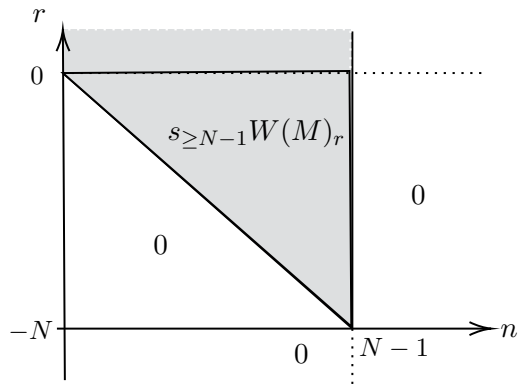
		Picture 9.1. The non-zero elements of the complexes $s_{\geq N-1}W(M)_r$.
	\end{center}
\end{rem}

\subsection*{The complexes $\cw_r$.}
Now, again let $k$ be an arbitrary field, $C/k$ be a smooth curve and $f:X\to C$ be a proper semi-stable degeneration. Let $\cO^*_{X_s}\hookrightarrow \cM_{X_s}$ be the pullback of the log structure on $X_s.$ Again, the map $f:(X_s,\cM_{X_s})\to pt_{log}$ defines the global section $t\in \Gamma(X_s,\cM_{X_s}).$ 

Let us abuse the notation and denote by $\cM$ the extension $\rm{Ex^{log}}\cM$ on the smooth étale site $Sm/X_0$. Further, we will denote by $f_n$ the rightmost differential\footnote{which is induced by the inclusion $\cO^*\otimes (\cM)^{\otimes n}\hookrightarrow (\cM)^{\otimes n+1}$.} in the Koszul complex
\begin{equation}\label{''KC}
	S^{n+1}(\cO^*\to \cM)=...\to \cO^*\otimes \Lambda^n(\cM)\to \Lambda^{n+1}(\cM).
\end{equation}
Let $\cw_r$ be the complex 
$$...\ra \underset{-1}{0} \ra \underset{0}{\Lambda^r \mathcal M^{gr}} \sa{\wedge t} \underset{1}{\Lambda^{r+1} \mathcal M^{gr}}\sa{\wedge t} \underset{2}{\Lambda^{r+2} \mathcal M^{gr}} \sa{\wedge t} ... $$
with differentials given by the multiplication on $t$ and 
$\cw^n_r=0$ for $n<0$ and $n>N$. 

Note that the compositions 
$$f_{j+i}\cdot f_{j+i+1} \cdot ... \cdot f_j: \cO^{*\otimes i}\otimes \Lambda^j\cM\ra \Lambda^{i+j}\cM$$
induce the maps 
$$\alpha_{r,r+i}:\cO^{*\otimes i}\otimes \cw_r \ra \cw_{r+i}.$$

\subsection*{The complexes $\tilde \cw (M)_r$.}
Let $\tilde \cw (M)_r$ be the complex 
$$...\to \underset{-1}{0} \to \underset{0}{\cO^{*\otimes N-1}\!\!\otimes\Lambda^{r+1} \mathcal M^{gr}}\to \underset{1}{\cO^{*\otimes N-2}\!\!\otimes\Lambda^{r+3} \mathcal M^{gr}}\,...\, \to \underset{N-1}{\Lambda^{r+2(N-1)+1} \mathcal M^{gr}} \to \underset{N}{0}\,\,...  $$
where the differential given by the composition of 
the multiplication by $t$
$$\cO^{*\otimes i}\otimes \Lambda^j\cM\sa{{\rm id}\otimes \wedge t}\cO^{*\otimes i}\otimes \Lambda^{j+1}\cM$$
with the differential of the Koszul complex 
$$(\cO^{*\otimes i-1})\otimes \cO^*\otimes \Lambda^{j+1}\cM \sa{{\rm id}\otimes f_{j+1}}\cO^{*\otimes i-1}\otimes\Lambda^{j+2}\cM$$
We define $\tilde \cw (M)_r^n=0$ for $n<0,\, n>N-1$ and 
$\tilde \cw (M)_r^n=\cO^{*\otimes N-1-n}\otimes \Lambda^{r+2n+1}\cM $ otherwise.
\\Note that the diagram 
$$\begin{tikzcd}
	\cO^{*\otimes i+j}\otimes \Lambda^k\cM \arrow[rr, "\wedge t"] \arrow[d, "(id\otimes f_{k+j-1})\cdot...(id \otimes f_k)"] &  & \cO^{*\otimes i+j}\otimes \Lambda^{k+1}\cM \arrow[d, "(id\otimes f_{k+j+1})\cdot...(id \otimes f_{k+1})"] \\
	\cO^{*\otimes i}\otimes \Lambda^{k+j}\cM \arrow[rr, "d"]                                                                 &  & \cO^{*\otimes i-1}\otimes \Lambda^{k+j+2}\cM                                                             
\end{tikzcd}$$
is commutative. So we have the maps 
$$\psi_{r,r+i}: \cO^{*\otimes N+i}\otimes \cw_r \ra \tilde \cw(M)_{r+i}.$$
%In particular, the map 
%$\psi_{N,N}: \cO^{*\otimes N}\otimes \cw_N \ra \tilde \cw(M)_{2N}$
%is an isomorphism of complexes (is not true!).

On the other hand, one can define the morphisms of complexes
$$\varphi_r: \cO^*\otimes \tilde \cw(M)_r\ra \tilde \cw_{r+1}(M),$$
using the commutative diagram 
$$\begin{tikzcd}
	\cO^{*\otimes N}\otimes \Lambda^{r+1}\cM \arrow[d, "f_{r+1}"] \arrow[r] & \cO^{*\otimes N-1}\otimes \Lambda^{r+3}\cM \arrow[d, "f_{r+3}"] \arrow[r] & \cO^{*\otimes N-2}\otimes \Lambda^{r+5}\cM \arrow[r] \arrow[d, "f_{r+5}"] & ... \\
	\cO^{*\otimes N-1}\otimes \Lambda^{r+2}\cM \arrow[r]                & \cO^{*\otimes N-2}\otimes \Lambda^{r+4}\cM \arrow[r]                    & \cO^{*\otimes N-3}\otimes \Lambda^{r+6}\cM \arrow[r]                      & ...
\end{tikzcd}$$
\subsection*{Motivic monodromy filtration}
Let $\cw(M)_{N-1}$ be the cokernel of the map 
$${\rm id}\otimes \psi_{0,N}:\cO^{*\otimes 2N-1}\otimes \cw_0 \ra \tilde \cw(M)_{N-1}.$$
We define the increasing filtration 
$$...\,0\subset \cw(M)_{-N+1} \subset \cw(M)_{-N+2} \subset ...\subset \cw(M)_{N-2} \subset \cw(M)_{N-1}$$
with 
$$\begin{tikzcd}
	\cw(M)_r:={\rm im}(\cO^{*\otimes N-1-r}\otimes \tilde \cw(M)_r \arrow[rrr, " \varphi_{N-2}\cdot...\cdot(({\rm id}\otimes\varphi_r)"] &  &  & \cw(M)_{N-1})
\end{tikzcd}$$
%We will call this filtration \textit{the shifted weight filtration}.
\begin{rem}\label{R9.4}
	Note that ${\rm Gr}^{\cw(M)}_r(\cw(M)_{N-1})=\bigoplus_{n\geq 0,-r} \cw(M)_r^n/\cw(M)_{r-1}^n [-n]$. Indeed, same as in the Steenbrink construction, the n-th differential of $\cw(M)_r$ factors through the map $\cO^*\otimes\cw(M)_{r-1}^{n+1}\sa{\varphi_{r-1}}  \cw(M)_{r}^{n+1}.$ So all differentials of the complex $\rm Gr^{\cw(M)}_r(\cw(M)_{N-1})$ are zero.
	On the other hand, by the same argument as in the proof of Lemma \ref{L9.1}, we get $(\rm Gr^{\cw(M)}_r(\cw(M)_{N-1}))^n=0$ for $n<-r.$
\end{rem}

Now, let us define the stable motivic sheaves $[W(M)_r]\in \rm DA_{et}(\mathit X_0,\bQ)$ by the formula 
$$[W(M)_r]:=\Sigma^\infty_{\Gm}L_{\Al}(\cw(M)_r)(-2(N-1))[-2(N-1)]$$
and define 
$$[\mathrm{Gr}^{W(M)}_r(W(M)_{N-1})]:=\Sigma^\infty_{\Gm}L_{\Al}(\mathrm{Gr}^{\cw(M)}_r(\cw(M)_{N-1}))(-2(N-1))[-2(N-1)].$$
Note that we have an exact sequences
$$ [W(M)_{r-1}]\ra [W(M)_r] \ra [\mathrm{Gr}^{W(M)}_r(W(M)_{N-1})].$$
One can think about sheaves $[W(M)_r]$ as the motivic analog of the filtration ${W(M)_{N-1}\supset W(M)_{N-2}...}$ . Now, we will find out some fundamental properties of the motives $[W(M)_r]$ and $[\mathrm{Gr}^{W(M)}_r(W(M)_{N-1})].$  
\begin{lem}\label{L9.5} The map $\cO^{*\otimes i}\otimes \Lambda^j\cM\to \mathrm{im}( \cO^{*\otimes i}\otimes \Lambda^j\cM \to \Lambda^{i+j}\cM)$ is a $\Al$-homotopy equivalence for any $i$ and $j$.
\end{lem}
\proof It suffice to check that $L_{\Al}(\mathrm{ker}(f_{i+j-1}\cdot...\cdot f_j))=0$. If $i=1$ we can consider the Koszul complex 
$$S^{j+1}(\cO^*\to \cM)=...\to S^{k+1}\cO^*\otimes \Lambda^{j-k}\cM\sa{\partial_{j-k}}S^{k}\cO^*\otimes \Lambda^{j-k+1}\cM \to...$$
Then  we have the exact sequences
$$\mathrm{im}(\partial_{j-k})\to S^{k}\cO^*\otimes \Lambda^{j-k+1}\cM \to \ker(\partial_{j-k-1})$$
and 
$$\mathrm{im}(\partial_{j-1})\to S^2\cO^*\otimes \Lambda^{j-1}\cM \to \ker(f_j)$$
But $L_{\Al}S^k\cO^*=\Lambda^k(\bQ)(k)[k]=0.$
\\If $i>1$ note that $f_{i+j+1}...f_j$ factors through the map 
$$\lambda_{i,j}:\Lambda^i\cO^*\otimes \Lambda^j \cM$$
$$\eta  \otimes \omega \mapsto \eta \wedge \omega.$$
Applying the snake lemma to the diagram 
$$
\begin{tikzcd}
	\mathrm{ker}(f_{i+j-1}...f_j) \arrow[r] \arrow[d] & \cO^{*i}\otimes \Lambda^j\cM \arrow[d, two heads] \arrow[r] & \mathrm{im}(f_{i+j-1}...f_j) \arrow[d, no head, Rightarrow] \\
	{\mathrm{ker}(\lambda_{i,j})} \arrow[r]           & \Lambda^i\cO^{*}\otimes \Lambda^j\cM \arrow[r]              & {\mathrm{im}(\lambda_{i,j})}                               
\end{tikzcd}
$$ 
we conclude that the left vertical arrow is a $\Al$-homotopy equivalence. 
Now, let us consider two Koszul complexes 
$$S^{i+j}(\cO^*\sa{\mathrm{id}+i}\cO^*\oplus \cM)=...\to S^{k+1}\cO^*\otimes \Lambda^{j-k}(\cO^*\oplus \cM)\sa{\partial_{j-k}'}S^{k}\cO^*\otimes \Lambda^{j-k+1}(\cO^*\oplus\cM)\to...;\footnote{Here $i$ is the canonical inclusion $\cO^*\hookrightarrow \cM$.}$$
$$S^{i+j}(\cO^*\sa{\mathrm{id}}\cO^*)=...\to S^{k+1}\cO^*\otimes \Lambda^{j-k}(\cO^*)\sa{\partial_{j-k}''}S^{k}\cO^*\otimes \Lambda^{j-k+1}(\cO^*M)\to...;$$
Then the map $\lambda_{i,j}$ can be given as the composition
$$ \begin{tikzcd}
	\Lambda^i\cO^*\otimes \Lambda^j\cM \arrow[r, hook] & \Lambda^{i+j}(\cO^*\oplus \cM) \arrow[d] \\
	& \Lambda^{i+j}\cM                        
\end{tikzcd}$$
So the kernel of $\lambda_{i,j}$ is the intersection of $\mathrm{im}(\cO^*\otimes \Lambda^{i+j-1}(\cO^*\oplus \cM)\sa{\partial_{i+j-1}'}\Lambda^{i+j}(\cO^*\oplus \cM))$ with $\Lambda^i\cO^*\otimes \Lambda^j\cM.$ Note that the preimage of $\Lambda^i\cO^*\otimes \Lambda^j\cM$ is the direct sum
$$(\cO^*\otimes\Lambda^{i-1}\cO^*\otimes \Lambda^j\cM)\oplus(\cO^*\otimes\Lambda^{i}\cO^*\otimes \Lambda^{j-1}\cM).$$
So $\ker(\lambda_{i,j})$ is the sum of kernels $\ker(\partial_{i}'')=\ker(\cO^*\otimes\Lambda^{i}\cO^*\otimes \Lambda^{j-1}\cM \ra \Lambda^{i+1}\cO^*\otimes \Lambda^{j-1}\cM)$ and $\ker(f_j)=\ker(\cO^*\otimes\Lambda^{i-1}\cO^*\otimes \Lambda^{j}\cM \ra \Lambda^{i-1}\cO^*\otimes \Lambda^{j+1}\cM).$
But both of them are $\Al$-homotopy equivalent to zero by the previous. $\square$

Now, let us describe the graded pieces  $[\mathrm{Gr}^{W(M)}_r(W(M)_{N-1})]$ explicitly. We will denote by $D_i$ the irreducible components of $X_0=D_1\cup ... D_N.$ Let us introduce the notation
$$D_J:
=\bigcap_{j\in J}D_j \,\,\,\,\,\,\,\,\,\,\, D\{m\}:=\bigsqcup_{|J|=m}D_J$$
and denote by $a_m$ the canonical inclusions $D\{m\}\hookrightarrow X_0.$

\begin{prop} For any $r$ we have $$[\mathrm{Gr}^{W(M)}_r(W(M)_{N-1})]\iso\bigoplus_{n\geq 0,-r}a_{r+2n+1*}\bQ_{D\{r+2n+1\}}(-n-r)[-2n-r].$$
\end{prop} 
\proof Firstly, observe that $\mathcal{\bar M}^{gr}=a_{1*}\bQ_{D\{1\}}$ and $\Lambda^j(a_{1*}\bQ_{D\{1\}})=a_{j*}\bQ_{D\{j\}}.$ So by Remark \ref{R9.4} it suffice to construct an $\Al$-homotopy equivalence 
\begin{equation}\label{F9.4}
	\cO^{*\otimes 2(N-1)-n-r}\otimes \Lambda^{r+2n+1} \bar{\mathcal M}^{gr}\ra (\cw(M)_r^n/\cw(M)_{r-1}^n).
\end{equation} 
Let $\beta$ be the canonical map $\tilde \cw(M)^n_r\to \cw(M)_r^n.$
Let us consider the diagram 
$$\begin{tikzcd}
	\mathrm{ker}(\beta\cdot \varphi_{N-2}...\varphi_{r-1}) \arrow[d, "\tilde \varphi"] \arrow[r] & \cO^{*\otimes N-1-r+1}\otimes\tilde \cw(M)_{r-1}^n \arrow[r] \arrow[d, "\mathrm{id}\otimes \varphi_{r-1}"] & \cw(M)_{r-1}^n \arrow[d, "\mu"] \\
	\mathrm{ker}(\beta\cdot \varphi_{N-2}...\varphi_{r}) \arrow[r]                            & \cO^{*\otimes N-1-r}\otimes\tilde \cw(M)_{r}^n \arrow[r]                                              & \cw(M)_r^n                     
\end{tikzcd}
$$
Then \ref{F9.4} can be defined as the map between cokernel of ${\rm id}\otimes \varphi_{r-1}$ and $\mu.$ So it suffice to check that $\tilde \varphi$ is an $\Al$-homotopy equivalence. Let us consider the diagram 
$$
\begin{tikzcd}
	\mathrm{ker}(\beta\varphi_{N-2}...\varphi_{j}) \arrow[r] \arrow[d] & \cO^{*\otimes N-1-j}\otimes \tilde \cw(M)^n_j \arrow[d, "\varphi_{N-2}...\varphi_{j}"] \arrow[r, two heads] & \mathrm{im}(\beta\varphi_{N-2}...\varphi_{j}) \arrow[d, hook] \\
	\mathrm{ker}(\beta) \arrow[r]                                      & \tilde \cw(M)^n_{N-1} \arrow[r, two heads]                                                                  & \cw(M)^n_{N-1}                                               
\end{tikzcd}
$$
Suppose we have $\alpha\in \tilde \cw(M)_{N-1}^n$ and $\beta(\alpha)=0$. Then $\alpha \in \mathrm{im}(\cO^{*\otimes2N-1}\otimes \Lambda^n\cM\to...)$.
Hence $\alpha \in \mathrm{im}(\cO^{*\otimes 2(N-1)-j-n}\otimes \Lambda^{j+2n+1}\cM\to...)$ and there is some $\tilde \alpha$ such that $\varphi_{N-2}...\varphi_{j}(\tilde \alpha)=\alpha$.
Then using commutativity of the diagram we conclude that  $\tilde \alpha\in \mathrm{ker}(\beta\varphi_{N-2}...\varphi_{j})$. So by the snake lemma the left vertical arrow is a surjection. Then for any $j$ we get the diagram with exact lines
$$
\begin{tikzcd}
	\mathrm{ker}(\varphi_{N-2}...\varphi_{j}) \arrow[r, hook] \arrow[d]     & \mathrm{ker}(\beta\varphi_{N-2}...\varphi_{j}) \arrow[r, two heads] \arrow[d] & \mathrm{ker}(\beta) \arrow[d, no head, Rightarrow] \\
	\mathrm{ker}\mathrm{ker}(\varphi_{N-2}...\varphi_{j+1}) \arrow[r, hook] & \mathrm{ker}(\beta\varphi_{N-2}...\varphi_{j+1}) \arrow[r, two heads]         & \mathrm{ker}(\beta)                               
\end{tikzcd}$$
Note that $\ker(\varphi_{N-2}...\varphi_j)=\ker(f_{N-1+2n}...f_{j+2n+1})$. So the Proposition follows from Lemma \ref{L9.5}.$\square$

\subsection*{Relative projection formula.} Suppose $f:X\ra Y$ is a proper morphism of schemes. Let $B\to Rf_*C$ and $B\to A$  be homomorphisms of motivic $\mathbb E_\infty$-algebras. Then by Remark \ref{MonAdj} and monoidality of $f^*$ the maps $f^*B\to C$ and $f^*B\to f^*A$ are also homomorphisms of motivic $\mathbb E_\infty$-algebras. So we can forms two relative tensor products $A\otimes_B Rf_*C$ and $f^*A\otimes_{f^*B}C$.
\begin{prop} We have the canonical weak equivalence
	\begin{equation}\label{F9.5}
		A\otimes_B Rf_*C \to Rf_*(f^*A\otimes_{f^*B}C).
	\end{equation}
\end{prop}
\proof Note that $f^*$ is monoidal so $f^*Rf_*C$ has the natural structure of $f^*B$-algebra. Moreover, the counit of adjunction $\varepsilon_C:f^*Rf_*C\to C$ is a morphism of $f^*B$-algebras. Indeed, by Remark \ref{MonAdj}, $\varepsilon_C$ is a homomorphism of algebras. The homomorphism  $f^*B\to C$ is adjoint to the homomorphism $B\to Rf_*C$. So the diagram 
$$
\begin{tikzpicture}[baseline= (a).base]
	\node[scale=0.8] (a) at (0,0){
		\begin{tikzcd}
			C                          &                                       \\
			& f^*Rf_*C \arrow[lu, "\varepsilon_C"'] \\
			f^*B \arrow[uu] \arrow[ru] &                                      
		\end{tikzcd}
	};
\end{tikzpicture}
$$
is commutative. 
Now, let us define the map \ref{F9.5} as the composition 
$$A\otimes_B Rf_*C\sa{\eta} Rf_*f^*(A\otimes_B Rf_*C)\iso Rf_*(f^*A\otimes_{f^*B}f^*Rf_*C)\sa{Rf_*({\rm id}\otimes \varepsilon)} Rf_*(f^*A\otimes_{f^*B}C).$$
Note that $f$ is proper so $Rf_*$ commutes with geometric realizations. Each maps in the composition can be lifted on the level of the corresponding simplicial sets. So we get the map 
\begin{equation}\label{F9.6}
	(A\otimes_B Rf_*C)^\bullet\ra Rf_*(f^*A\otimes_{f^*B}C)^\bullet 
\end{equation}
For fixed $n$ we have
$$
\begin{tikzcd}
	A\otimes B^{\otimes n}\otimes Rf_*C \arrow[r, "\eta"] & Rf_*f^*(A\otimes B^{\otimes n}\otimes Rf_*C) \arrow[d, "\wr"]                                     &                                             \\
	& Rf_*(f^*A\otimes f^*B^{\otimes n}\otimes f^*Rf_*C) \arrow[r, "Rf_*({\rm id}\otimes \varepsilon)"] & Rf_*(f^*A\otimes f^*B^{\otimes n}\otimes C)
\end{tikzcd}
$$ 
Note that this maps induces by \textit{the exchange transformation} $Ex(f^*_*,\otimes)$ (see \cite{'key-22}, Section 1.1.31). So, using projection formula,  we conclude that (\ref{F9.6}) is a weak equivalence of simplicial objects. $\square$

\subsection*{The limit log scheme.} By the relative projection formula we get 
$$\mathrm{LM}_f\iso Rf_*(\bQ\otimes_{\mathrm{Sym^*}(\bQ\oplus \bQ(-1)[-1])}\mathrm{Sym^*(M^{gr}(-1)[-1])})$$
But $\mathrm{Sym^*}$ is monoidal and commutes with colimits. So 
$$\bQ\otimes_{\mathrm{Sym^*}(\bQ\oplus \bQ(-1)[-1])}\mathrm{Sym^*(M^{gr}(-1)[-1])} \iso \mathrm{Sym^*(\tilde M^{gr}(-1)[-1])}$$
where $\rm \tilde M^{gr}:=\bQ\oplus_{\bQ(1)[1]\oplus\bQ}M^{gr}.$ Note that ${\rm \tilde M^{gr}}=\Sigma^\infty_{\Gm}L_{\Al}\mathrm{Ex^{log}}(\mathcal{\tilde M}^{gr}\otimes \bQ)$ where $\mathcal{\tilde M}^{gr}$ is the colimit of the diagram 
$$\begin{tikzcd}
	\cO^*                                    &                 \\
	\cO^*\oplus \bZ \arrow[r, "t"] \arrow[u] & \mathcal M^{gr}.
\end{tikzcd}$$
Moreover, one can note that the sheaf $\tcM$ corresponds to the fiber product 
$$\begin{tikzcd}
	& {(X_s,\cO^*\hookrightarrow \cM)} \arrow[d, "f"] \\
	X_s \arrow[r] & X_s\times pt_{log}                             
\end{tikzcd}.$$
in the category of virtual log schemes. Let us denote the fiber product by $X_{\rm lim}$. %We will call $X_{\rm lim}$ \textit{the limit log scheme}. 

\begin{rem}
	Note that the fiber product make sense only in the context of virtual log geometry. Indeed, there is no morphisms $X_s\to X_s\times pt_{log}$ in the classical category of logarithmic schemes.
\end{rem}

By the previous we get 
\begin{prop}\label{P16.7}
	The relative projection formula induces the canonical weak equivalence 
	$${\rm LM}_f^*\iso [X_{\rm lim}]^{log}.$$ 
	In particular, any inclusion $k\subset \bC$ induces the isomorphism 
	$$H^*(X^{log}_{\rm lim},\bQ)\iso H^*(X_t,\bQ).$$ 
\end{prop}
Observe that the sheaf $\tcM$ is closely related to the complexes $\cw_r$.

\begin{lem}\label{L9.8}
	For any $n$ we have 
	$$\mathrm{im}(\Lambda^n\cM\sa{\wedge t} \Lambda^{n+1} \cM)\iso \mathrm{ker}(\Lambda^{n+1}\cM \sa{\wedge t} \Lambda^{n+2} \cM )\iso \Lambda^n \tcM.$$
\end{lem} 
As a consequence, 
$$H^0(\cw_r)= \Lambda^r \tcM,\,\,\,H^{N-1}(\cw_r)= \Lambda^{r+N-1}\tcM \,\,\, \mathrm{and} \,\, H^n(\cw_r)=0\,\,\, \mathrm{for} \,\,\, n\neq0,N-1.$$

\begin{rem} \label{R9.9}
	Suppose we have two complexes $A,\, B \in Comp^\bullet(Sh_{et}(Sm/X_0,\bQ))$ together with filtrations $F^\bullet A,\, F^\bullet B$. Let $g:A\to B$ be a filtered morphism. If the induced maps ${\mathrm{Gr}^{FA}_r(A)\to \mathrm{Gr}^{FB}_r(B)}$ are $\Al$-homotopy equivalences then so is $g$. In particular, $L_{\Al}A\iso0$ if $L_{\Al}\mathrm{Gr}^{FA}(A)\iso0$.
\end{rem}
Now, let us construct a quasi-isomorphism ${\rm Sym^*(\tilde M^{gr}(-1)[-1])}\iso [W(M)_{N-1}].$ Firstly, let us consider the diagram 
$$\begin{tikzcd}
	\Lambda^{N-1}\tilde{ \mathcal M}^{gr} \arrow[r] & \cw_N \arrow[r]                           & \Lambda^{2N-1}\tilde{ \mathcal M}^{gr}                   \\
	0 \arrow[u] \arrow[r]                         & \cO^{*\otimes N}\otimes \cw_0 \arrow[u] \arrow[r] & \cO^{*\otimes N}\otimes \Lambda^{N-1}\tilde {\mathcal M}^{gr} \arrow[u]
\end{tikzcd}$$
The lines of the diagram is exact by Lemma \ref{L9.8}. The right vertical arrow is a composition of $\Al$-homotopy equivalences. So the map $\Lambda^{N-1}\tilde {\mathcal M}^{gr}\ra\mathrm{Cone}(\cO^{*\otimes N}\otimes \cw_0\to \cw_N )$ is a $\Al$-homotopy equivalence. 

On the other hand, note that the map 
$$\cO^{*\otimes N-1}\otimes \cw_N \ra \tilde \cw(M)_{N-1}$$
is a $\Al$-homotopy equivalence. Indeed, let we consider the silly filtration on $\cO^{*\otimes N-1}\otimes \cw_N$ and $\tilde \cw(M)_{N-1}$. Then on the n-th graded pieces we get 
$$\mathrm{id}\otimes (f_{N+2n-1}...f_{N+n}): \cO^{*\otimes N-1}\otimes \Lambda^{N+n}\cM \ra \cO^{*\otimes N-1-n} \otimes \Lambda^{N+2n}\cM.$$
Again, this map is a composition of $\Al$-homotopy equivalences and we can use Remark \ref{R9.9}. So we proved
\begin{prop}\label{''P15.10}
	We have a chain of homotopy equivalences 
	$$\begin{tikzcd}
		\cO^{*N-1}\otimes\Lambda^{N-1}\tilde {\mathcal M}^{gr} \arrow[r]                & \mathrm{Cone}(\cO^{*\otimes 2N-1}\otimes \cw_0\to \cw_N) \arrow[ld] \\
		\mathrm{Cone}(\cO^{*\otimes 2N-1}\otimes \cw_0\to \tilde \cw(M)_{N-1}) \arrow[r] & \cw(M)_{N-1}                                                      
	\end{tikzcd}$$
	which induces the quasi-isomorphism 
	$${\rm Sym^*(\tilde M^{gr}(-1)[-1])}\iso [W(M)_{N-1}].$$
\end{prop}

\subsection*{Motivic logarithm of monodromy}
Firstly, let us recall the classical picture. In \cite{PS08}, Steenbrink defined the morphism $\nu:\psi^{Hdg}_f\bQ\to \psi^{Hdg}_f\bQ(-1)$ which maps $W(M)_r$ to $W(M)_{r-2}$. On the level of the complex $Tot(C^{\bullet,\bullet})$ the construction is the following. 

Let $K^{\bullet,\bullet}$ be a bicomplex. Note that the operation $K^{p,q}\mapsto K^{p+1,q-1}$ does not change the totalization. Let us define the map of bicomplexes by the rule
$$\begin{tikzcd}
	0 \arrow[d] \arrow[r]         & {S^{r+1}(exp)[1]} \arrow[r] \arrow[d, "\mathrm{id}"] & {S^{r+3}(exp)[2]} \arrow[d, "\mathrm{id}"] \arrow[r] & ... \\
	{S^{r-2+1}(exp)[0]} \arrow[r] & {S^{r-2+3}(exp)[1]} \arrow[r]                        & {S^{r-2+5}(exp)[2]} \arrow[r]                        & ...
\end{tikzcd}$$
So we get the map 
$$\mathrm{Tot}(\tilde W(M)_r)\to \mathrm{Tot}(\tilde W(M)_{r-2}).$$ 
We also have the map of bicomplexes
$$\begin{tikzcd}
	0 \arrow[d] \arrow[r]     & {S^{0}(exp)[1]} \arrow[r] \arrow[d, hook] & {S^{1}(exp)[2]} \arrow[r] \arrow[d, hook] & ... \\
	{S^{0}(exp)[0]} \arrow[r] & {S^{1}(exp)[1]} \arrow[r]                 & {S^{2}(exp)[2]} \arrow[r]                 & ...
\end{tikzcd}$$
which induces the map 
$$\mathrm{Tot}(W_0)\to \mathrm{Tot}(W_0).$$
Then, taken limit by $r$, we get the endomorphism 
$$Tot(C^{\bullet,\bullet})\sa{\nu}Tot(C^{\bullet,\bullet})$$
such that $\nu(W(M)_r)\subset W(M)_{r-2}.$ Steenbrink proved that $2\pi i\nu$ acts on $\mathbb H^k(X_0,Tot(C^{\bullet,\bullet}))\iso H^k(X_t,\bQ)$ as a logarithm of monodromy $\log (\gamma).$

It should be noted that $\nu$ is nilpotent. %Indeed, $$\nu^{N}(Tot(C^{\bullet,\bullet}))\iso \nu^{m}(Tot(W(M)_{N-1}))\subset \nu^{N-1}Tot(W(M)_{N-3})\subset... \subset W(M)_{-N}$$ and $W(M)_{-N}$ is quasi-isomorphic to zero by the previous.
Indeed, $\nu^{r+p+1}$ maps $S^{r+2p+1}(exp)/S^{p}(exp)\,\,$ to ${S^{r+2p+1}(exp)/S^{r+2p+1}(exp)=0.}$ But $S^{r+2N+1}(exp)/S^N(exp)$ is quasi-isomorphic to 0. So $\nu^{r+N+1}=0$ in the homotopy category.

Now, let us define a motivic analog of $\nu.$ Note that we have the map of complexes 
$$\tilde{\cw}(M)_r\to \cO^*\otimes \tilde{\cw}(M)_{r-2}[1]$$
which can be defined using the diagram
$$\begin{tikzcd}
	0 \arrow[d] \arrow[r]                                    & \cO^{*N-1}\otimes \Lambda^{r+1}\mathcal M^{gr} \arrow[r] \arrow[d, "\mathrm{id}"] & \cO^{*N-2}\otimes \Lambda^{r+3}\mathcal M^{gr} \arrow[r] \arrow[d, "\mathrm{id}"] & ... \arrow[d] \\
	\cO^{*N}\otimes \Lambda^{r-2+1}\mathcal M^{gr} \arrow[r] & \cO^{*N-1}\otimes \Lambda^{r-2+3}\mathcal M^{gr} \arrow[r]                        & \cO^{*N-2}\otimes \Lambda^{r-2+5}\mathcal M^{gr} \arrow[r]                        & ...          
\end{tikzcd}.$$
Let us rewrite this as the map 
$$\cO^*\otimes (\cO^{*N-1-r}\otimes \tilde{\cw}(M)_r)\to (\cO^{*N-1-r+2}\otimes \tilde{\cw}(M)_{r-2})[1].$$
We also have the map 
$$\cO^*\otimes (\cO^{*2(N-1)}\otimes \cw_0)\to (\cO^{*2(N-1)}\otimes \cw_0)[1]$$
defined by the rule
$$
\begin{tikzcd}
	0 \arrow[d] \arrow[r]   & \cO^{*2(N-1)}\otimes \cO^* \arrow[d] \arrow[r] & \cO^{*2(N-1)}\otimes (\cO^*\otimes\mathcal M^{gr}) \arrow[d] \arrow[r] & ...  \\
	\cO^{*2(N-1)} \arrow[r] & \cO^{*2(N-1)}\otimes \mathcal M^{gr} \arrow[r] & \cO^{*2(N-1)}\otimes \Lambda^2 \mathcal M^{gr} \arrow[r]              & ...
\end{tikzcd}.$$
Here the vertical arrows induced by the rightmost differential in the Koszul complex (\ref{''KC}). So we get the morphism 
$$\cO^*\otimes \cw(M)_{N-3}\sa{\nu} \cw(M)_{N-1}[1]$$
which maps $ \cO^*\otimes \cw(M)_{r}$ to $\cw(M)_{r-2}[1].$ It suffice to apply $\Sigma^\infty_{\Gm}$, $L_{\Al}$ and $-\otimes \bQ(-1)[-1]$ and get the maps 
$$[W(M)_r]\sa{\nu} [W(M)_{r-2}](-1).$$
Moreover, using Proposition \ref{''P15.10} we can define the \textit{motivic logarithm of monodromy} 
$$\mathrm{Sym^*(\tilde M^{gr}(-1)[-1]))}\sa{\nu}\mathrm{Sym^*(\tilde M^{gr}(-1)[-1]))}(-1)$$
such that all the diagrams
$$\begin{tikzcd}
	{\mathrm{Sym^*(\tilde M^{gr}(-1)[-1]))}} \arrow[r, "\nu"] & {\mathrm{Sym^*(\tilde M^{gr}(-1)[-1]))}(-1)} \\
	{[\cw(M)_r]} \arrow[u] \arrow[r, "\nu"]                   & {[\cw(M)_{r-2}](-1)} \arrow[u]              
\end{tikzcd}$$
is commutative. Notice that $\nu$ is nilpotent by the same arguments as in the classical case.

\section{Relations with  Ayoub's nearby cycles}
\subsection*{Specialization systems}
Let $B$ be a scheme together with closed and open subsets $Z, \, U\subset B$.  Following Ayoub \cite{'key-18}, a \textit{specialization system} is a collection of functors
$$Sp_f:DA_{et}(X_U,\bQ)\to DA_{et}(X_Z,\bQ),$$ 
one for each $B$-scheme $f:X\to B$, satisfied a certain list of axioms. In particular, for any composition $Y\sa{g}X\sa{f}B$ we have the natural transformations 
$$\alpha_g:g^*Sp_f\to Sp_{fg}g^*, \,\,\,\,\,\, \beta_g: Sp_fRg_*\to Rg_*Sp_{fg}$$
such that $\alpha_g$ is an equivalence when $g$ is smooth and $\beta_g$ is an equivalence when $g$ is proper.

Let $\varphi:B'\to B$ be a morphism of scheme. Let $Z':=\varphi^{-1}(Z)$ and $U':=\varphi^{-1}(U).$ Using $f$, one can restrict a specialization system from $(Z,B,U)$ to $(Z',B',U')$ by an obvious way.

Now, let us consider the following  category: 
\begin{itemize}
	\item the objects given by pairs $(C,D)$ where $C$ is curve over k and $D$ is a divisor;
	\item the morphisms between $(C,D)$ and $(C',D')$ given by morphisms $f:C\to C'$ such that $f^{-1}(D')=D$.
\end{itemize}
\begin{lem}
	Suppose we choose for each $(D,C,C\setminus D)$ a specialization system $Sp^{(C,D)}$ such that the collection $\{Sp^{(C,D)}\}$ is coherent under restrictions. Then the collection can be recovered using $Sp^{(\Al,0)}$.
\end{lem}
\proof Firstly, note that the collection uniquely defined by $Sp^{(C,D)}$ with affine $C.$ Indeed, let $W_j\hookrightarrow C$ be an open affine cover of an arbitrary $C$ and $f:X\to C$ be a $C$-scheme. Let us denote by $\nu_j$ the inclusions $X_{W_j}\hookrightarrow X$. Then 
$$\nu_j^*Sp^{(C,D)}_f\iso Sp^{(C,D)}_{f\nu_j} \nu_j^*=Sp^{(W_j,D_j)}_{f\nu_j} \nu_j^*.$$ 
So we can recover $Sp^{(C,D)}$ by gluing together $Sp^{(W_j,D_j)}_{f\nu_j} \nu_j^*.$  

On the other hand, for affine $C$ we can choose $g:C\to \Al$ such that $D=g^{-1}(0).$ Then $ Sp^{(C,D)}$ is the restriction of $Sp^{(\Al,0)}$.$\square$

\subsection*{Motivic nearby cycles}

Let $(C,s)$ be a curve with a marked point. In \cite{'key-18}, Ayoub constructed \textit{unipotent motivic nearby cycles} $\Upsilon_*$ as a specialization system over $(s,C,\mathring C).$ By the previous it suffice to define $\Upsilon_*$  for $C=\Al,\,\,\,s=0.$

Let $f:X\to \mathbb A^1$ be an $\Al$-scheme. Let us consider the diagram 
$$
\begin{tikzcd}
	{[\Gm \sa{id} \Gm]} \arrow[d, "\Delta"] &                                        \\
	{[\Gm\times \Gm \ra \Gm]}               & {[\Gm \sa{id} \Gm]} \arrow[l, "1"']
\end{tikzcd}.$$
This is a diagram of coalgebras in $Sm/\Gm.$ Note that the canonical functor ${Sm/\Gm \to DA_{et}(\Gm,\bQ)}$ is monoidal. So we get the diagram of $\mathbb E_\infty$-coalgebras in $DA_{et}(\mathring X,\bQ)$
$$\begin{tikzcd}
	\bQ_{\mathring X} \arrow[d, "\Delta"]                                                                   &                   \\
	{\bQ_{\mathring X}\oplus \bQ_{\mathring X}(1)[1]}  & \bQ_{\mathring X} \arrow[l]
\end{tikzcd}$$
and the cosimplicial object $\rm coBar_{\bQ_{\mathring X}\oplus \bQ_{\mathring X}(1)[1]}(\bQ_{\mathring X}, \bQ_{\mathring X})^\bullet.$ Let $\mathcal F$ be a motivic sheaf over $\mathring X$. Then, following Ayoub, 
$$\Upsilon_f\mathcal F:=|i^*Rj_*(\underline{Hom}(\rm coBar_{\bQ_{\mathring X}\oplus \bQ_{\mathring X}(1)[1]}(\bQ_{\mathring X}, \bQ_{\mathring X})^\bullet,\mathcal F))|.$$
\begin{rem}\label{''rem16.2}
	Let $U$ be a Zariski neighborhood of $s.$ Let us consider the composition $X_U\sa{g}X\sa{f}C.$ Then $\Upsilon_f(\mathcal F)= g^*\Upsilon_f(\mathcal F)\iso \Upsilon_{fg}(g^*\mathcal F).$ So we get well defined functor 
	$$\Upsilon_f:DA_{et}(X_\eta,\bQ)\to DA_{et}(X_\eta,\bQ)\footnote{here $\eta$ is a generic point of $C.$}.$$
	
\end{rem}

\subsection*{The equivalence of construction}
Let $f:X\to C$ be a proper semi-stable degeneration. Now, we want to construct a natural weak equivalence between the limit motive $\mathrm{LM}_f$ and $\Upsilon_{\mathrm{id}}([X_\eta]^*).$ Let us start from the following
\begin{lem}
	Let $(X,D)$ be a log smooth scheme with smooth $X$. Let $j:X-D\hookrightarrow X$ be the open immersion. Then there is a natural weak equivalence $$\mathrm{Sym}^{*}(\mathrm{M^{gr}_{(X,D)}}(-1)[-1])\iso Rj_*\bQ.$$
\end{lem}
\proof We have the natural isomorphism of algebras $j^*\mathrm{Sym}^{*}(\mathrm{M^{gr}_{(X,D)}}(-1)[-1])\iso \bQ.$ So we get the homomorphism $\mathrm{Sym}^{*}(\mathrm{M^{gr}_{(X,D)}}(-1)[-1])\ra Rj_*\bQ.$ Let $Y$ be a scheme smooth over $X.$  Then, by Corollary \ref{''C15.6} %, $H^{p,q}_{log}(Y,\bQ)\iso H^{p,q}(Y^*,\bQ)$
and  the construction of  functoriality of $H^{p,q}_{log}$, the map
$$\mathrm{Hom}(Y(q)[p],\mathrm{Sym}^{*}(\mathrm{M^{gr}_{(X,D)}}(-1)[-1]))\ra \mathrm{Hom}(Y(q)[p],Rj_*\bQ)$$
is a quasi-isomorphism for any $p$ and $q \,\,\,\,\square.$

By the previous, the algebra $\mathrm{Sym}^{*}(\mathrm{\tilde M^{gr}}(-1)[-1])$ is equivalent to the bar construction associated with the diagram 
$$\begin{tikzcd}
	{\mathrm{Sym}^*(\mathrm{M^{gr}_{X_0}}(-1)[-1])}              &           \\
	{\bQ_{X_0}\oplus \bQ_{X_0}(-1)[-1]} \arrow[r] \arrow[u, "\mu"] & \bQ_{X_0}
\end{tikzcd}$$
Note that $f$ comes from the global section of $\mathcal M^{gr}_{(X,D)}.$ So the diagram can be extended to the diagram 
$$\begin{tikzcd}
	Rj_*\bQ_{\mathring X}                                    &         \\
	{\bQ_{X}\oplus \bQ_{X}(-1)[-1]} \arrow[r] \arrow[u, "\mu"] & \bQ_{X}
\end{tikzcd}$$
So we get 
$$\mathrm{Sym}^{*}(\mathrm{\tilde M^{gr}}(-1)[-1])\iso i^*(Rj_*\bQ_{\mathring X}\otimes_{[\Gm]_X^*}\bQ_X)\iso|i^*(Rj_*\bQ_{\mathring X}\otimes_{[\Gm]_X^*}\bQ_X^\bullet)|.$$
Moreover, 
$$i^!((Rj_*\bQ_{\mathring X})\otimes_{[\Gm]_X^*}\bQ_X^n)=i^!(Rj_*\bQ_{\mathring X}\otimes(\bQ\oplus \bQ(-1)[-1])^{\otimes n}\otimes \bQ)\iso \oplus i^!Rj_*\bQ(-k)[-k]=0.$$ 

So $$(Rj_*\bQ_{\mathring X})\otimes_{[\Gm]_X^*}\bQ_X^\bullet \iso Rj_*j^*((Rj_*\bQ_{\mathring X})\otimes_{[\Gm]_X^*}\bQ_X^\bullet )\iso Rj_*(\bQ\otimes_{[\Gm]_{\mathring X}^*}\bQ)^\bullet$$
and 
\begin{equation}\label{''ner}
	\mathrm{Sym}^{*}(\mathrm{\tilde M^{gr}}(-1)[-1])\iso |i^*Rj_*(\bQ\otimes_{[\Gm]_{\mathring X}^*}\bQ)^\bullet|
\end{equation}
Here the relative tensor product is the product associated with the diagram 
\begin{equation}\label{''F16.2}
	\begin{tikzcd}
		\bQ_{\mathring X}                                                            &                   \\
		{\bQ_{\mathring X}\oplus \bQ_{\mathring X}(-1)[-1]} \arrow[r] \arrow[u, "\mu"] & \bQ_{\mathring X}
	\end{tikzcd}.
\end{equation}

\begin{thm}\label{vanish_cycl_thm}
	Suppose that $k$ satisfies the Beilinson-Soulé vanishing conjecture. Then  
	$${\mathrm{Sym}^{*}(\mathrm{\tilde M^{gr}}(-1)[-1])\iso \Upsilon_f(\bQ)}$$
	and, consequently, 
	$$\mathrm{LM}_f\iso \Upsilon_{\mathrm{id}}([X_\eta]^*).$$
\end{thm} 
\proof It suffice construct the first equivalence. By Remark \ref{''rem16.2} we may assume that $C$ is affine. Let us choose a function $g$ on $C$ with the property $g^{-1}(0)=s$ and denote by $\varphi$ the composition $X\sa{f} C\sa{g}\Al.$ Note that $\Upsilon_f\bQ$ can be given by the formula $|i^*Rj_*(\rm Bar_{\bQ_{\mathring X}\oplus \bQ_{\mathring X}(-1)[-1]}(\bQ_{\mathring X}, \bQ_{\mathring X})^\bullet)|$ where the bar construction associated with the diagram 
\begin{equation}\label{''F16.3}
	\begin{tikzcd}
		\bQ_{\mathring X}                                                                   &                   \\
		{\bQ_{\mathring X}\oplus \bQ_{\mathring X}(-1)[-1]} \arrow[u, "\Delta^*"] \arrow[r] & \bQ_{\mathring X}
	\end{tikzcd}
\end{equation} 
Then, according \ref{''ner}, it enough to check that the diagrams \ref{''F16.2} and \ref{''F16.3} define equivalent bar constructions.

Note that the map $\Delta^*$ is dual to the global section $\varphi=g\cdot f$ of sheaf $\cO^*$ on $\mathring X.$ On the other hand, the map $\mu$ corresponding to the map $(id,\varphi):\cO^*\oplus \bZ \to \cO^*.$ The both map can be represented as the pullbacks from $\Gm$ under $\varphi.$ So we should compare the diagrams
$$\begin{tikzcd}
	\bQ_{\Gm}                                                         &         &  & \bQ_{\Gm}                                                                           &         \\
	{\bQ_{\Gm}\oplus \bQ_{\Gm}(-1)[-1]} \arrow[r] \arrow[u, "{(1,t)^*}"] & \bQ_{\Gm} &  & {\bQ_{\Gm}\oplus \bQ_{\Gm}(-1)[-1]} \arrow[r] \arrow[u, "{(1,t)\otimes \bQ(-1)[-1]}"] & \bQ_{\Gm}
\end{tikzcd}.$$

\begin{rem}\label{''rem16.5}
	Let $\mathcal C^\otimes$ be a monoidal $\infty$-category. Let $A\sa{\alpha}A'$ be a weak equivalence of $\mathbb E_\infty$-algebras.
	Suppose we have the diagrams 
	$$\begin{tikzcd}
		B                               &   \\
		A \arrow[u, "f"] \arrow[r, "g"] & C
	\end{tikzcd}
	\mathrm{;}\,\,\,\, \,\,\,\,\,\,\,\,\,\,\,\,\,\,\,\,\,\,\,\,\,\,\,\,
	\begin{tikzcd}
		B'                               &   \\
		A' \arrow[u, "f'"] \arrow[r, "g'"] & C
	\end{tikzcd}$$ 
	of $\mathbb E_\infty$-algebras such that $g'\alpha=g$, $B'=B\otimes_AA'$ and $f'=f\otimes_AA'.$ Then $B'\otimes_A'C\iso B\otimes_AC$ as objects of $\mathcal C^\otimes.$ Indeed, by results of (\cite{key-1}, Section 4.4.3), 
	$$B'\otimes_A'C=(B\otimes_AA')\otimes_A'C\iso B\otimes_A(A'\otimes_A'C)\iso B\otimes_AC.$$
\end{rem}

\begin{rem} \label{''rem16.6}
	Let $B$ and $B'$ be as above. Then $B\iso B'$ in $\mathcal C^\otimes.$ In fact, $\alpha:A\to A'$ induce the weak equivalence $B\iso A\otimes_AB\to A'\otimes_AB.$  
\end{rem}
\begin{rem}\label{''rem16.7}
	Note that the unit of adjunction $\mathrm{Sym}^*([\Gm]_{\Gm}^*)\sa{\epsilon} [\Gm]_{\Gm}^*$ is a weak equivalence of algebras. Indeed, by conservativity of the forgetful functor it suffice to check that $\epsilon$ is a weak equivalence of the motivic sheaves. But the natural inclusion $[\Gm]_{\Gm}^*\to \mathrm{Sym}^*([\Gm]_{\Gm}^*)$ is an equivalence and we can use the triangle identities.
\end{rem}
Let us consider the automorphism 
$$(1,-1):\bQ\oplus \bQ(-1)[-1]\ra \bQ\oplus \bQ(-1)[-1].$$
By the Remark \ref{''rem16.7} it is an automorphism of algebras. For any module $M$ let us denote by $M\otimes_{(1,-1)}[\Gm]^*$ the tensor product $M\otimes_{\bQ\oplus \bQ(-1)[-1]}\bQ\oplus \bQ(-1)[-1]$ with respect to $(1,-1)$. Then $\bQ\otimes_{(1,-1)}[\Gm]^*\iso \bQ$ and $(1,t)\otimes_{(1,-1)}[\Gm]^*=(1,-t)$ in the homotopy category of $\mathbb E_\infty$-algebras. Indeed, the first equivalence follows from Remark \ref{''rem16.6}. To prove the second equality, let us consider the diagram 
$$\begin{tikzcd}
	\bQ \arrow[r, "a"]                                                & {\bQ\iso \bQ\otimes_{(1,-1)}[\Gm]^*} \\
	{\bQ\oplus\bQ(-1)[-1]} \arrow[r, "{(1,-1)}"] \arrow[u, "{(1,t)}"] & {\bQ\oplus\bQ(-1)[-1]} \arrow[u, "{(b,c)}"].                      
\end{tikzcd}$$ 
Note that $\bQ\otimes_{(1,-1)}[\Gm]^*$ is an algebra with unit. So $a=b=1.$ Then, by commutativity, $c=-t.$

Now, one can apply Remark \ref{''rem16.5} to prove the Theorem \ref{vanish_cycl_thm}. Namely, we should check that the homomorphisms $\bQ_{\Gm}\oplus \bQ_{\Gm}(-1)[-1]\sa{(1,t)^*}\bQ_{\Gm}$ and ${\bQ_{\Gm}\oplus \bQ_{\Gm}(-1)[-1]\sa{(1,-t)}\bQ_{\Gm}}$ are coincide in the homotopy category of $\mathbb E_\infty$-algebras.
Since the algebra $[\Gm]_{\Gm}^*$ is free it suffice to check that $(1,t)^*=(1,-t)$ in the homotopy category of motivic sheaves $DA_{et}(\Gm,\bQ)_{\setminus \bQ}.$ Then the Theorem follows from the Lemma.

\begin{lem} 
	Suppose that $k$ satisfies the Beilinson-Soulé vanishing conjecture. Let ${f: \bQ\ra \bQ(1)[1]}$ be a morphism in $DA_{et}(\Gm,\bQ).$ Let $\check f$ and $f(-1)[-1]$ be the images of $f$ in $\mathrm{Hom}(\bQ(-1)[-1],\bQ)$ under the functors $\mathrm{Hom}(-,\bQ)$ and $-\otimes \bQ(-1)[-1].$ Then $\check f=-f(-1)[-1].$
\end{lem}
\proof Let $DTM(\Gm)_\bQ\subset DA_{et}(\Gm,\bQ)$ be the triangulated subcategory of Tate motives. By result of \cite{Lev93}, there is weight t-structure on $DTM(\Gm)_\bQ$ together with tensor exact functors
$$gr^W_q:DTM(\Gm)_\bQ \ra T_q$$
where $T_q\iso D(Vect_\bQ)$ is a full triangulated subcategory of $DTM(\Gm)_\bQ\subset DA_{et}(\Gm,\bQ)$ generated by $\bQ(q)[n]$, $n\in \bZ.$

By assumption, $\Gm$ is also satisfies Beilinson-Soulé conjecture. Indeed, 
$$H^p(\Gm,\bQ(q))=H^{p-1}(k,\bQ(q-1))\oplus H^p(k,\bQ(q)).$$ 
%and for $q=1$, $p\leq 0$ we have $H^p(\Gm,\bQ(1))=H^{p-1}(k,\bQ)=0.$
So we can apply Theorem 1.4 of \cite{Lev93} and get the following: 
\begin{itemize}
	\item There is the non-degenerate t-structure on $DTM(\Gm)_\bQ$ with the heart $TM(\Gm)_\bQ$ which consists $M\in DTM(\Gm)_\bQ$ such that all $gr^W_i(A)\iso \bQ(-i)^{n_i}$ for some $n_i\in \mathbb N;$
	\item  The category $TM(\Gm)_\bQ$ is closed under extensions and $$\mathrm{Ext}^1_{TM(\Gm)_\bQ)}(B,A)\iso \mathrm{Hom}_{DA_{et}(\Gm,\bQ)}(B[-1],A);$$
	\item
	The fiber functor $\rm \oplus_i gr^{W}_i:TM(\Gm)_\bQ\ra Vect_\bQ$ is a faithful exact tensor functor. By the same arguments as \cite{DG05}, the group $\rm Aut(\oplus_i gr^{W}_i)$ is isomorphic to $\Gm\ltimes U$ where $U$ is  a pro-unipotent algebraic group.
\end{itemize}
Now, let $\bQ(1)\to A\to \bQ$ be the extension associated with $f$ and let 
$$\bQ\to \check A\to \bQ(-1); \,\,\,\,\,\,\,\,\,\, \bQ\to A(-1)\to \bQ(-1)$$
be the extensions associated with $\check f$ and $f(-1)[-1].$ It suffice to construct an appropriate map of the extensions in the category of representations of $G=\Gm\ltimes U$. Let $(\lambda,u)$ be an element of $\Gm\ltimes U$. Let $(v_1,v_2)$ be the basis of $A$ such that $v_1$ corresponds to the inclusion $\bQ(1)\to A.$ By construction, the action of $G$ preserves the weight filtration so 
$$\rho(\lambda,u)=
\begin{pmatrix}
	\lambda & a_{(\lambda,u)} \\
	0 & 1 
\end{pmatrix}.  $$
Then the action of $(\lambda,u)$ on $A(-1)$ is given by
$ \begin{pmatrix}
	1 & a_{(\lambda,u)}\lambda^{-1} \\
	0 & \lambda^{-1} 
\end{pmatrix}.  $
Finally, the action on $\check A$ given by 
$$(\rho(\lambda,u)^{-1})^{\mathrm T}=\begin{pmatrix}
	\lambda^{-1} & 0 \\
	-a_{(\lambda,u)}\lambda^{-1} & 1 
\end{pmatrix}$$
where $\delta_2=\begin{pmatrix}
	0 \\
	1 
\end{pmatrix} $ corresponds to the inclusion $\bQ\to \check A.$ So we have the map of extensions
$$\begin{tikzcd}
	\bQ \arrow[d, "-1"] \arrow[r] & A(-1) \arrow[d,"B"] \arrow[r] & \bQ(-1) \arrow[d, "1"] \\
	\bQ \arrow[r]                & \check A \arrow[r]                                                           & \bQ(-1)                
\end{tikzcd}$$
where $B=\begin{pmatrix} 0 & 1 \\ -1 & 0  \end{pmatrix}$ with respect to the basis $v_1,v_2$ and $\delta_1,\delta_2 \,\,\, \square.$

\appendix

\section{The inverse image of an algebraic group }\label{A}
Suppose that $\mathcal{G}$ is an abelian algebraic group over $\mathrm{k}$.
Then the for each $X\in\mathrm{Sch/k}$ the group scheme $\mathcal{G}\times_{\mathrm{k}}X$
defines the sheaf $\mathcal{G}_{X}\in\mathrm{Sh_{et}(}\mathrm{Sm}/X,\mathbb{Z})$.
Suppose that there is a morphism $f\!:\!Y\longrightarrow X$. For
any scheme $U$ smooth over $X$ and $\varphi\in\mathcal{G}(U)$ the
composition
\[
(U\overset{\varphi}{\longrightarrow}\mathcal{G})\longmapsto(U\times_{X}Y\longrightarrow U\overset{\varphi}{\longrightarrow}\mathcal{G})
\]
define the map 
\[
f_{\mathcal{G}}:\mathcal{G}_{X}\longrightarrow f_{*}\mathcal{G}_{Y}
\]
where $f_{*}:\mathrm{Sh_{et}(}\mathrm{Sm}/Y,\mathbb{Z})\longrightarrow\mathrm{Sh_{et}(}\mathrm{Sm}/X,\mathbb{Z})$.
Let $f_{\mathcal{G}}^{\natural}$ be the adjoint map.
\begin{thm}\label{AA.1}
	Suppose that $\mathcal{G}$ is smooth over $k$. Then for any morphism
	$f\!:\!Y\longrightarrow X$  the map $f_{\mathcal{G}}^{\natural}:f^{-1}\mathcal{G}_{X}\longrightarrow\mathcal{G}_{Y}$
	is an isomorphism.
\end{thm}
\textit{Proof:} Let $G$ be an abelian group. Let $I\subset\mathbb{Z}[G]$
be the subgroup generated by $[0]$ and $J\subset\mathbb{Z}[G\times G]$
be the subgroup generated by all $[(0,g)]$ and $[(g,0)]$ for any
$g\in G$. Let us define the complex $Q_{\geq-1}^{*}(G)$:
\[
\underset{-1}{\mathbb{Z}[G\times G]/J}\stackrel{d}{\longrightarrow}\underset{0}{\mathbb{Z}[G]/I}
\]
where $d([(g_{1},g_{2})])=[g_{1}]+[g_{2}]-[g_{1}+g_{2}]$.

Now, let $S$ be a scheme. Then we can define the complex of presheaves
$\boldsymbol{Q}_{\geq-1}(\mathcal{G}_{X})$ on $\mathrm{Sm/}S$ by
the rule 
\[
U\longmapsto Q_{\geq-1}^{*}(\mathcal{G}(U)).
\]
Note that $H^{0}(Q_{\geq-1}^{*}(G))=G$ for any $G$. So $H^{0}(\boldsymbol{Q}_{\geq-1}(\mathcal{G}_{X}))=\mathcal{G}_{X}$.
For any $f:Y\longrightarrow X$ we have the canonical map 
\begin{equation}\label{A.1}
	\boldsymbol{Q}_{\geq-1}(\mathcal{G}_{X})\longrightarrow f_{*}\boldsymbol{Q}_{\geq-1}(\mathcal{G}_{Y})
\end{equation}
induced by $\mathbb{Z}[f_{\mathcal{G}}]$ and $\mathbb{Z}[f_{\mathcal{G}\times\mathcal{G}}]$.
The map $f_{\mathcal{G}}$  is obtained from \ref{A.1}
by applying of the functor $H^{0}$. So it suffice to check that the
map adjoint to \ref{A.1} is an isomorphism. 

Observe that $f^{-1}\mathbb{Z}[\mathcal{G}_{X}]=\mathbb{Z}[\mathcal{G}_{Y}]$
by definition of $f^{-1}$. Moreover, $\mathbb{Z}[f_{\mathcal{G}}]$
and $\mathbb{Z}[f_{\mathcal{G}\times\mathcal{G}}]$ are the units
of adjunction. The map
\begin{equation}\label{A.2}
	\mathbb{Z}[\mathcal{G}_{X}]/I_{X}\longrightarrow f_{*}(\mathbb{Z}[\mathcal{G}_{Y}]/I_{Y})
\end{equation}
is obtained by applying of the functor $H^{0}$ to the morphism of
complexes
\[
\xymatrix{\mathbb{Z}[pt_{X}]\ar[r]^{\mathbb{Z}[f_{pt}]}\ar[d]^{[0]} & f_{*}\mathbb{Z}[pt_{Y}]\ar[d]^{[0]\,\,\,\,\,\,\,\,\,\,\,\,\,\,.}\\
	\mathbb{Z}[\mathcal{G}_{X}]\ar[r]^{\mathbb{Z}[f_{\mathcal{G}}]} & f_{*}\mathbb{Z}[\mathcal{G}_{Y}]
}
\]
The horizontal map are the units of adjunction. So the adjoint morphism
is an identity and \ref{A.2} is an isomorphism. Analogically,
the map
\begin{equation}\label{A.3}
	\mathbb{Z}[\mathcal{G}_{X}\times\mathcal{G}_{X}]/J_{X}\longrightarrow f_{*}(\mathbb{Z}[\mathcal{G}_{Y}\times\mathcal{G}_{Y}]/J_{Y})
\end{equation}
is obtained from the morphism of complexes
\[
\xymatrix{\mathbb{Z}[\mathcal{G}_{X}]\oplus\mathbb{Z}[\mathcal{G}_{X}]\ar[rr]^{\mathbb{Z}[f_{\mathcal{G}}]\oplus\mathbb{Z}[f_{\mathcal{G}}]}\ar[d]^{i_{(0,1)}+i_{(1,0)}} &  & f_{*}(\mathbb{Z}[\mathcal{G}_{Y}]\oplus\mathbb{Z}[\mathcal{G}_{Y}])\ar[d]^{\,\,\,\,\,\,\,\,\,\,\,\,\,\,\,\,\,\,\,\,\,\,\,\,\,\,\,\,\,.}\\
	\mathbb{Z}[\mathcal{G}_{X}\times\mathcal{G}_{X}]\ar[rr]^{\mathbb{Z}[f_{\mathcal{G}\times\mathcal{G}}]} &  & f_{*}\mathbb{Z}[\mathcal{G}_{Y}\times\mathcal{G}_{Y}]
}
\]
So the map adjoint to \ref{A.3} is an isomorphism by
the same arguments. $\square$
\begin{cor}\label{CA.2}
	For a scheme $S$ let us denote by $\mathcal{O}_{\mathrm{Sm/}S}^{*}$
	the étale sheaf of invertible functions on $\mathrm{Sm/}S$ . Let $f\!:\!Y\longrightarrow X$
	be a morphism and $f^{\natural}:f^{-1}\mathcal{O}_{\mathrm{Sm/}X}^{*}\longrightarrow\mathcal{O}_{\mathrm{Sm/}Y}^{*}$
	be the canonical map. Then the map $f^{\natural}$ is an isomorphism.
\end{cor}
%\begin{rem}
%In paper \cite{'key-6} Eilenberg and MacLane define for any abelian group the complex $Q^*(G).$ (cubical construction). By construction the first truncation of $Q^*(G)=??$ coincides with the complex and the cohomology $H^*(Q^*(G))$ coincides with the stable cohomology of Eilenberg-MacLane spaces. One can compute that $H^m(K(G,n),\bQ)=0$ for any abelian group $G$. So the complex of sheaves $Q(\mathcal{G}_{X}) \bQ$ is an acyclic resolution of $\mathcal{G}_{X}\otimes \bQ\in\mathrm{Sh_{et}(}\mathrm{Sm}/X,\mathbb{Q})$. Then using the same methods as above one can check that the natural map $f^-1\to Q$ is an isomorphism. This is also implise Theorem \ref{AA.1}.
%\end{rem}

\section{($\infty$,1)-Grothendieck construction}\label{B}

In classical category theory the Grothendieck construction provide
the equivalence of categories between the category of Grothendieck
fibration over the category $C$ and the category of pseudofunctors
$\mathrm{Fun}(C^{\mathrm{op}},\mathbf{Cat})$ . The similar pictures
exists in the word of infinity categories. The notation of Cartesian
fibrations of simplicial sets is a direct analog of Grothendieck fibrations
of 1-categories. Moreover, for any simplicial set $S$ we
have the canonical equivalence of $\infty$-categories $\mathrm{Fun}_{\infty}(S^{\mathrm{op}},\mathbf{Cat}_{(\infty,1)})\simeq\mathbf{Cart}_{S}$
where $\mathbf{Cart}_{S}$ is the category of Cartesian fibrations
over $S$. 
%\begin{rem}
%Strictly speaking $\mathbf{Cart}_{S}\!\overset{\mathrm{def}}{=}\!N((Set_{\triangle}^{+})_{/S})$
%where $(Set_{\triangle}^{+})_{/S}$ is the category of marked simplicial
%set over $S$ together with (simplicial) Cartesian model structure
%( see HTT 3.2).
%\end{rem}
This equivalence can be describe using the straightening and unstraightening
functors:
\[
\mathrm{Un^{+}}\!:\!\mathrm{Fun}_{\infty}(S^{\mathrm{op}},\mathbf{Cat}_{(\infty,1)})\rightleftarrows\mathbf{Cart}_{S}\!:\!\mathrm{St^{+}}
\]
The interested reader can find all details in Chapter 3.2 of \cite{key-3}. We will
abuse the notation and denote by $\int_{C}F$ the value of functor
$\mathrm{Un^{+}}$ on the presheaf $F$. In fact, we will only consider the case when $S$ is the nerve of a 1-category $C$.
%In fact, we will only use the ustraightening in the case then the simplicial set $S$ is the nerve of a 1-category $C$. 
%\begin{rem}
%Let $A$ be a simplicial model category (for example $SSet_{Joal}$)
%and $C$ be a 1-category. Then there is an equivalence between $\mathrm{Fun}_{\infty}(N(C),N(A))$
%and $N(Fun(C,A))$ where $\mathrm{Fun}(C,A)$ is a simplicial model
%category of functors with projective model structure.
%\end{rem}
\subsection*{Relative nerve construction.}

Let $F:C\longrightarrow SSet$ be a 1-functor. We define the simplicial  
set $N_{F}(C)$ as follows. Let $J$ be a linear order  set with $J=|n|$. One can consider $J$ as a category. We will denote by $\Delta^J$ the associated simplicial set. The n-cell of $N_{F}(C)$ contained the following  data:
\begin{itemize}
	\item [{(RN1)}] A functor $\sigma$ from $J$ to $C$. 
	\item [{(RN2)}] For every nonempty subset $J'\subseteq J$ having a maximal
	element $j_{max}'$ , a map ${\tau(J'):\,\triangle^{J'}\rightarrow F(\sigma(j_{max}'))}$ 
\end{itemize}
together with the condition
\begin{itemize}
	\item [{(RN3)}] For nonempty subsets $J''\subseteq J'\subseteq J$ , with
	maximal elements $j_{max}''\in J''$,$j_{max}'\in J',$ the diagram
	\[
	\xymatrix{\triangle^{J''}\ar[r]\ar[d] & F(\sigma(j_{max}''))\ar[d]^{F(\sigma(j_{max}'')\longrightarrow(\sigma(j_{max}))}\\
		\triangle^{J'}\ar[r] & F(\sigma(j_{max}'))
	}
	\]
	is required to commute. 
\end{itemize}
The simplicial set $N_{F}(C)$ is called \textit{the relative nerve of} $F$. Suppose $F$ takes value in weak Kan complexes. Then $N_{F}(C)$ can be identified with $\int_{C}F$ (\cite{key-3}, Section 3.2.5).
\begin{rem}
	Note that a 0-simplex of $N_F(C)$ is a pair $(c,Y)$ where $c\in Ob(C)$ and $Y$ is a 0-simplex of $F(c)$. Moreover, a 1-simplex of $N_F(C)$ is also a pair which contains a morphism $\alpha:c_1\to c_2$ of $C$ and a 1-simplex $F(\alpha)(Y_1)\to Y_2$. So the relative nerve compatible with the classical Grothendieck construction. 
\end{rem}
\begin{example}
	Let $S$ be a simplicial set and $\underline{S}:C\longrightarrow SSet$
	be the constant functor. Then $N_{\underline{S}}(C)$ is the Cartesian
	product $S\times N(C)$. Indeed, the maps from the condition (RN2)
	is just a point of $S_{|J'|}$. By (RN3) all such maps unique defined
	by the one map $\tau(J)$. %
	\begin{comment}
	\textcolor{blue}{You should wright about marked simplicial sets and
	describe the functor $N_{F}^{+}$(see nLab and HTT 3.2)}
	\end{comment}
\end{example}

\subsection*{Presentable fibrations.}

Following \cite{key-3} we will call a map of simplicial
sets ${p:X\to S}$ a \textit{presentable fibration}
if it is a Cartesian and cocartesian fibration such that
each fiber ${X_{s}=X\times_{S}\{s\}}$ is a presentable \ensuremath{\infty}-category. 
\begin{prop}
	(Proposition 5.5.3.3. of \cite{key-3})
\end{prop}
\begin{itemize}
	\item Let $p\!:\!X\rightarrow S$ be a Cartesian fibration of simplicial
	sets, classified by a map ${\chi:S^{\mathrm{op}}\longrightarrow\mathbf{Cat}_{(\infty,1)}}$
	Then $p$ is a presentable fibration if and only if $\chi$ factors
	through $\mathbf{Pr}^{\mathrm{R}}\subseteq\mathbf{Cat}_{(\infty,1)}$. 
	\item Let $p\!:\!X\rightarrow S$ be a cocartesian fibration of simplicial
	sets, classified by a map ${\chi:S\longrightarrow\mathbf{Cat}_{(\infty,1)}}$
	Then $p$ is a presentable fibration if and only if $\chi$ factors
	through $\mathbf{Pr}^{\mathrm{L}}\subseteq\mathbf{Cat}_{(\infty,1)}$
\end{itemize}
\begin{cor}
	(Corollary 5.5.3.4. of \cite{key-3},) For every simplicial set S, there
	is a canonical bijection $[S,\mathbf{Pr}^{\mathrm{L}}]\simeq[S^{\mathrm{op}},\mathbf{Pr}^{\mathrm{R}}]$
	where $[S,C]$ denotes the collection of equivalence classes of objects
	of $\mathbf{Fun}_{\infty}(S,C)$. 
\end{cor}
%Proof: According to Proposition 5.5.3.3, both $[S,\mathbf{Pr}^{\mathrm{L}}]$and
%$[S^{\mathrm{op}},\mathbf{Pr}^{\mathrm{R}}]$ can be identified with
%the collection of equivalence classes of presentable fibrations $X\rightarrow S$.
%$\Square$ 

For $F:C\longrightarrow\mathbf{Pr}^{\mathrm{R}}$ let us denote by
$F^{\rm adj}$ the corresponding functor $C^{\mathrm{op}}\longrightarrow\mathbf{Pr}^{\mathrm{L}}.$ Note that $F^{\rm adj}$ can be described as follows:
\begin{itemize}
	\item $F^{\rm adj}(c)=F(c)$ for any $c\in C$;
	\item $F^{\rm adj}(\alpha)$ is right-adjoint to $F(\alpha)$ for any $\alpha: c_1\to c_2$. Further, we will denote $F^{\rm adj}(\alpha)$ and $F(\alpha)$ by $\alpha^*$ and $\alpha_*$. 
\end{itemize}

According to the Corollary, $\int_{C}(F^{\mathrm{op}})$
and $(\int_{C^{\mathrm{op}}}F^{\rm adj})^{\mathrm{op}}$ are
equivalent as $\infty$-categories.
\begin{rem}\label{RB.5}
	Let $S=\Delta^1$. Then the data $F:\Delta^1\to \mathbf{Pr}^{\rm R}$ is equivalent to the pair of categories $F(0)$ and $F(1)$ together with a pair of adjoint functors $\alpha^*$ and $\alpha_*:$
	$$\alpha^*:F(1) \to F(0):\alpha_*.$$ 
	Then the corresponding relative nerves can be described as follows: 
	\begin{itemize}
		\item the objects of $\int_{\Delta^1}(F^{\mathrm{op}})$ are the disjoint union of objects $F(0)$ and $F(1)$. For any $x,y\in Ob(\int_{C}(F^{\mathrm{op}}))$ we have
		$$\mathrm{Hom}(x,y)^{\bullet}=\begin{cases}
			\varnothing & \mathrm{if}\,\,x\in F(1),y\in F(0)\\
			\mathrm{Hom}_{F(i)}(y,x)^\bullet & \mathrm{if}\,\,x,y\in F(i)\\
			\mathrm{Hom}_{F(1)}(y,\alpha_*(x))^{\bullet} & \mathrm{if}\,\,x\in F(0),y\in F(1)
		\end{cases}.$$
		\item the objects of $(\int_{{\Delta^1}^{\mathrm{op}}}F^{\rm adj})^{\mathrm{op}}$ are the disjoint union of objects $F(0), \, F(1)$ and 
		$$\mathrm{Hom}(x,y)^{\bullet}=\begin{cases}
			\varnothing & \mathrm{if}\,\,x\in F(1),y\in F(0)\\
			\mathrm{Hom}_{F(i)}(y,x)^\bullet & \mathrm{if}\,\,x,y\in F(i)\\
			\mathrm{Hom}_{F(0)}(\alpha^*(y),x)^{\bullet} & \mathrm{if}\,\,x\in F(0),y\in F(1)
		\end{cases}.$$
	\end{itemize}
	Note that for any $x\in F(0)$ and $y\in F(1)$ we have the canonical map 
	\begin{equation}\label{FB.1}
		\mathrm{Hom}^\bullet_{F(1)}(y,\alpha_*(x))\to\mathrm{Hom}^\bullet_{F(0)}(\alpha^*(y),x)
	\end{equation}
	which is a weak equivalence of simplicial sets. Moreover, the maps \ref{FB.1} define the functor
	\[
	\!\int_{\Delta^1}(F^{\mathrm{op}})\ra(\int_{\Delta^{1\,\mathrm{op}}}F^{\rm adj})^{\mathrm{op}}
	\] 
	which is an equivalence of $\infty$-categories.
\end{rem}

\subsection*{Construction of the equivalence.}
Let $C$ be a 1-category and $F:N(C)\to \mathbf{Pr}^{\rm R}$ be a functor. 
Let us construct the canonical equivalence of categories $\int_{C}(F^{\mathrm{op}})$ and $(\int_{C^{\mathrm{op}}}F^{\rm adj})^{\mathrm{op}}.$
%\[
%\varphi\!:\!\int_{C}(F^{\mathrm{op}})\overset{\sim}{\longrightarrow}(\int_{C^{\mathrm{op}}}F^{\vee})^{\mathrm{op}}
%\]
%at least in the case when $C$ is a 1-category. 
First of all, observe that by (RN3) the functors $\tau(J')$ factor
through the simplicial subsets 
\[
\mathrm{Hom_{\mathit{F^{op}(\sigma(j'_{max}))}}^{\bullet}\!(}\alpha_{*}(\tau(j_{min}^{'})),\tau(j_{max}^{'}))\subset F(\sigma(j_{max}')).
\]
Here $\alpha_{*}\!\overset{\mathrm{def}}{=}\!F^{op}(\sigma(j'_{min}\longrightarrow j'_{max})$
and 
\[
\tau(j_{min}^{'})\!:\!\triangle^{0}\longrightarrow F^{op}(\sigma(j'_{min}))
\]
\[
\tau(j_{max}^{'})\!:\!\triangle^{0}\longrightarrow F^{op}(\sigma(j'_{max}))
\]
are the morphisms corresponded to the inclusions $j_{min}^{'}\subset J'\subset J$,
$j_{max}^{'}\subset J'\subset J$. 

With this in mind we can define the morphisms 
\[
\varphi_{n}\!:\!N_{F^{op}}(C)_{n}\longrightarrow N_{F^{\rm adj}}(C^{op})_{n}.
\]
For $n=0$ this is identity map. For $n>0$ we map the data (RN1)
and (RN2) to the following data:
\begin{itemize}
	\item [{$\rm (RN1^{adj}) $}] the same functor $\sigma$ from $J^{op}$ to $C^{op}$. 
	\item [{$\rm (RN2^{adj}) $}] For every nonempty subset $J'^{op}\subseteq J^{op}$ having
	the\textit{ maximal} element $j_{min}'$ , the map $\tau(J'^{op})^{\vee}:\triangle^{J'}\longrightarrow F(\sigma(j_{min}'))$
	which define as the composition
	\[
	\xymatrix{\triangle^{J'}\ar[d]\ar[rr]^{\!\!\!\!\!\!\!\!\!\!\tau(J'^{op})^{\vee}} &  & F(\sigma(j_{min}'))\\
		\mathrm{Hom_{\mathit{F(\sigma(j'_{max}))}}^{\bullet}\!(}\tau(j_{max}^{'}),\alpha_{*}\tau(j_{min}^{'})\!)\ar[rr] &  & \mathrm{Hom_{\mathit{F(\sigma(j'_{min}))}}^{\bullet}}(\alpha^{*}\tau(j_{max}^{'}),\tau(j_{min}^{'})\!)\ar@{^{(}->}[u]
	}
	\]
\end{itemize}
Here $\alpha^{*}\!\overset{\mathrm{def}}{=}\!F^{\rm adj}(\sigma(j'_{max}\longrightarrow j'_{min}))$
(recall that we replace $J$ with $J^{op}$). Notice that the condition
(RN3) for data ($\rm RN1^{adj} $) and ($\rm RN2^{adj} $) automatically satisfied.
%\begin{rem}
%	We will not use that $\varphi$ is the equivalence of category. We
%only care about the construction of $\varphi$ itself for comparison
%with \textcolor{red}{(8.7)}. %
%	\begin{comment}
%\textcolor{blue}{can you prove that this is an equivalence of categories?
%	It is not enough to prove that this is a weak homotopy equivalence
%	of simplicial sets because we work with the Joal model structure!}
%	\end{comment}
%\end{rem}
\begin{prop}\label{PB.4}
	The map $\varphi$ is an equivalence of $\infty$-categories.
\end{prop}
\textit{Proof:} It suffice to check that for any pair $(c_1,Y_1)$ and $(c_2,Y_2)$ the map 
$$\varphi: \mathrm{Hom}^\bullet_{\int_{C}(F^{\mathrm{op}})}((c_1,Y_1),(c_2,Y_2))\ra \mathrm{Hom}^\bullet_{(\int_{C^{\mathrm{op}}}F^{\rm adj})^{\mathrm{op}}}((c_1,Y_1),(c_2,Y_2))$$ is a weak equivalence of simplicial sets. Note that we have the the commutative diagram of spaces 
$$\begin{tikzcd}
	{|\mathrm{Hom}^\bullet_{\int_{C}(F^{\mathrm{op}})}((c_1,Y_1),(c_2,Y_2))|} \arrow[rd, "\pi_1"] \arrow[rr, "\varphi"] &                           & {|\mathrm{Hom}^\bullet_{(\int_{C^{\mathrm{op}}}F^{\rm adj})^{\mathrm{op}}}((c_1,Y_1),(c_2,Y_2))|} \arrow[ld, "\pi_2"'] \\
	& {\mathrm{Hom}_C(c_1,c_1)} &                                                                                                                       
\end{tikzcd}$$
and $\varphi$ maps $\pi_1^{-1}(\alpha)$ to $\pi_2^{-1}(\alpha)$ for any $\alpha$. Moreover, $\pi_i^{-1}(\alpha)$ and $\pi_i^{-1}(\alpha')$ lie in different connection components for different $\alpha$ and $\alpha'.$ Hence it suffice to prove that $\varphi$ induces the weak equivalence 
$$\Delta^1\times_C\mathrm{Hom}^\bullet_{\int_{C}(F^{\mathrm{op}})}((c_1,Y_1),(c_2,Y_2))\to \Delta^1\times_C\mathrm{Hom}^\bullet_{(\int_{C^{\mathrm{op}}}F^{\rm adj})^{\mathrm{op}}}((c_1,Y_1),(c_2,Y_2))$$
for any $\alpha:\Delta^1\to C.$
\\On the other hand, for any 1-functor $C'\to C$ the fiber product $C'\times_C \int_CF$ 
compatible with the inverse image $C'\to C\overset{F}{\to}\mathbf{Pr}^{\rm R}.$ So we may assume that $C=\Delta^1$. In this case the statement follows from Remark \ref{RB.5}. $\square$

\end{sloppypar}

%\include{Sh1}
%\include{LM}
%\bibliography{referenses}
%\include{appendix}
%\include{referenses}

\textsc{Skolkovo Institute of Science and Technology, 3 Nobelya, Moscow,
	121205}

\textit{\negmedspace{}\negmedspace{}\negmedspace{}\negmedspace{}\negmedspace{}}\textsc{National
	Research University Higher School of Economics, Department of Mathematics,
	6 Usacheva st, Moscow, 119048 ; }

\textit{\negmedspace{}\negmedspace{}\negmedspace{}\negmedspace{}\negmedspace{}$\,\,\,\,\,\,\,\,$E-mail
	address:} \textbf{Shuklin.ium@yandex.ru}

\end{document}